\documentclass[a4paper,12pt]{article}
\usepackage{dcolumn}
\usepackage{bm}
\usepackage[dvips]{graphicx}
\usepackage{amsmath}
\usepackage{hyperref}
\usepackage{amssymb}
\usepackage{amsfonts}
\usepackage{ucs}
\usepackage[utf8x]{inputenc}
\usepackage[T1]{fontenc}
\usepackage{url}
\usepackage{psfrag}
\usepackage{dsfont}
\usepackage{color}

\def\del  {\partial}
\def\eps{\epsilon}

\def\R{\mathbb{R}}

\def\N{\mathbb{N}}

 \def\dx{{\rm d}x}
 \def\dy{{\rm d}y}

\def\dz{{\rm d}z}

\def\dv{{\rm d}v}
\def\ds{{\rm d}s}

\def\dist{\operatorname{dist}}

\def\dtau{{\rm d}\tau}

\newtheorem{proposition}{\textbf{Proposition}}

\newtheorem{corollary}{\textbf{Corollary}}
\newtheorem{remark}{\textbf{Remark}}
\newtheorem{theorem}{\textbf{Theorem}}
\newtheorem{lemma}{\textbf{Lemma}}
\newtheorem{definition}{\textbf{Definition}}

\begin{document}

\title{Short time heat diffusion in compact domains with discontinuous transmission boundary conditions}

\author{CLAUDE BARDOS\footnote{Laboratory Jacques Louis Lions, University of Paris 6, Pierre et Marie Curie, 4 place Jussieu, Paris, France,
claude.bardos@gmail.com}, DENIS GREBENKOV\footnote{Laboratoire de Physique de la Mati\`ere Condens\'ee, CNRS -- Ecole Polytechnique,
Palaiseau, France,
denis.grebenkov@polytechnique.edu},\\ ANNA ROZANOVA-PIERRAT\footnote{Laboratory Applied Mathematics and Systems, CentraleSup\'elec Paris, Grande Voie des Vignes,
Ch\^atenay-Malabry, France,
anna.rozanova-pierrat@centralesupelec.fr}}

\maketitle

\begin{abstract}
We consider a heat problem with discontinuous diffusion coefficients
and discontinuous transmission boundary conditions with a resistance
coefficient.  For all compact $(\eps,\delta)$-domains
$\Omega\subset\R^n$ with a $d$-set boundary (for instance, a
self-similar fractal), we find the first term of the small-time
asymptotic expansion of the heat content in the complement of
$\Omega$, and also the second-order term in the case of a regular
boundary.  The asymptotic expansion is different for the cases of
finite and infinite resistance of the boundary.  The derived formulas
relate the heat content to the volume of the interior Minkowski
sausage and present a mathematical justification to the de Gennes'
approach.  The accuracy of the analytical results is illustrated by
solving the heat problem on prefractal domains by a finite elements
method.
\end{abstract}



\section{Introduction}

\begin{figure}[!htb]
\begin{center}
\psfrag{Op}[l]{$\Omega_+=\Omega$ ``hot''}
  \psfrag{Om}{$\Omega_-=\R^n\setminus\overline{\Omega}$ ``cold''}
  \psfrag{dO}{$\del \Omega$}
\includegraphics[width=7cm]{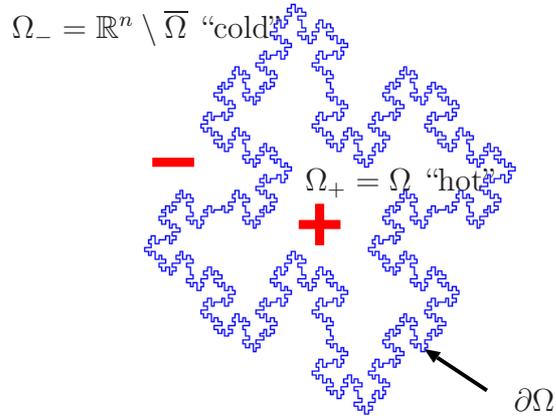}
\end{center}
\caption{Illustration of the heat content problem for a planar domain $\Omega$
with prefractal boundary $\del \Omega$ presented by the third
generation of the Minkowski fractal (of fractal dimension $3/2$).
This boundary splits the plane into two complementary regions.  At
time $t=0$, the inner region $\Omega_+=\Omega$ is ``hot'' (functions
on $\Omega_+$ are denoted with subscript $+$), while the outer region
$\Omega_-=\R^n\setminus\overline{\Omega}$ is ``cold'' (functions on
$\Omega_-$ are denoted with subscript $-$).}\label{FigOmega}
\end{figure}

We consider a compact domain $\Omega \subset \R^n$ with boundary
$\del\Omega$ that splits $\R^n$ into ``hot'' and ``cold'' media,
$\Omega_+ = \Omega$ and $\Omega_- = \R^n \setminus\overline{\Omega}$,
characterized by (distinct) heat diffusion coefficients $D_+$ and
$D_-$ (Fig.~\ref{FigOmega}).  On the boundary $\del\Omega$ is also
defined a function $0\le\lambda(x)\le \infty$ which describes the
resistivity to heat exchange through the boundary.

We are interested in propagation of the heat content associated with
the following problem:
\begin{eqnarray}
  && \del_t u_\pm- D_{\pm} \Delta u_\pm =0 \quad x\in \Omega_\pm, \; t>0,\label{prb1}\\
 && u_+|_{t=0}= 1, \quad u_-|_{t=0}=0,\label{prb1i}\\
 && D_- \frac{\del u_-}{\del n}|_{\del \Omega}=\lambda(x)(u_--u_+)|_{\del \Omega},\label{prb1c} \\
 && D_+ \frac{\del u_+}{\del n}|_{\del \Omega}=D_- \frac{\del u_-}{\del n}|_{\del \Omega},\label{prb1end}
\end{eqnarray}
where $\partial/\partial n$ is the normal derivative directed outside
the domain $\Omega$.

A rigorous analysis of the problem~(\ref{prb1})--(\ref{prb1end}) for
irregular boundaries requires its variational formulation in
appropriate functional spaces (see Section~\ref{secWP}).  The
variational problem is shown to have a unique weak solution with the
desired trace properties on the boundary $\del \Omega$ (see
Section~\ref{secWP}).  The variational problem is equivalent to the
problem~(\ref{prb1})--(\ref{prb1end}) for a piecewise Lipschitz $\del
\Omega$ according to the classical trace theorem.  In turn, extensions
of the trace theorem have to be used for fractal boundaries or, more
precisely, $d$-sets (see Subsection~\ref{SubsExFrac}).

Once a unique solution $u_{\pm}$ of the
problem~(\ref{prb1})--(\ref{prb1end}) is established, we study the
asymptotic expansion of the heat content as $t\to 0$
\begin{equation}\label{Nt}
 N(t)=\int_{\R^n\setminus \Omega} u_-(x,t)\dx=\operatorname{Vol}(\Omega)-\int_\Omega u_+(x,t)\dx.
\end{equation}

Eqs.~(\ref{prb1})--(\ref{prb1end}) describe heat exchange between two
media prepared initially at different temperatures and separated by a
partially isolating boundary\cite{Carslaw,Crank}.  In fact, $u(x,t)$
can describe how the distribution of (normalized) temperature evolves
with time.  The transmission boundary conditions (\ref{prb1c}),
(\ref{prb1end}) impose the continuity of the temperature flux across
the boundary, and relate this flux to the temperature drop at the
boundary due to thermal isolation.  The growth rate of the heat
content with time characterizes the efficiency of thermal isolation.
Understanding this problem is relevant to improve heat exchangers,
e.g., cooling of metallic radiators or thermal isolation of pipes and
buildings.  Depending on application, cooling rate has to be either
enhanced (e.g., in the case of microprocessors or nuclear reactors),
or slowed down (e.g., in the case of pipes and buildings).  For these
purposes, one can either modify the thermal isolation (i.e., the
resistivity $\lambda$), or the shape of the exchange boundary.  It is
therefore crucial to understand how the shape of the boundary
influences heat exchange.  In particular, would an irregular (e.g.,
fractal) boundary with a very large exchange area significantly speed
up cooling?

Similar equations can describe molecular diffusion between two media
across semi-permeable membranes~\cite{Tanner78,Powles92}.  In that
case, $u(x,t)$ represents the (normalized) concentration of molecules,
while Eqs.~(\ref{prb1})--(\ref{prb1end}) can model the leakage of
molecules from a cell ($\Omega_+$) to the extracellular space
($\Omega_-$) or, more generally, the diffusive exchange between two
compartments (e.g., oxygen or carbon dioxide exchange between air and
blood across the alveolar membrane in the lungs).  The resistance
$\lambda$ is related to the cellular membrane permeability.  As for
heat exchange, one may need to enhance or to slow down the molecular
leakage, and the shape of the boundary may play an important role.

The discontinuity of the initial condition, of the diffusion
coefficient, and of the solution $u(x,t)$ across the boundary between
two domains constitutes one of the mathematical difficulties to be
treated.  From a physical point of view, such discontinuities might
appear unrealistic.  For instance, the diffusive flux at the boundary
at time $t=0$ is infinite.  For any physical setting of heat or
molecular diffusion, there would be an intermediate layer between two
media in which the material properties would change rapidly but
continuously.  When the thickness of this intermediate layer is much
smaller than the size of the domain, the physical problem with
continuously varying parameters can be approximated by the heat
problem~(\ref{prb1})--(\ref{prb1end}).  Such an approximation is
applicable starting from a small cut-off time while understanding the
heat exchange at smaller time scales would need either restituting an
intermediate layer, or introducing nonlinear terms into the heat
equation.  Throughout this paper, we focus on the mathematical
problem~(\ref{prb1})--(\ref{prb1end}).

The physical properties of the two media $\Omega_+$ and $\Omega_-$ are
supposed to be different: $D_+ \ne D_-$.  This implies the
discontinuity of the metric on $\del \Omega$.  The case of continuous
metric ($g_-|_{\del \Omega}=g_+|_{\del \Omega}$) on smooth compact
$n$-dimensional Riemannian manifolds with a smooth boundary $\del
\Omega$ was considered in Ref.~\cite{GK}.  The case of continuous
transmission boundary conditions for the expansion of the heat kernel
on the diagonal was treated in~Ref.~\cite{PirN} (see also
Ref.~\cite{Vassiliev} for a survey of results on asymptotic
expansion of the heat kernel for different boundary conditions).  The
heat content asymptotic expansion with Dirichlet boundary condition
was found
\begin{itemize}
\item
up to the third-order term for a compact connected domain
$\Omega\subset \R^n$ with a regular boundary $\del \Omega\in C^3$
(Refs.~\cite{BGD,BGD2});

\item
up to an exponentially small error for a compact connected domain
$\Omega\subset \R^2$ with a polygonal $\del \Omega$
(Ref.~\cite{BS}) and for $\Omega\subset \R^2$ with $\del \Omega$
given by the triadic Von Koch snowflake (Ref.~\cite{Fleckinger});

\item
up to the second-order term for the general case of self-similar
fractal compact connected domains in $\R^n$ (Ref.~\cite{Levitin}).
\end{itemize}
 
In general, the boundary between two media can have some resistance to
heat exchange, described by the function $\lambda(x)\geq 0$ ($x\in
\del \Omega$) that may account for partial thermal isolation.  We
outline three cases of boundary conditions according to $\lambda$:
\begin{enumerate}
\item 
If $0<\lambda(x) <\infty$ for all $x\in \del \Omega$, $u$ is discontinuous
on $\del \Omega$ and we have: 
\begin{equation*}
\left(\lambda(x)u_--D_- \frac{\del u_-}{\del n}\right)|_{\del
\Omega}=\lambda(x)u_+|_{\del \Omega}, \quad D_+
\frac{\del u_+}{\del n}|_{\del \Omega}=D_- \frac{\del u_-}{\del
n}|_{\del \Omega}.
\end{equation*}

\item 
If $\lambda=+\infty$ for all $x\in \del \Omega$, $u$ is continuous on
$\del \Omega$ due to the transmission condition and in this case 
\begin{equation*}
u_+|_{\del \Omega}=u_-|_{\del \Omega},\quad D_+ \frac{\del u_+}{\del n}|_{\del
\Omega}=D_- \frac{\del u_-}{\del n}|_{\del \Omega}.
\end{equation*}

\item 
If $\lambda=0$ for all $x\in \del \Omega$, we have the Neumann
boundary condition 
\begin{equation*}
\frac{\del u_-}{\del n}|_{\del \Omega}=\frac{\del u_+}{\del n}|_{\del \Omega}=0
\end{equation*}
that models the complete thermal isolation of $\del \Omega$ and
implies the trivial solution given by $u_-(x,t) = 0$ and $u_+(x,t) =
1$ for all time $t\ge 0$.
\end{enumerate}

The main goal of the article is to develop the preliminary
study\cite{Anna} and especially to formalize the seminal approach by
de Gennes\cite{deGennes82}.  In the case $\lambda=+\infty$, de Gennes
argued that as $t\to +0$, $N(t)$ is proportional to the volume
$\mu(\del \Omega,\sqrt{D_+ t})$ of the interior Minkowski sausage of
$\del \Omega$ of the width equal to the diffusion length $\sqrt{D_+
t}$:
\begin{equation*}
\mu(\del \Omega,\ell) = \operatorname{Vol}\bigl(\{ x\in \Omega|\operatorname{dist}(x,\del \Omega)<\ell\}\bigr)
\end{equation*}
(see also Ref.~\cite{Levitin}).  In particular,
\begin{itemize}
\item 
for a regular boundary $\del \Omega$, $N(t)$ is proportional to
$\operatorname{Vol}(\del \Omega) \sqrt{D_+t}$;

\item 
for a fractal boundary $\del \Omega$ of the Hausdorff dimension $d$,
$N(t)$ is proportional to $(D_+t)^{\frac{n-d}{2}}$.
\end{itemize}
The de Gennes scaling argument was further investigated in
Ref.~\cite{Anna}, both experimentally and numerically.  It was
shown that irregularly shaped passive coolers rapidly dissipate at
short times, but their efficiency decreases with time.  The de Gennes
scaling argument was shown to be only a large scale approximation,
which is not sufficient to describe adequately the temperature
distribution close to the irregular frontier.

In the present paper, we provide a mathematical foundation and further
understanding for the de Gennes approach.  We obtain three results
valid for all compact $(\eps,\delta)$-domains $\Omega$ in $\R^n$ with
connected boundary $\del \Omega$, presented by a closed $d$-set (see
Section~\ref{SubsExFrac} for the definitions of
$(\eps,\delta)$-domains and $d$-sets): the well-posedness of the
problem~(\ref{prb1})--(\ref{prb1end}), the continuity of the solution
on $\lambda$ (see Section~\ref{secWP}), and the asymptotic expansion
of the heat content~(\ref{Nt}).  In particular, these results hold for
domains with a self-similar fractal boundary.

We show in Theorem~\ref{FinLem} that the heat content $N(t)$ is
approximated by the volume of the interior Minkowski sausage of $\del
\Omega$ of the radius $\sqrt{4 D_+ t}$:
\begin{eqnarray}
N(t)=\tau_\lambda\left[ C_\lambda(\del \Omega)\mu\left(\del \Omega,\sqrt{4D_+t}\right)
+O\left(\mu^2\left(\del \Omega,\sqrt{4D_+t}\right)\right)\right],\label{EqEtoile}
\end{eqnarray}
where $\tau_\lambda$ is equal to $1$ if $\lambda=\infty$ and
$\sqrt{t}$ if $\lambda>0$ is finite.  Here $C_\lambda(\del \Omega)$ is
a constant depending only on the shape of $\del \Omega$ and finiteness
of $\lambda$ (see Theorem~\ref{FinLem} for the exact formulas).
Formula~(\ref{EqEtoile}) is the first approximation of
Eqs~(\ref{NilocFint2C}),~(\ref{NilocFint2InfC}) given in
Theorem~\ref{FinLem}, which allows to find $N(t)$ up to terms of the
order $\tau_\lambda O\left(\sqrt{t}~\mu\left(\del
\Omega,\sqrt{4D_+t}\right)\right)$.

Moreover, the asymptotic relation \eqref{EqEtoile} remains valid even
for mixed boundary conditions for three disjoint boundary parts, i.e.
when $\lambda = \infty$ on one part of the boundary, $\lambda = 0$ on
another part, and $0 < \lambda < \infty$ on the remaining boundary
(see Theorem~\ref{ThMixed}).  However, changes of the type of the
boundary condition should be continuous (see Theorem~\ref{Thulamc})
such that $u$ remains a continuous function of $\lambda$.  In this
more general case, the coefficient $C_\lambda(\del \Omega)$ in
Eq.~(\ref{EqEtoile}) is given either by Eq. (\ref{NFint2Ca}) for $0
<\lambda<\infty$, or by Eq. (\ref{NFint2InfCa}) for $\lambda =
\infty$, or is equal to $0$ for $\lambda = 0$ (the boundary with
$\lambda = 0$ does not contribute to the short-time asymptotics of the
heat content).  Finding the asymptotics for mixed boundary conditions
with a discontinuous jump from a finite $\lambda$ to $\lambda=\infty$
is still an open problem.

As expected, the resistivity of the boundary to heat transfer makes
heat diffusion {\it slower} due to the presence of the coefficient
$\tau_\lambda=\sqrt{t}$.

For a fractal boundary we replace $\mu\left(\del
\Omega,\sqrt{4D_+t}\right)$ by the volume of the interior Minkowski
sausage which scales as $(4D_+ t)^{(n-d)/2}$, where $d$ is the fractal
dimension\cite{Levitin}.  In the fractal case the integral over $\del
\Omega$ should be understood by using the Hausdorff measure
(see~Ref.~\cite{Kigami,Giona,Grisvard}).
\begin{figure}[!htb] 
\begin{center}
\includegraphics[width=68mm]{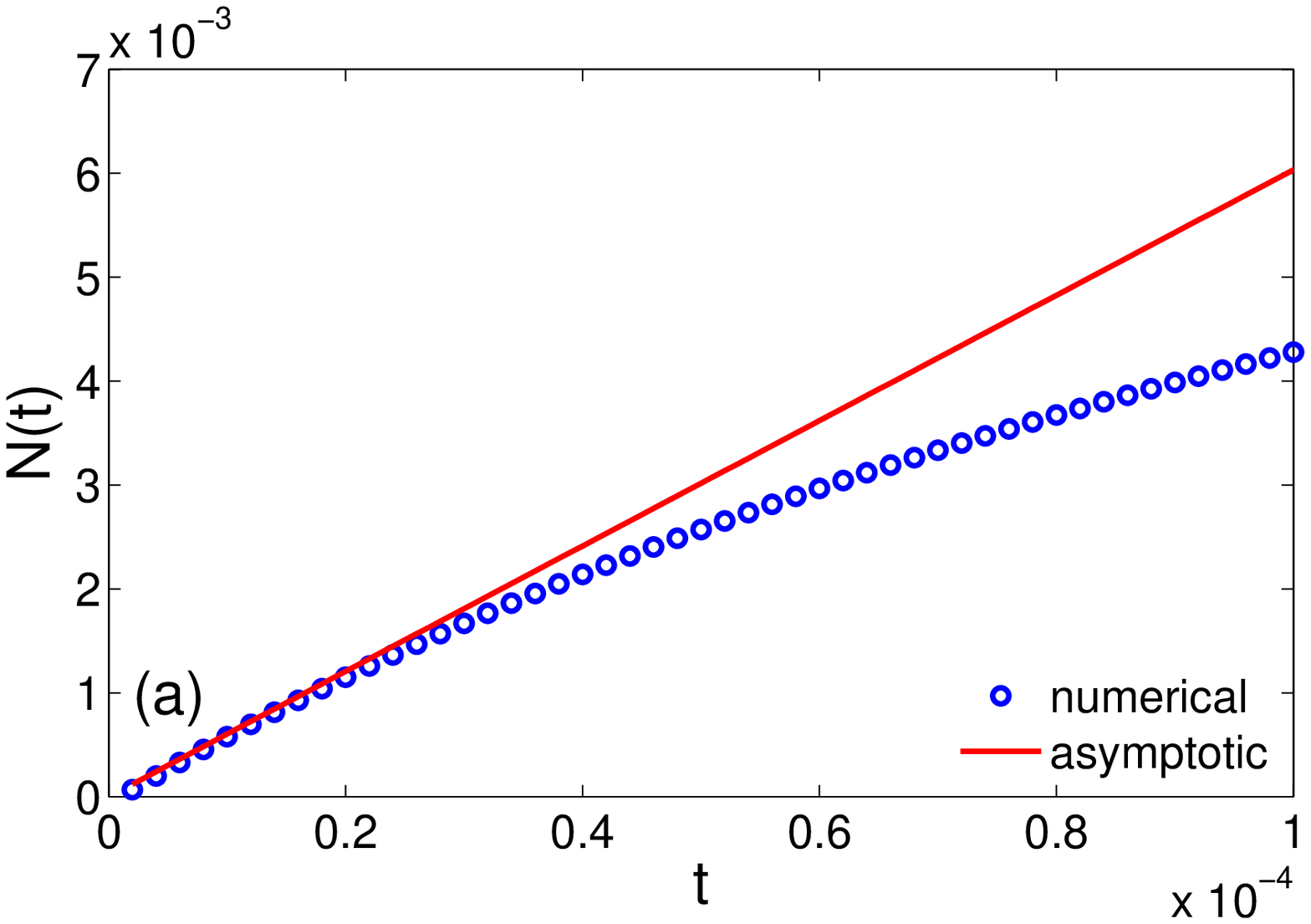}   
\includegraphics[width=68mm]{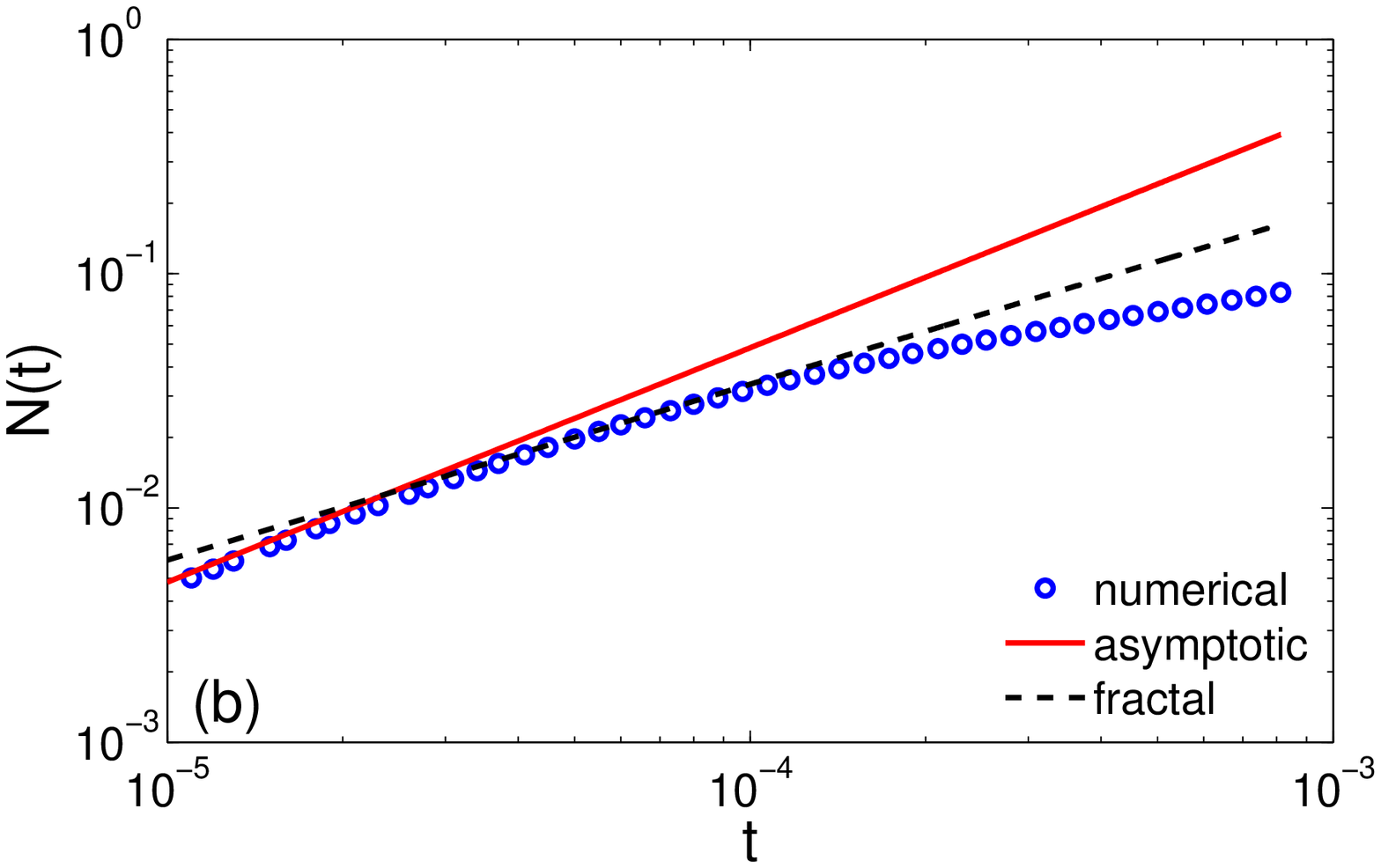}   
\end{center}
\caption{
Comparison between the asymptotic formula~(\ref{EqEtoile}) (solid
line) and a FreeFem++ numerical solution of the
problem~(\ref{prb1})--(\ref{prb1end}) (circles) for two domains: (a)
the unit square ($\operatorname{Vol}(\del \Omega)=4$) and (b) the
third generation of the Minkowski fractal ($\operatorname{Vol}(\del
\Omega)=2^3\cdot4$), with $D_+=1/100$, $D_-=1$, and $\lambda=17$.
Since the Hausdorff dimension of the boundaries of these domains is
$1$ (even for the {\it prefractal} case), Eq. (\ref{EqEtoile}) for a
constant $\lambda$ is reduced, according to Theorem \ref{FinLem}, to
$N(t)=2\sqrt{t} C_0 \lambda \mu(\del\Omega,\sqrt{4D_+t}) +
O(t^{\frac32})$ with $\mu(\del\Omega,\sqrt{4D_+t})\simeq\sqrt{4D_+t}
\operatorname{Vol}(\del\Omega)$ and $C_0$ given by
Eq. (\ref{eq:C0}).  For plot (b), dashed line shows the fractal
asymptotic (that would be exact for the infinite generation of the
fractal) with de Gennes approximation of $\mu\left(\del
\Omega,\sqrt{4D_+t}\right)$ in Eq.~(\ref{EqEtoile}) by
$(4D_+t)^\frac{1}{4}$.  This approximation is valid for intermediate
times.}\label{FigLam}
\end{figure}

\begin{figure}[!htb]
\begin{center}
\includegraphics[width=68mm]{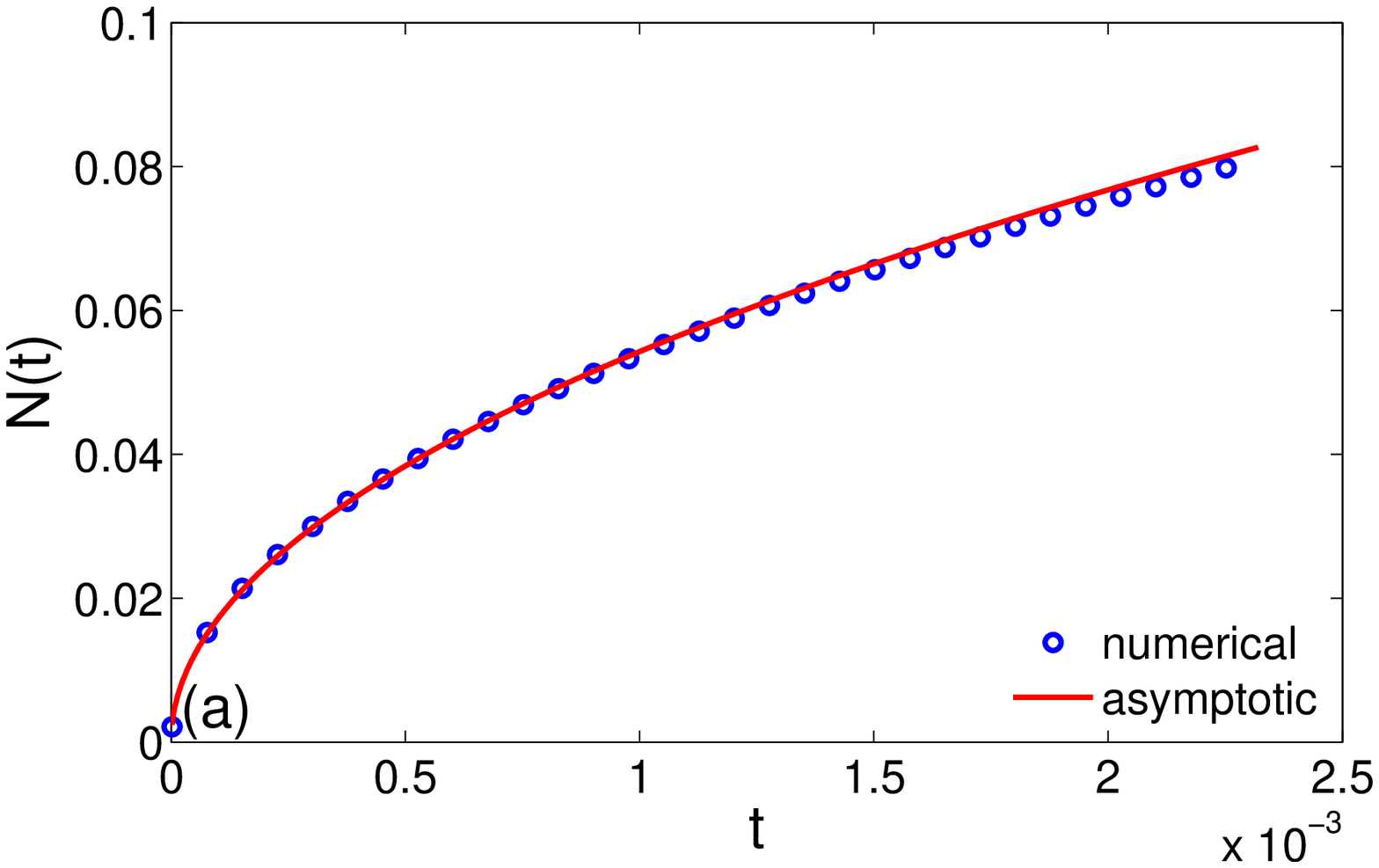}
\includegraphics[width=68mm]{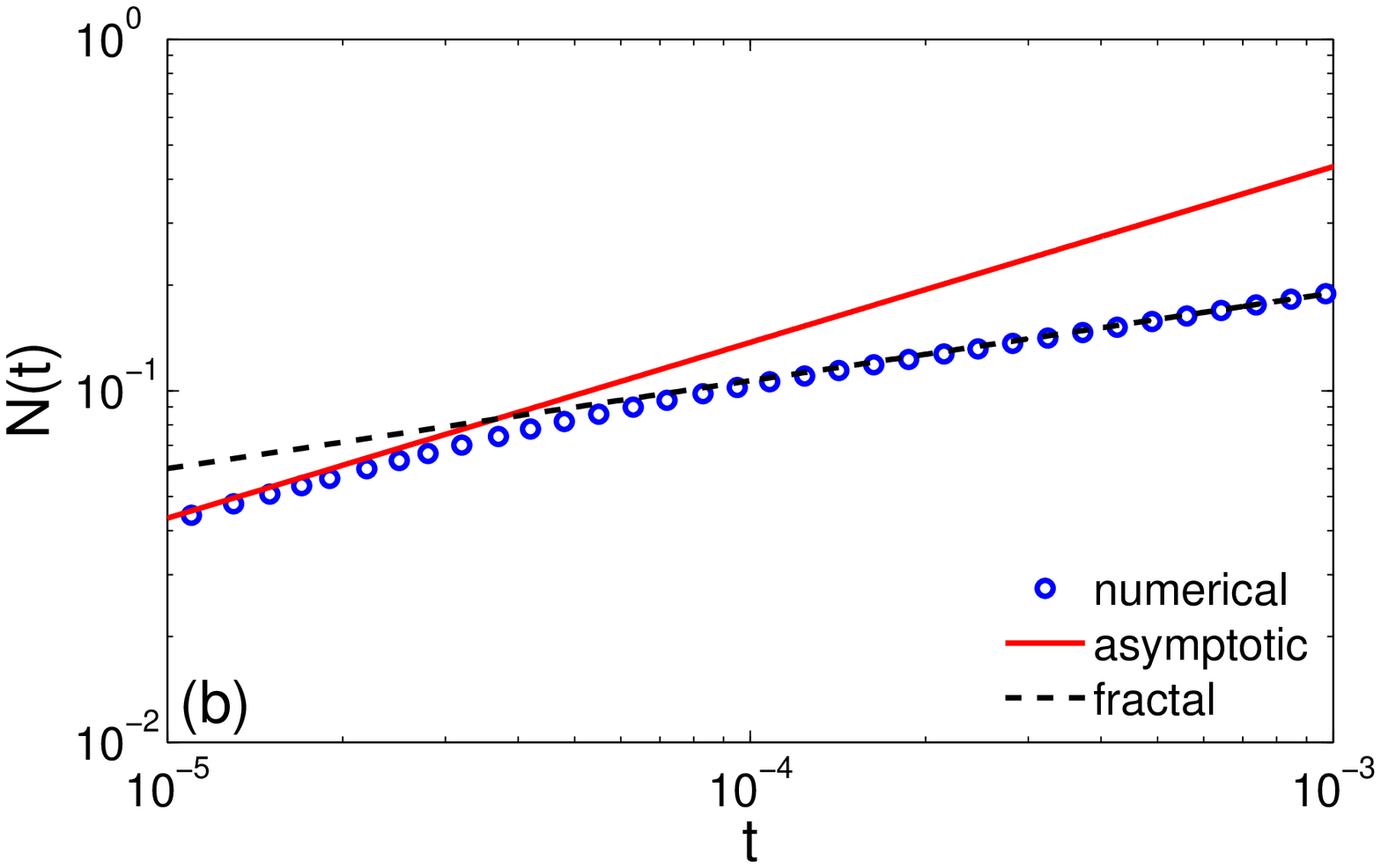}   
\end{center}
\caption{
Comparison between the asymptotic formula~(\ref{EqEtoile}) (solid
line) and a FreeFem++ numerical solution of the
problem~(\ref{prb1})--(\ref{prb1end}) (circles) for two domains: (a)
the unit square ($\operatorname{Vol}(\del \Omega)=4$), and (b) the
third generation of the Minkowski fractal ($\operatorname{Vol}(\del
\Omega)=2^3\cdot4$), with $D_+=0.4$, $D_-=1$, and $\lambda = \infty$.
Since the {\it prefractal} boundary $\del\Omega$ has the Hausdorff
dimension $1$, Eq. (\ref{EqEtoile}) is reduced to Eq. (\ref{NtregLI}),
i.e., $N(t) \propto \sqrt{t}$.  In turn, dashed line shows the fractal
asymptotic (that would be exact for the infinite generation of the
fractal) with de Gennes approximation of $\mu\left(\del
\Omega,\sqrt{4D_+t}\right)$ in Eq.~(\ref{EqEtoile}) by
$2.5(4D_+t)^\frac{1}{4}$.  This approximation is valid for
intermediate times. }\label{FigLamInfF3}
\end{figure}

The comparison between the asymptotic formula~(\ref{EqEtoile}) and a
numerical solution of the problem~(\ref{prb1})--(\ref{prb1end}) for
the unit square and a prefractal domain is shown in Fig.~\ref{FigLam}
for a finite $\lambda$ and in Fig. \ref{FigLamInfF3} for
$\lambda=+\infty$.  The numerical solution was obtained in FreeFem++
by a finite elements method with the implicit $\theta$-schema, also
known as Crank-Nicolson schema, for the time discretization with
$\theta=\frac{1}{2}$ and $\Delta t=10^{-6}$.  The domain $\Omega$ was
centered in a ball $B$ of diameter (at least) twice bigger than the
diameter of $\Omega$.  The Neumann boundary condition was imposed on
the boundary of the ball.  According to the principle ``not feeling
the boundary''\cite{Fleckinger} (see also Section~\ref{secLoc}), the
heat content propagation in $\R^2$ with a prescribed boundary $\del
\Omega$ can be very accurately approximated at small times by
the heat content propagation computed in $B$.  The accuracy of this
approximation can also be checked by changing the diameter of the
ball.  In the case of the square domain $\Omega$, the ball was
replaced by a square with four times bigger edge.  Each pre-fractal
edge was discretized with 27 space points while 57 points were used in
the external boundary of the ball.  The mesh size was varied to check
the accuracy of the presented numerical solutions.  For the case of
the discontinuous solution on the boundary (when $0<\lambda<\infty$)
we apply the domain decomposition method and match the boundary values
of the respective solutions on $\del \Omega$ by a Picard fixed point
method.  We consider therefore the numerical solution of heat
propagation for small times as a reference, to which asymptotic
formulas are compared with.  In particular, deviations between the
numerical solution and the asymptotic formulas observed at longer
times illustrate the range of validity of the short-time expansion.

For the regular case $\del \Omega\in C^3$, we obtain the heat content
approximation up to the third-order term.  The formulas are given in
Theorem~\ref{ThFinalr}.  For the case $\lambda<\infty$, the
coefficient in front of the second-order term ($t^\frac{3}{2}$) in the
asymptotic expansion depends on the mean curvature.  In turn, for
$\lambda=\infty$, the second-order term (here, $t$) in the asymptotic
expansion vanishes:
\begin{equation}
N(t)=2\frac{1-e^{-4}}{\sqrt{\pi}} ~
\frac{\sqrt{D_-D_+}}{\sqrt{D_+}+\sqrt{D_-}}\operatorname{Vol}(\del
\Omega) \sqrt{t} + O(t^{\frac{3}{2}}).
\end{equation}

The rest of the paper is organized as follows.  In
Section~\ref{secWP}, we describe the class of irregular boundaries and
prove the well-posedness of the model relying on the variational
formulation of the problem.  The boundary conditions are treated in
the weak sense by generalizing the trace operator and the Green
formula to fractals using fractal Besov spaces, $B^{2,2}_\beta(\del
\Omega)$ and $B^{2,2}_{-\beta}(\del \Omega)$
($\beta=1-\frac{n-d}{2}>0$ for a $d$-dimensional $\del \Omega$)
defined in~\ref{ApBS}.  In Section~\ref{secWP} we also establish the
continuity of $u$ as a function of $\lambda$.  In Section~\ref{secLoc}
we prove that the problem to find $N(t)$ can be replaced by a heat
problem localized in $O(\sqrt{t})$-interior Minkowski sausage of the
boundary by a variant of the principle ``not feeling the
boundary''\cite{Fleckinger} in the general case in $\R^n$.  This
allows, due to the continuity of $u$ on $\lambda$, to establish
Theorem~\ref{ThMixed} for a mixed boundary condition including zero,
finite, or infinite values of $\lambda$.  Considering a regular $\del
\Omega$ (at least in $C^3$) and using the localization properties from
Section~\ref{secLoc}, we rewrite in Section~\ref{secLC} the formula
for $N(t)$ in the terms of the local coordinates.
Section~\ref{secAprox} gives the approximation of the heat problem
solution through the solution of one-dimensional constant coefficient
problem.  The heat content is calculated in terms of the volume of the
interior Minkowski sausage of the boundary in Section~\ref{subs2}.
Firstly, to illustrate the technique of the proof on a simple case, we
give the proof for the case of continuous diffusion coefficients
$D_+=D_-$, just with discontinuity of the initial condition.  In this
case, all formulas given in Section~\ref{subs2} are valid for all
types of the boundary introduced in Subsection~\ref{SubsExFrac}.  The
calculation relies on the Green function of the problem with constant
coefficients for $\Omega$ being a half-space (see~\ref{SecGreen}).  We
also obtain the Green function used in Section~\ref{SecR} for the
proof of the asymptotic heat expansion up to the third-order term for
a regular $\del \Omega\in C^3$.

\section{Well-posedness of the model}\label{secWP}

Let $\Omega$ be an open connected bounded subset of $\R^n$ such that
$\del \Omega$ is closed with $\operatorname{Vol}(\Omega)<\infty$.  We
denote by $\Omega_+=\Omega$ and
$\Omega_-=\R^n\setminus\overline{\Omega}$ (Fig.~\ref{FigOmega}).

We are looking for the solution of the
problem~(\ref{prb1})--(\ref{prb1end}), where $D_+\ne D_-$, $D_+>0$ and
$D_->0$, $\lambda(x)\ge0$ for all $x\in \del \Omega$.  The boundary
$\del \Omega$ is divided into two disjoint parts:
$\Gamma_\infty=\{x\in
\del \Omega|\; \lambda(x)=+\infty\}$ and $\del
\Omega\setminus\Gamma_\infty=\{x\in \del \Omega|\;0\le
\lambda(x)<+\infty\}$.  Each of the parts can be the empty set.  We
thus assume that $\lambda\in L^\infty(\del
\Omega\setminus\Gamma_\infty)$.

\subsection{Regular boundary: at least piecewise Lipschitz}\label{SubsecWPLip}

Firstly, we consider the case when $\del \Omega$ is regular (at least
piecewise Lipschitz) and $\Gamma_\infty$ is the empty set.

To prove the existence, the uniqueness, and the stability of a
solution of the problem~(\ref{prb1})--(\ref{prb1end}), we proceed with
its variational formulation.

We introduce the space $H=L^2(\R^n)$ and the space 
\begin{equation*}
V=\{f\in H|\;  f_+=f|_{\Omega_+}\in H^1(\Omega_+),\;\hbox{and } f_-=f|_{\Omega_-}\in
H^1(\Omega_-)\}
\end{equation*}
of functions $f=f_+\mathds{1}_{\Omega_+}+f_-\mathds{1}_{\Omega_-}$
defined on $\Omega_+\cup\Omega_-$ such that their restrictions
$f_+=f|_{\Omega_+}$ and $f_-=f|_{\Omega_-}$ belong to $H^1$.  We equip
$V$ with the norm:
\begin{equation*}
 \|u\|_V^2=D_+\int_{\Omega_+} |\nabla u_+|^2\dx+D_-\int_{\Omega_-} |\nabla u_-|^2 \dx+\int_{\Omega_+\cup\Omega_-} |u|^2 \dx.
\end{equation*}
We notice that $V$ is a Hilbert space, $V\subset L^2(\Omega)$, and $V$
is dense in $L^2(\Omega)$.  In addition, $V\subset L^2(\R^n) \subset
V'$, where $V'$ is the dual space to $V$.  Finally, since $\del
\Omega$ is regular, the inclusion $V\subset L^2(\R^n)$ is compact.

Applying the usual trace theorem under the assumptions that $\Omega$ is bounded and $\del
\Omega$ is at least piecewise Lipschitz, the bilinear form
\begin{multline}\label{VarForm1}
 a(u,v)=D_+\int_{\Omega_+} \nabla u_+ \nabla v_++D_-\int_{\Omega_-} \nabla u_- \nabla v_-+
 \int_{\del \Omega_+}\lambda(x)(u_+-u_-)(v_+-v_-) d \sigma
\end{multline}
is continuous,
\begin{equation*}
|a(u,v)|\le C(\|\lambda\|_{L^\infty(\del \Omega)},D_+,D_-,\Omega_+) \|u\|_V\|v\|_V \quad (\hbox{for a constant } C>0),
\end{equation*}
and coercive on $V\times V$, i.e.,
\begin{eqnarray*}
 &&a(u,u)=D_+\int_{\Omega_+} |\nabla u_+|^2\dx+D_-\int_{\Omega_-} |\nabla u_-|^2 \dx+\int_{\del \Omega_+}\lambda(x)|u_+-u_-|^2 d \sigma\\
 &&\ge \|u\|^2_V-\|u\|^2_{L^2(\R^n)}>0.
\end{eqnarray*}

Thus we conclude\cite{LionsM} that the bilinear form $a(u,v)$ defines
an operator $A:V\to V'$ by $a(u,v)= \langle Au,v\rangle$.  Moreover,
$-A|_{L^2(\R^n)}$ with $D(A)=\{u\in V|\; Au\in L^2(\R^n)\}$ generates
an analytical semigroup.

\begin{remark}
When $\Gamma_\infty$ is not empty, the variational
form~(\ref{VarForm1}) is well adaptable to the case where $u$ is
continuous across the part $\Gamma_\infty\subset \del
\Omega$ of the interface.  By convention we put on this part
$\lambda(x)=\infty$ which implies $u_+=u_-$ on $\Gamma_\infty$ (see
also Theorem~\ref{Thulamc}).

For $\Gamma_\infty \ne \emptyset$, we introduce $V$ as the space of
functions $u\in L^2(\R^n)$ such that
\begin{equation*}
u_+=u|_{\Omega_+}\in H^1(\Omega_+),\quad u_-=u_-|_{\Omega_-}\in H^1(\Omega_-),\quad u_+|_{\Gamma_\infty}=u_-|_{\Gamma_\infty},
\end{equation*}
and, therefore, we consider the bilinear continuous and coercive form
on $V\times V$
\begin{multline}
 a(u,v)=D_+\int_{\Omega_+} \nabla u_+ \nabla v_++D_-\int_{\Omega_-} \nabla u_- \nabla v_-+
 \int_{\del \Omega\setminus\Gamma_\infty }\lambda(x)(u_+-u_-)(v_+-v_-) d \sigma.
\end{multline}
In particular, for $\Gamma_\infty=\del \Omega$, we get
$V=H^1(\Omega_+\cup\Omega_-)$ and
\begin{equation*}
a(u,v)=D_+\int_{\Omega_+} \nabla u_+ \nabla v_++D_-\int_{\Omega_-} \nabla u_- \nabla v_-.
\end{equation*}
\end{remark}

\subsection{Extension to $d$-sets (fractal case for $d>n-1$)}\label{SubsExFrac}

Let us define a class of fractal domains to be considered.  We will
see that the existence and uniqueness results of a weak solution of
the problem~(\ref{prb1})--(\ref{prb1end}) hold for a class of bounded
$(\eps,\delta)$-domains\cite{Jones,JW,Wallin} $\Omega_+$ such that
$\del \Omega$ is a $d$-set\cite{JW}:
\begin{definition}\label{Defdset}\textbf{($d$-set\cite{JW,Wallin,JWBesov})} 
Let $\Gamma$ be a closed subset of $\R^n$ and $0<d\le n$.  A positive
Borel measure $m_d$ with support $\Gamma$ is called a $d$-measure of
$\Gamma$ if, for some positive constants $c_1$, $c_2>0$,
\begin{equation*}
c_1r^d\le m_d(\Gamma\cap U_r(x))\le c_2 r^d, \quad \hbox{ for  } ~ \forall~x\in \Gamma,\; 0<r\le 1,
\end{equation*}
where $U_r(x)\subset \R^n$ denotes the Euclidean ball centered at $x$
and of radius~$r$.

The set $\Gamma$ is a $d$-set if there exists a $d$-measure on $\Gamma$.
\end{definition}
As it is known from Ref.~[\cite{MEas}, p.30], any two $d$-measures
on $\Gamma$ are equivalent.
\begin{definition}\textbf{($(\eps,\delta)$-domain\cite{Jones,JW,Wallin})}
An open connected subset $\Omega$ of $\R^n$ is an
$(\eps,\delta)$-domain, $\eps > 0$, $0 < \delta \leq \infty$, if
whenever $x, y \in \Omega$ and $|x - y| < \delta$, there is a
rectifiable arc $\gamma\subset \Omega$ with length $\ell(\gamma)$
joining $x$ to $y$ and satisfying
\begin{enumerate}
 \item $\ell(\gamma)\le \frac{|x-y|}{\eps}$ and
 \item $d(z,\del \Omega)\ge \eps |x-z|\frac{|y-z|}{|x-y|}$ for $z\in \gamma$. 
\end{enumerate}
\end{definition}
In particular, a Lipschitz domain $\Omega$ is an
$(\eps,\delta)$-domain and also a $n$-set\cite{Wallin} (i.e., a
$d$-set with $d=n$).  Self-similar fractals (e.g., von Koch's
snowflake domain) are examples of $(\eps,\delta)$-domains with the
$d$-set boundary\cite{Capitan,Wallin}, $d>n-1$.

In order to describe irregular boundaries of fractal dimension
$d>n-1$, we define sets preserving Markov's inequality
(Ref.~\cite{JW} Ch. II):
\begin{definition}
A closed subset $V$ in $\R^n$ preserves Markov's inequality if for
every fixed positive integer $k$, there exists a constant $c=c(V,n,k)
> 0$, such that 
\begin{equation*}
\max_{V\cap U_r(x)} |\nabla P | \le \frac{c}{r}\max_{V\cap U_r(x)}|P|
\end{equation*}
for all polynomials $P \in \mathcal{P}_k$ and all closed balls
$U_r(x)$, $x \in V$ and $0 < r \le 1$.
\end{definition}
Examples of sets that preserves Markov's inequality are $d$-sets in
$\R^n$, where $d > n -1$, and self-similar sets that are not a subset
of any $(n-1)$-dimensional subspace of $\R^n$ (see
Refs. \cite{Wallin,Bodin}).

To extend the variational formulation introduced in
Subsection~\ref{SubsecWPLip} to fractal boundaries of the type of
$d$-sets, we use the existence of the $d$-dimensional Hausdorff
measure $m_d$ on $\del \Omega$ (the $d$-measure from
Definition~\ref{Defdset}) and the theorem which generalizes the usual
trace theorem and the Green formula.

For example, for $d=n-1$ and a Lipschitz $\del \Omega$, we
know\cite{LionsM,Marsc} that the trace operator is linear continuous
and surjective from $H^1(\Omega)$ onto $H^\frac{1}{2}(\del \Omega)$,
and the formula
\begin{equation}
  \int_\Omega v\Delta u\dx=\langle \frac{\del u}{\del \nu}, \mathrm{Tr}v\rangle _{((H^\frac{1}{2}(\del \Omega))', H^\frac{1}{2}(\del \Omega))}
-\int_\Omega \nabla v \nabla u \dx,\nonumber
\end{equation}
holds whatever $u\in H^1(\Omega)$ such that $\Delta u\in L^2(\Omega)$
and $v\in H^1(\Omega)$.  

To generalize the trace operator and the Green formula to fractal
boundaries, one introduces the Besov space $B_\beta^{2,2}(\del
\Omega)$ with $\beta=1-\frac{n-d}{2}>0$ (see~\ref{ApBS}).  Note that
for $d=n-1$, one has $\beta=\frac{1}{2}$ and
\begin{equation*}
B_\frac{1}{2}^{2,2}(\del \Omega)=H^\frac{1}{2}(\del \Omega) ,
\end{equation*}
i.e., one recovers the above relations.  In general, 
\begin{enumerate}
\item For an arbitrary open set $\Omega$ of $\R^n$, the trace operator
$\mathrm{Tr}$ is defined\cite{JW,Bodin,Lancia} for $u\in
L^1_{loc}(\Omega)$ by
\begin{equation}\label{TraceF}
 \mathrm{Tr} u(x)=\lim_{r\to 0} \frac{1}{m(\Omega\cap U_r(x))}\int_{\Omega\cap U_r(x)}u(y)dy,
\end{equation}
where $m$ denotes the Lebesgue measure.  The trace operator
$\mathrm{Tr}$ is considered for all $x\in\overline{\Omega}$ for which
the limit exists.
\item If $\Omega$ is a bounded $(\eps,\delta)$-domain in $\R^n$
such that its boundary $\del \Omega$ is a closed $d$-set preserving
Markov's inequality, then\cite{JW,Wallin}
 \begin{enumerate}
  \item the trace operator
$\mathrm{Tr}: H^1(\Omega) \to B^{2,2}_{\beta}(\del \Omega)$  is linear continuous and surjective;
\item the Green formula holds (see also Refs.~\cite{Lancia,Capitan2} for
the von Koch case in $\R^2$):
\begin{equation}\label{FracGreen}
 \int_\Omega v\Delta u\dx=\langle \frac{\del u}{\del \nu}, 
\mathrm{Tr}v\rangle _{((B^{2,2}_{\beta}(\del \Omega))', B^{2,2}_{\beta}(\del \Omega))}-\int_\Omega \nabla v \nabla u \dx,
\end{equation}
where the dual Besov space $(B^{2,2}_{\beta}(\del
\Omega))'=B^{2,2}_{-\beta}(\del \Omega)$ is introduced in
Ref.~\cite{JWBesov} (see~\ref{ApBS}).
\end{enumerate}
 \end{enumerate}
Let us also notice that the Green's formula~(\ref{FracGreen}) still holds
whatever $u\in H^1(\Omega)$ such that $\Delta u\in L^2(\Omega)$ and
$v\in H^1(\Omega)$.

\subsection{Well-posedness}

The above preliminaries allow us to prove the following Proposition:
\begin{proposition}\label{PropExis}
\begin{enumerate}
\item 
Let $\Omega=\Omega_+$ be a bounded domain in $\R^n$ with a closed
piecewise Lipschitz boundary $\del \Omega$ and $0<\lambda(x)\le
+\infty$ be a given function defined on $\del \Omega$.  By
$\Gamma_\infty$ is denoted the part of $\del \Omega$ such that
\begin{equation*}
\forall x\in \Gamma_\infty \quad \lambda(x)=+\infty,
\end{equation*}
in the such way that $\lambda\in L^\infty(\del \Omega\setminus
\Gamma_\infty)$.  Then the bilinear form
\begin{equation*}
a(u,v)=D_+\int_{\Omega_+} \nabla u_+ \nabla v_++D_-\int_{\Omega_-} \nabla u_- \nabla v_-+
 \int_{\del \Omega\setminus\Gamma_\infty }\lambda(x)(u_+-u_-)(v_+-v_-) d \sigma
\end{equation*}
is continuous and coercive on $V\times V$ with
\begin{multline}
 V=\{u\in L^2(\R^n)|\; u_+=u|_{\Omega_+}\in H^1(\Omega_+),\;u_-=u|_{\Omega_-}\in H^1(\Omega_-),\;\\
 u_+=u_- \hbox{ on } \Gamma_\infty\}. \label{spaceV}
\end{multline}

\item Let 
$\Omega=\Omega_+$ be a bounded $(\eps,\delta)$-domain in $\R^n$ with a
closed $d$-set boundary $\del \Omega$ and $\lambda\in C(\del \Omega)$
be a positive continuous function defined on $\del \Omega$.  By $m_d$
is denoted the $d$-measure on $\del \Omega$ (see
Definition~\ref{Defdset}).  Then the bilinear form
\begin{multline*}
   a(u,v)=D_+\int_{\Omega_+} \nabla u_+ \nabla v_++D_-\int_{\Omega_-} \nabla u_- \nabla v_-\\
   +
 \int_{\del \Omega}\lambda(x)\mathrm{Tr}(u_+-u_-)\mathrm{Tr}(v_+-v_-) d m_d
\end{multline*}
is continuous and coercive on $V\times V$ ($V$ is defined in
Eq.~(\ref{spaceV})).

\item 
Let $\Omega=\Omega_+$ be a bounded $(\eps,\delta)$-domain in $\R^n$
with a closed $d$-set boundary $\del \Omega$ and $\lambda(x)=+\infty$
for all $x\in \del \Omega$.  Then the bilinear form
\begin{equation*}
a(u,v)=D_+\int_{\Omega_+} \nabla u_+ \nabla v_++D_-\int_{\Omega_-} \nabla u_- \nabla v_-
\end{equation*}
is continuous and coercive on $V\times V$ with $V=H^1(\R^n)$.
\end{enumerate}

\end{proposition}
Consequently, we obtain the following theorem:
\begin{theorem}\label{ThWP}\textbf{(Well-posedness)}
In all cases from Proposition~\ref{PropExis} for all $u_0\in
H=L^2(\R^n)$ there exists a unique solution $u\in
C(\R^+_t,L^2(\R^n))\cap L^2(\R^+_t,V)$ of the variational problem
\begin{eqnarray}
 &&\forall v\in V\quad \frac{d}{dt}\langle  u,v \rangle_H + a(u,v)=0,\quad u(x,0)=u_0\in L^2(\R^n),\label{varpr}
\end{eqnarray}
where by $\langle \cdot,\cdot \rangle_H$ is denoted the inner product
in $H$.  In addition, this solution verifies the energy equality:
\begin{eqnarray}
 &&\frac{1}{2}\int_{\R^n}|u(t)|^2\dx+\int_0^t a(u,u) \ds=\frac{1}{2}\int_{\R^n}|u_0(x)|^2 \dx.\label{zvezda}
\end{eqnarray}
\end{theorem}

\begin{remark}
On one hand, any ``smooth enough'' solution of the
problem~(\ref{prb1})--(\ref{prb1end}) gives the solution of
Theorem~\ref{ThWP}.  On the other hand, any solution from
Theorem~\ref{ThWP} satisfies the relations~(\ref{prb1})--(\ref{prb1i})
and, in a weak sense (in the sense of the duality presented above),
satisfies the relations~(\ref{prb1c})--(\ref{prb1end}).
\end{remark}

Finally, we prove

\begin{theorem}\label{Thulamc}\textbf{(Continuity of $u_\lambda$ on $\lambda$ and the case $\lambda=\infty$)} 
Let $(\lambda_k)_{k\in \N}$ be a positive sequence converging to
$\lambda^*$ in $L^\infty(\del \Omega)$. Then the corresponding
sequence of the solutions $(u_{\lambda_k})_{k\in \N}$ of the
system~(\ref{prb1})--(\ref{prb1end}) converges strongly to
$u_\lambda^*$ in $C(\R^+_t,L^2(\R^n))\cap L^2(\R^+_t,V)$, $i.e.$,
$u_\lambda$ is continuous as a function of $\lambda$.
 
If $\lambda_k\to \infty$ in $L^\infty(\del \Omega)$, then
$u_{\lambda_k}\to u_\infty$ in $C(\R^+_t,L^2(\R^n))\cap L^2(\R^+_t,V)$
with $(u_\infty)_+=(u_\infty)_-$ on $\del \Omega$.  In this case,
$u_\infty\in C(\R^+_t,L^2(\R^n))\cap L^2(\R^+_t,H^1(\R^n))$ solves
\begin{equation}\label{infeq}
  \forall v\in H^1(\R^n)\quad \int_{\R^n}\del_t u_\infty v \dx + 
\int_{\R^n}D(x) \nabla u_\infty \nabla v \dx =0,   \quad u_\infty (x,0)=u_0\in L^2(\R^n),
\end{equation}
with $D(x)=\mathds{1}_{\Omega_+}D_++\mathds{1}_{\Omega_-}D_-$.
\end{theorem}
\textbf{Proof.}
Firstly we suppose that $\lambda^*$ is a finite bounded function on
$\del \Omega$ ($\|\lambda^*\|_{L^\infty(\del \Omega)}<\infty$).  Since
$u(0)=u_0$ does not depend on $\lambda$, the equality (\ref{zvezda})
implies that the sequence $(u_{\lambda_k})$ is bounded in
$C(\R^+_t,L^2(\R^n))\cap L^2(\R^+_t,V)$.
 
Therefore, due to the unicity of the solution for $\lambda^*$ and
the unicity of the weak limit, the convergence $\lambda_k \to
\lambda^*$ in $L^\infty(\del \Omega)$ implies
$u_{\lambda_k}\rightharpoonup u_{\lambda^*}$.  Since $u_k|_{\del
\Omega}\in B^{2,2}_{\beta}(\del \Omega)$ with
$\beta=1-\frac{n-d}{2}>0$, with the help of~(\ref{varpr}) and the
coercive behavior of $a(u,u)$,
\begin{equation*}
a(u,u)>\alpha\|u\|^2_V, \quad \hbox{for } \alpha>0,
\end{equation*}
we conclude that $u_{\lambda_k}\to u_{\lambda^*}$ in
$C(\R^+_t,L^2(\R^n))\cap L^2(\R^+_t,V)$.
 
In the case $\|\lambda_k\|_{L^\infty(\del \Omega)}\to +\infty$, we
find from~(\ref{zvezda}) that
\begin{eqnarray*}
  &&\forall k\in \N \quad \int_{\del \Omega}\lambda_k(x)\mathrm{Tr} ((u_k)_+-(u_k)_-)^2d m_d <\infty,\\
  &&\left(\int_{\del \Omega}\mathrm{Tr}((u_k)_+-(u_k)_-)^2 d m_d\right)^\frac{1}{2}\le \frac{1}{2\sqrt{\|\lambda_k\|_{L^\infty(\del \Omega)}}}\int_{\R^n}|u_0(x)|^2 \dx.
\end{eqnarray*}
Therefore, we obtain in this case that $(u_\infty)_+=(u_\infty)_-$ on
$\del \Omega$, where by $u_\infty$ we denote the limit of $u_k$ as
$\lambda\to +\infty$.  In addition, $u_\infty(t,\cdot)\in H^1(\R^n)$
and it is the solution of~(\ref{infeq}).
$\Box$

\section{Heat content localization to a small neighborhood of the  boundary}\label{secLoc}
 
s the initial condition is zero in $\R^n\setminus \Omega$, we have
\begin{equation}\label{NOmG}
N(t)=\int_\Omega(1-u(x,t))\dx= \operatorname{Vol}(\Omega)-\int_\Omega u(x,t)\dx,
\end{equation}
or equivalently, in terms of the Green function of the
problem~(\ref{prb1})--(\ref{prb1end}),
\begin{equation*}
 N(t)=\operatorname{Vol}(\Omega)-\int_{\Omega} \int_\Omega G(x,y,t)\dy \dx.
\end{equation*}
Let us show that it is sufficient to integrate only on a small
neighborhood of the boundary $\del \Omega$ to obtain the desired heat
content with an exponentially small error:

\begin{lemma}\label{Bound}
Let $F\subset \Omega$ be a non-empty open bounded set in $\R^n$, such
that $\operatorname{dist}(F,\del \Omega)=\eps>0$.  Then for $t\to +0$
and $u=u_+\mathds{1}_{\Omega}+u_-\mathds{1}_{\R^n\setminus\Omega}$ the
solution of~(\ref{prb1})--(\ref{prb1end}), associated with the Green
function $G(x,y,t)$,
\begin{enumerate}
\item  it holds
\begin{equation}\label{errF}
\begin{split}
 \int_F(1-u_+(x,t))\dx & =\int_F \left(1- \int_\Omega G(x,y,t)\dy\right) \dx  \\
& =O\left(\left(\frac{\eps}{\sqrt{4D_+t}} \right)^{n-2}e^{- \eps^2/(4D_+t)}\right).  \\
\end{split}
\end{equation}
\item for $\eps>2\sqrt{D_+ t}$ such that $\eps=O(\sqrt{t})$, 
there exists $\delta>0$ (a constant independent on time)   such that
 the heat content
$N(t)$ can be expressed as
\begin{equation}\label{ErrOm}
\begin{split}
  N(t) & =\int_{\R^n\setminus\Omega} u_-(x,t)\dx  \\  
  & =\int_{\Omega_\eps}\left(1-\int_{\Omega_\eps} G(x,y,t)\dy\right)\dx
+ O\left(e^{-\frac{1}{t^{\delta}}}\right),  \\
\end{split}
\end{equation}
where $\Omega_\eps$ is the $\eps$-neighborhood of $\del \Omega$.
\end{enumerate}
\end{lemma}

\textbf{Proof.} 
As it was shown, the problem~(\ref{prb1})--(\ref{prb1end}) has a
unique solution
$u=u_+\mathds{1}_{\Omega}+u_-\mathds{1}_{\R^n\setminus\Omega}$.  Let
$G(x,y,t)$ be the Green function so that
\begin{equation*}
u(x,t)=\int_\Omega G(x,y,t)\dy.
\end{equation*}
Thus, using the properties of $G$ such as $G\ge 0$ for all $(x,y,t)\in
\R^n\times \R^n\times\R_+$ and $\int_{\R^n}G(x,y,t)\dy=1$, we easily
see that 
\begin{equation*}
0\le \int_\Omega G(x,y,t)\dy=u(x,t)\le \int_{\R^n}G(x,y,t)\dy=1.
\end{equation*}

We notice that, by the assumption, $\lambda(x)>0$ is a regular function
on $\del \Omega$ and all other coefficients are constant.  By
definition $u_+$ is the solution of the system
\begin{eqnarray*}
 &&\left(\del_t -D_+ \Delta \right) u_+=0,\quad x\in \Omega\subset\R^n,\\
&&u_+|_{t=0}=1, \\
&&u_+|_{\del \Omega}=\left(u_--\frac{D_-}{\lambda}\frac{\del u_-}{\del n}\right)|_{\del \Omega}, \quad \lambda>0,
\end{eqnarray*}
which can be reformulated for $\hat{v}=1-u_+$
\begin{eqnarray*}
 &&\left(\del_t -D_+\Delta\right) \hat{v}=0,\quad x\in \Omega\subset\R^n,\\
&&\hat{v}|_{t=0}=0, \\
&&(1-\hat{v})|_{\del \Omega}=\left(u_--\frac{D_-}{\lambda}\frac{\del u_-}{\del n}\right)|_{\del \Omega},
\end{eqnarray*}
where $0 \le u_-\le 1$ for all $t$. Moreover, as $0\le \hat{v} \le 1$,
it follows that
\begin{equation*}
0\le \left(u_--\frac{D_-}{\lambda}\frac{\del u_-}{\del n}\right)|_{\del \Omega}\le 1
\end{equation*}
and, as $u_-$ is increasing in time on $\del \Omega$, then $\hat{v}$
is decreasing in time on $\del \Omega$. Therefore, $\hat{v} \le v$,
where $v$ is the solution of the following problem:
\begin{eqnarray*}
 &&\left(\del_t -D_+\Delta\right) v=0,\quad x\in \Omega\subset\R^n,\\
&&v|_{t=0}=0, \\
&&v|_{\del \Omega}=1,
\end{eqnarray*}
Thus, as in Ref.~\cite{Falconer} (p.231 Lemma 12.7) for $n=2$, but
now in $\R^n$ ($n\ge 2$), we find that for the ball $\Omega=U_r(z)$
centered at $z$ and of radius $r$, the solution satisfies as $t\to +0$
\begin{equation*}
v(z,t)\le C\left(\frac{r}{\sqrt{4D_+t}}\right)^{n-2} \exp\left(-\frac{r^2}{4D_+t}\right) ,
\end{equation*}
with a constant $C>0$ depending only on $n$ ($C$ can be explicitly
obtained by the integration by parts in the generalized spherical
coordinates in $\R^n$, where the coefficient
$\left(\frac{r}{\sqrt{4D_+t}}\right)^{n-2}$ corresponding to the leading
term as $t\to +0$, appears from the integral
$\int_{\frac{r}{\sqrt{4Dt}}}^{+\infty}e^{-w^2}w^{n-1}dw$).
Consequently (see Ref.~\cite{Falconer} Corollary 12.8 p.232), for
$z\in \textrm{int} \{ \Omega\}$ and $t\to +0$ we find
\begin{equation*}\label{vest}
v(z,t)\le C \left(\frac{\dist(z,\del \Omega)}{\sqrt{4D_+t}}\right)^{n-2}\exp\left(-\frac{\dist(z,\del \Omega)^2}{4D_+t}\right).
\end{equation*}
Then we immediately obtain Eq.~(\ref{errF}) by integration.

For $n=2$ we obtain directly the exponential decay in Eq.~(\ref{errF})
for all $\eps>0$.  If $n>2$, we still have the exponential decay for a
small constant $\alpha>0$ depending only on $\eps$:
\begin{equation*}
O\left(\left(\frac{\eps}{\sqrt{4D_+t}}
\right)^{n-2}e^{- \eps^2/(4D_+t)}\right)=O(e^{- \alpha/t}).
\end{equation*}

Note that $O\left(e^{-\eps^2/(4D_+t)}\right)$ gives an exponentially
small remaining term iff $\eps =2\sqrt{D_+}~t^{\frac{1}{2}-\delta_0}$
for a constant $\delta_0>0$.  For small enough $\delta_0$    we have
$\eps=O(\sqrt{4D_+t})$, also knowing that $\eps>\sqrt{4D_+t}$.

So, for this $\eps$, we split $\Omega$ in two parts:
$\Omega_{\eps}$, the neighborhood of $\del \Omega$ such that
$\dist(x,\del \Omega)\le\eps$, and $\Omega\setminus
\Omega_{\eps}$.  For all $F\subseteq\Omega\setminus
\Omega_{\eps}$, $\dist(F,\del \Omega)>\eps>2\sqrt{D_+t}$, we have
\begin{equation*}
\int_F \left(1- \int_\Omega G(x,y,t)\dy\right) \dx=O\left(e^{-c(F)/t^{\delta(F)}}\right),
\end{equation*}
where $c(F)$ and $\delta(F)$ are positive constants depending only on
the distance between $F$ and $\del\Omega$ and the dimension $n$.

To complete the proof of the second statement, we first find that
\begin{eqnarray*}
  &&N(t)=\int_{\R^n\setminus\Omega} u_-(x,t)\dx\\
&&=\int_{\R^n\setminus\Omega}\int_\Omega G(x,y,t)\dy \dx =\int_{\R^n}\int_\Omega G(x,y,t)\dy \dx-\int_\Omega \int_\Omega G(x,y,t)\dy \dx\\
&&= \operatorname{Vol}(\Omega)-\int_\Omega \int_\Omega G(x,y,t)\dy \dx.
\end{eqnarray*}
For $\Omega=\Omega_\eps \cup (\Omega \setminus\Omega_\eps)$ we can
write
\begin{eqnarray*}
 &&\int_\Omega \int_\Omega G (x,y,t)\dy \dx=\left(\int_{\Omega_\eps}\int_{\Omega}+\int_{\Omega \setminus\Omega_\eps} \int_{\Omega} \right) G (x,y,t)\dy \dx\\
&&=\operatorname{Vol}(\Omega \setminus\Omega_\eps)+\int_{\Omega_\eps}\int_{\Omega}G (x,y,t)\dy \dx+O\left(e^{-\frac{1}{t^{\delta}}}\right)\\
&& =\operatorname{Vol}(\Omega)-\operatorname{Vol}(\Omega_\eps)+\int_{\Omega_\eps}\int_{\Omega_\eps}G (x,y,t)\dy \dx \\
&&+ \int_{\Omega_\eps}\int_{\Omega \setminus\Omega_\eps}G (x,y,t)\dy \dx+O\left(e^{-\frac{1}{t^{\delta}}}\right).
\end{eqnarray*}
Moreover, 
\begin{eqnarray*}
 &&\int_{\Omega_\eps}\int_{\Omega \setminus\Omega_\eps}G (x,y,t)\dy \dx=\int_{\Omega}\int_{\Omega \setminus\Omega_\eps}G (x,y,t)\dy \dx
-\int_{\Omega \setminus\Omega_\eps}\int_{\Omega \setminus\Omega_\eps}G (x,y,t)\dy \dx\\
&&=\operatorname{Vol}(\Omega\setminus\Omega_\eps)+O\left(e^{-\frac{1}{t^{\delta}}}\right)-
\int_{\Omega \setminus\Omega_\eps}\int_{\Omega \setminus\Omega_\eps}G (x,y,t)\dy \dx
\end{eqnarray*}
and since
\begin{eqnarray*}
 &&\int_{\Omega \setminus\Omega_\eps}\int_{\Omega \setminus\Omega_\eps}G (x,y,t)\dy \dx-\operatorname{Vol}(\Omega\setminus\Omega_\eps)\\
&&\le  \int_{\Omega \setminus\Omega_\eps}\int_{\Omega}G (x,y,t)\dy \dx-\operatorname{Vol}(\Omega\setminus\Omega_\eps)=O\left(e^{-\frac{1}{t^{\delta}}}\right),
\end{eqnarray*}
we conclude that
\begin{eqnarray*}
&&-\int_{\Omega_\eps}\int_{\Omega \setminus\Omega_\eps}G (x,y,t)\dy \dx=O\left(e^{-\frac{1}{t^{\delta}}}\right)
\end{eqnarray*}
and finally
\begin{eqnarray*}
 &&N(t)=\operatorname{Vol}(\Omega_\eps)-\int_{\Omega_\eps}\int_{\Omega_\eps}G (x,y,t)\dy \dx+O\left(e^{-\frac{1}{t^{\delta}}}\right),
\end{eqnarray*}
that completes the proof.  $\Box$

A variant of Lemma~\ref{Bound} for $n=2$ can be found in
Ref.~\cite{Fleckinger}, where the heat localization near the
boundary is also called by the principle of ``not feeling the
boundary''.  In addition, we can consider the case of the distinct
parts of the boundary:
\begin{corollary}\label{Sledstvie}
Let $X$ and $Y$ be different closed parts of $\del \Omega$ such that
$\operatorname{dist}(X,Y)>2\epsilon$, where
$\eps=O(\sqrt{t})>2\sqrt{D_+}~t$.  Let $U_r(X)=\{x\in \R^n|\;
d(x,X)<r\}$ be the open neighborhood of $X$ of size $r>0$.  Consider
$u_+$ and $\hat{u}_+$ as the respective solutions of the following
systems:
\begin{eqnarray*}
 &&\del_t u_+ -D_+\triangle  u_+=0,\quad x\in \Omega\subset\R^n,\nonumber\\
&&u_+|_{t=0}=1,\nonumber \\
&&u_+|_{\del \Omega}=\left(u_--\frac{D_-}{\lambda}\frac{\del u_-}{\del n}\right)|_{\del \Omega}, \quad \lambda>0,\label{s1}
\end{eqnarray*}
\begin{eqnarray*}
 &&\del_t\hat{u}_+ -D_+\triangle \hat{u}_+ =0,\quad x\in \Omega\subset\R^n,\nonumber\\
&&\hat{u}_+|_{t=0}=1,\nonumber \\
&&\hat{u}_+|_{\del \Omega\cap \overline{U(X)}}=
\left(u_--\frac{D_-}{\lambda}\frac{\del u_-}{\del n}\right)|_{\del \Omega\cap \overline{U(X)}}, \quad \lambda>0\nonumber\\
&&\hat{u}_+|_{\del \Omega\setminus(\del \Omega\cap\overline{U(X)})}=1,\label{s2}
\end{eqnarray*}
where $U(X)$ is an open neighborhood of $X$ of a radius strictly
greater than $2\eps$: $U_{2\eps}(X)\subset U(X)$.

Then there exists  $\delta>0$ such that
\begin{equation*}
      \int_{U_\eps(X)}|u_+-\hat{u}_+|\dx=O\left(e^{-\frac{1}{t^{\delta}}}\right).
     \end{equation*}
Moreover, if $\tilde{u}_+$ is the solution of the system:
\begin{eqnarray*}
 &&\del_t\tilde{u}_+ -D_+\triangle  \tilde{u}_+ =0,\quad x\in \Omega\subset\R^n,\nonumber\\
&&\tilde{u}_+|_{t=0}=1, \nonumber\\
&&\tilde{u}_+|_{Y}=\left(u_--\frac{D_-}{\lambda}\frac{\del u_-}{\del n}\right)|_{Y}, \quad \lambda>0\nonumber\\
&&\tilde{u}_+|_{\del \Omega\setminus Y}=1,\label{s3}
\end{eqnarray*}
then \begin{equation*}
      \int_{U_\eps(X)}(1-\tilde{u}_+)\dx=\int_{\Omega\setminus U_\eps(Y)}(1-\tilde{u}_+)\dx=O\left(e^{-\frac{1}{t^{\delta}}}\right).
     \end{equation*}
\end{corollary}
The proof of Corollary~\ref{Sledstvie} follows from the proof of the
first statement of Lemma~\ref{Bound}.

Note that the continuity of $u$ on $\lambda$ (see
Theorem~\ref{Thulamc}) and the localization of the heat content near
the boundary allow one to consider mixed boundary conditions:
\begin{theorem}\label{ThMixed}
Let $\Omega$ be a bounded $(\eps,\delta)$-domain (see
Section~\ref{secWP}) with a closed connected $d$-set boundary $\del
\Omega=\Gamma_0\sqcup\Gamma_\lambda\sqcup\Gamma_\infty$.  Let
$\lambda\in C(\Gamma_0\sqcup\Gamma_\lambda)$ such that
\begin{equation*}
\lambda(x)=\left\{\begin{array}{ll}
                      0,& x\in \Gamma_0,\\
                      0<f(x)<\infty,&x\in \Gamma_\lambda,\\
                      +\infty,& x\in \Gamma_\infty
                     \end{array}\right.
\end{equation*}
and $\eps=O(\sqrt{t}) > \sqrt{4D_+t}$.  We assume that the connection
between different types of boundary is performed in the continuous way
(see Theorem~\ref{Thulamc}) such that the solution $u$ remains
continuous as a function of $\lambda$.

We split the $\eps$-interior Minkowski sausage of $\del \Omega$ into
disjoint subsets
\begin{equation*}
\Omega_\eps=\Omega_\eps^{\Gamma_0}\sqcup\Omega_\eps^{\Gamma_\lambda}\sqcup\Omega_\eps^{\Gamma_\infty}
\end{equation*}
such that each subset $\Omega_\eps^{\Gamma}$ is contained in the
$\eps$-interior Minkowski sausage of $\Gamma$ ($\Gamma\subset \del
\Omega$).  Then, for $\delta>0$ from Lemma~\ref{Bound}, the heat
content of the problem~(\ref{prb1})--(\ref{prb1end}),
\begin{equation*}
N(t)=\int_{\Omega} (1-u(x,t))\dx=\int_{\Omega_\eps}(1-u(x,t))\dx+O(e^{-\frac{1}{t^{\delta}}}),
\end{equation*}
can be found as a sum of two heat contents:
\begin{equation*}
N(t)=\int_{\Omega_\eps^{\Gamma_\lambda}}(1-u(x,t))\dx+ \int_{\Omega_\eps^{\Gamma_\infty}}(1-u(x,t))\dx+O(e^{-\frac{1}{t^{\delta}}}).
\end{equation*}
\end{theorem}

In order to locally approximate the solution of the
problem~(\ref{prb1})--(\ref{prb1end}) by considering the problem with
coefficients frozen on a fixed boundary point, according to
Corollary~\ref{Sledstvie}, we also obtain the following proposition:
\begin{proposition}\label{Local}
Let $\sigma$ be a fixed point of the boundary $\del \Omega$ and let
define
\begin{equation}\label{Blee}
B_{l\eps,\eps}=U_{l\eps}(\sigma)\cap(\Omega_\eps\cup\Omega_{-\eps}) \quad \hbox{ for }\quad l\in \N,
\end{equation}
where $U_{l\eps}(\sigma)\subset \R^n$ is a ball of radius $l\eps$
centered at $\sigma$, $\eps$ is defined in Corollary~\ref{Sledstvie}.
Let $\phi_\sigma\in \mathring{C}^\infty(B_{4\eps,\eps}(\sigma))$ be a
smooth cut-off function with a compact support on
$B_{4\eps,\eps}(\sigma)$:
\begin{equation}\label{phisigma}
 \phi_\sigma(x)=\left\{\begin{array}{ll}
                     1&x\in B_{3\eps,\eps}(\sigma),\\
\hbox{a smooth function }0\le \eta<1&x\in B_{4\eps,\eps}(\sigma)\setminus B_{3\eps,\eps}(\sigma),\\
0&x\in \Omega\setminus B_{4\eps,\eps}(\sigma)
                    \end{array}
 \right.
\end{equation}
If $u$ is the solution of the problem~(\ref{prb1})--(\ref{prb1end}),
then $\phi_\sigma u$ is the solution of the following problem:
\begin{eqnarray}
 &&\del_t\left(\phi_\sigma u_\pm \right)-D_\pm \triangle\left(\phi_\sigma u_\pm \right)=
\left\{\begin{array}{ll}
    -(1-u_\pm)D_\pm \triangle\phi_\sigma) &x\in B_{4\eps,\eps}(\sigma)\setminus B_{3\eps,\eps}(\sigma),\\
								  0& \hbox{elsewhere in } \Omega,
                                                                 \end{array} \right.\label{pfi1}\\
&&\left(\phi_\sigma u_\pm \right)|_{t=0}=\mathds{1}_{\Omega}(x)\phi_\sigma(x), \\
&& D_- \frac{\del(\phi_\sigma u_-)}{\del n}|_{\del \Omega}=\lambda(x)\phi_\sigma(x)(u_--u_+)|_{\del \Omega}, \\
 && D_+ \frac{\del(\phi_\sigma u_+)}{\del n}|_{\del \Omega}=D_- \frac{\del (\phi_\sigma u_-)}{\del n}|_{\del \Omega}.\label{pfiEnd}
\end{eqnarray}
Therefore, there exists $\delta>0$ such that
\begin{eqnarray*}
&&\int_{B_{2\eps,\eps}(\sigma) }|u-\phi_\sigma u|\dx=O\left(e^{-\frac{1}{t^{\delta}}}\right),
\end{eqnarray*}
and if $\phi_\sigma u^\sigma$ is the solution of the
problem~(\ref{pfi1})--(\ref{pfiEnd}) with frozen coefficients in the
boundary point $\sigma$, then
\begin{equation}
 \label{SmislBee}
 \int_{\Omega\setminus B_{\eps,\eps}(\sigma) }\phi_\sigma (1-u^\sigma)\dx=O\left(e^{-\frac{1}{t^{\delta}}}\right).
\end{equation}
\end{proposition}
 
\section{Local coordinates for a regular $\del \Omega\in C^3$}\label{secLC}

In order to prove Eq.~(\ref{EqEtoile}) for a large class of
$(\eps,\delta)$-compact connected domains $\Omega$ in $\R^n$, we first
prove it for the case of domains with regular boundary $\del \Omega\in
C^\infty$ or at least in $C^3$.  As $\Omega$ is compact, for all types
of connected $\del \Omega$, the volume of $\Omega$ is finite and,
therefore, the volume of the $\eps$-neighborhood of $\del \Omega$ in
$\Omega$ is also finite and can be approximated by a sequence of
volumes of Minkowski sausages with regular boundaries (the same
argument was used in Ref.~\cite{Fleckinger} p.378).

Let us consider the regular boundary $\del \Omega\in C^3$.

Given a positive $\eps>0$ provided in Lemma~\ref{Bound}, we denote by
$\Omega_\eps$ and $\Omega_{-\eps}$ the open $\eps$-neighborhoods of
$\del \Omega$ in $\Omega$ and in $\R^n\setminus\Omega$, respectively.

According to Eq. (\ref{ErrOm}) and the regularity of the boundary
$\del \Omega$, we can decompose $\Omega_{\eps}\cup\del
\Omega\cup\Omega_{-\eps}=\bigsqcup_{i=1}^I B_{i,\eps}$ ($I$ is a
finite integer because $\overline{\Omega}_+\cup\Omega_{-\eps}$ is a
compact domain) in such way that on each $B_{i,\eps}$ it is possible
to introduce the local coordinates.  In addition, we assume that for
all $i=1,\ldots,I$ there exists $\sigma_i\in \del \Omega \cap
B_{i,\eps}$ such that $B_{i,\eps}\subset B_{2\eps,\eps}(\sigma_i)$
(see Eq.~(\ref{Blee}) for the definition).  Due to
Proposition~\ref{Local}, the last assumption ensures that
\begin{equation*}
\int_{B_{i,\eps}}(1-u)\dx=\int_{B_{i,\eps}}\phi_{\sigma_i}(1-u)\dx+O\left(e^{-\frac{1}{t^{\delta}}}\right).
\end{equation*}

For all $i$ we perform the change of the space variables 
$(x_1,\ldots,x_n)\in B_{i,\eps}$ to the local
coordinates $(\theta_1,\ldots,\theta_{n-1},s)$ by the formula
\begin{equation}\label{localc}
  x=\hat{x}(\theta_1,\ldots,\theta_{n-1})-sn(\theta_1,\ldots,\theta_{n-1}) \;\quad
 \left\{\begin{array}{l}
 0<s<\eps\;\hbox{for } x\in B_{i,\eps}\cap\Omega_{\eps}\\
 -\eps<s<0 \; \hbox{for } x\in  B_{i,\eps}\cap\Omega_{-\eps}
                       \end{array}\right.,
\end{equation}
where $\hat{x}(\theta_1,\ldots,\theta_{n-1})\in \del \Omega$ and $x$,
$\hat{x}$ and $n$ are the vectors in $\R^n$ such that
\begin{equation*}
\left\{\frac{\del\hat{x} }{\del\theta_1 }, \ldots,\frac{\del\hat{x} }{\del\theta_{n-1} },n\right\}
\end{equation*}
is an orthonormal basis in $\R^n$.

In what follows we denote $B_{i,\eps}\cap\Omega_{\eps}$ by
$\Omega_{i,+\eps}$ and $B_{i,\eps}\cap\Omega_{-\eps}$ by
$\Omega_{i,-\eps}$ respectively.  In each of two regions,
$\Omega_{i,+\eps}$ and $\Omega_{i,-\eps}$, the change of variables
$(x_1,\ldots,x_n)\mapsto (\theta_1,\ldots,\theta_{n-1},s)$ is a local
$C^1$-diffeomorphism.

In local coordinates  $\del \Omega$ is described by $s=0$. 

Thus, Eq.~(\ref{ErrOm}) becomes
\begin{eqnarray}
  &&N(t)=\sum_{i=1}^I \int_{\Omega_{i,+\eps}} (1-u(x,t))\dx+O(e^{-\frac{1}{t^{\delta}}}). \label{nloc}
\end{eqnarray}
  
Denoting $\theta=(\theta_1,\ldots,\theta_{n-1})$, the integration
domain $\Omega_{i,+\eps}$ in~(\ref{nloc}) becomes
\begin{equation*}
\Omega_{i,+\eps}=\{0<s<\eps, \quad \theta\in \del \Omega\cap\overline{\Omega}_{i,+\eps}\},
\end{equation*}
which is actually a parallelepiped neighborhood
$(\del \Omega\cap\overline{\Omega}_{i,+\eps})\times]0,\eps[.$

For this change of variables we have
\begin{eqnarray*}
&&|\nabla_x s|^2=1,\quad \nabla_x s \nabla_x \theta_i=0, \quad \nabla_x \theta_j\nabla_x \theta_i=
\frac{\delta_{ij}}{(1-sk_i)^2} \quad \hbox{for } i,j=1,\ldots,n-1,\\
 &&\nabla u_\pm\nabla \phi_\pm=\frac{\del u_\pm}{\del s}\frac{\del \phi_\pm}{\del s}+
\sum_i^{n-1}\frac{\del u_\pm}{\del \theta_i}\frac{\del \phi_\pm}{\del \theta_i}\frac{1}{(1-sk_i)^2},
\end{eqnarray*}
and therefore, using twice the integration by parts and the notations
\begin{equation}\label{Jacobian}
 |J(s,\theta)|=\prod_{i=1}^{n-1}(1-sk_i)
\end{equation}
for the Jacobian and $k_i=k_i(\theta_1,\ldots,\theta_{n-1})$ of the
principal curvatures for $\del \Omega$ curving away the outward normal
$n$ to $\del \Omega$ like in the case of the sphere, we find that for
all test functions $\phi=(\phi_+,\phi_-)\in V|_{B_{i,\eps}}$
\begin{eqnarray*}
 &&\int_{B_{i,\eps}}\del_t u ~|J(s,\theta)|~\phi ~\ds d \theta_1\cdots d \theta_{n-1}
\\
&&-\int_{\Omega_{i,+\eps}} \hspace*{-1mm} \left[\frac{\del }{\del s}\left(D_+|J(s,\theta)|\frac{\del u_+ }{\del s} \right)
+\sum_{i=1}^{n-1}\frac{\del }{\del \theta_i}\left(\frac{D_+|J(s,\theta)|}{(1-sk_i)^2} 
\frac{\del u_+}{\del \theta_i}\right)\right]\phi_+ \ds d \theta_1\cdots d \theta_{n-1}\\
&&-\int_{\Omega_{i,-\eps}} \hspace*{-1mm} \left[\frac{\del }{\del s}\left(D_-|J(s,\theta)|\frac{\del u_- }{\del s} \right)+
\sum_{i=1}^{n-1}\frac{\del }{\del \theta_i}\left(\frac{D_-|J(s,\theta)|}{(1-sk_i)^2} 
\frac{\del u_-}{\del \theta_i}\right)\right]\phi_- \ds d \theta_1\cdots d \theta_{n-1}\\
&&+ \int_{s=0}\lambda(\theta)(u_+-u_-)(\phi_+-\phi_-) d\theta=0.
\end{eqnarray*}
The regularity of the boundary ensures that the principal curvatures
$k_i(\theta)$ are at least in $C^1(\del \Omega\cap \del B_{i,\eps}).$

Therefore, the problem~(\ref{prb1})--(\ref{prb1end}) locally becomes
\begin{eqnarray}
 &&\frac{\del}{\del t}u_+-D_+\left(\frac{\del^2}{\del s^2}+\sum_{i=1}^{n-1}\frac{\del^2}{\del \theta^2_i}\right) u_+=
 D_+\sum_{i=1}^{n-1}  \frac{sk_i(\theta)}{1-sk_i(\theta)}
\left( 1+\frac{1}{1-sk_i(\theta)}\right)\frac{\del^2u_+}{\del \theta^2}\nonumber\\
 && -D_+\left(\sum_{i=1}^{n-1}k_i(\theta)+s\sum_{i=1}^{n-1} \frac{k_i^2(\theta)}{1-sk_i(\theta)} \right)
\frac{\del u_+}{\del s}\nonumber\\
 &&+\frac{D_+}{|J(s,\theta)|}\sum_{i=1}^{n-1}\frac{\del}{\del \theta_i}
\left(\frac{|J(s,\theta)|}{(1-sk_i(\theta))^2} \right) \frac{\del u_+}{\del \theta_i},\;0<s<\eps, \quad 
\theta\in (\del \Omega\cap\overline{\Omega}_{i,+\eps})\label{eq47}\\
  &&\frac{\del}{\del t}u_--D_-\left(\frac{\del^2}{\del s^2}+\sum_{i=1}^{n-1}\frac{\del^2}{\del \theta^2_i}\right) u_-=
 D_-\sum_{i=1}^{n-1}  \frac{sk_i(\theta)}{1-sk_i(\theta)}\left( 1+\frac{1}{1-sk_i(\theta)}\right)
\frac{\del^2u_-}{\del \theta^2} \nonumber\\
&&-D_-\left(\sum_{i=1}^{n-1}k_i(\theta)+s\sum_{i=1}^{n-1} \frac{k_i^2(\theta)}{1-sk_i(\theta)} \right)
\frac{\del u_-}{\del s}\nonumber\\
&&+\frac{D_-}{|J(s,\theta)|}\sum_{i=1}^{n-1}\frac{\del}{\del \theta_i}\left(\frac{
|J(s,\theta)|}{(1-sk_i(\theta))^2} \right) \frac{\del u_-}{\del \theta_i}, \; -\eps<s<0, 
\quad \theta\in (\del \Omega\cap\overline{\Omega}_{i,+\eps}),\label{eq48}\\
&&u_+|_{t=0}=1,\; u_-|_{t=0}=0,\label{ini}\\
 && D_-\frac{\del u_-}{\del s}|_{s=-0}=\lambda(\theta) (u_--u_+)|_{s=0}, \label{boundar}\\
&& D_+ \frac{\del u_+}{\del s}|_{s=+0}=D_-\frac{\del u_-}{\del s}|_{s=-0}.\label{endsys}
\end{eqnarray}
We emphasize that the problem~(\ref{eq47})--(\ref{endsys}) should be
considered as the trace of Eqs.~(\ref{prb1})--(\ref{prb1end}) on
$B_{i,\eps}$ in the sense of the problem~(\ref{pfi1})--(\ref{pfiEnd})
with $\phi_{\theta_i}\equiv 1$ on $B_{i,\eps}.$

Therefore, we can rewrite~(\ref{nloc}) in new coordinates and use the
parallelepiped property of $\Omega_{i,+\eps}$ in the space of
variables $(s,\theta)$:
\begin{eqnarray*}
  &&N(t)=\sum_{i=1}^I\int_{\Omega_{i,+\eps}}(1-u(s,\theta,t))|J(s,\theta)|dsd\theta+O(e^{-\frac{1}{t^{\delta}}})\\
  &&=\sum_{i=1}^I\int_{\del \Omega\cap\overline{\Omega}_{i,+\eps}}d\theta \int_{[0,\eps]}ds (1-u(s,\theta,t))|J(s,\theta)|+O(e^{-\frac{1}{t^{\delta}}}).\\
\end{eqnarray*}
Since this local representation holds for all $i$ (the form of the
problem~(\ref{eq47})--(\ref{endsys}) is the same for all $i$) and
$\sum_{i=1}^I\int_{\del
\Omega\cap\overline{\Omega}_{i,+\eps}}d\theta=\int_{\del \Omega} d
\theta$, we can formally write
\begin{eqnarray}
   &&N(t)=\int_{\del \Omega}d\theta \int_{[0,\eps]}ds (1-u(s,\theta,t))|J(s,\theta)|+O(e^{-\frac{1}{t^{\delta}}}),\label{nnloc}
\end{eqnarray}
where $u$ is the solution of~(\ref{eq47})--(\ref{endsys}) in
$]-\eps,\eps[\times \del \Omega$ in the local sense, as explained
previously.

\section{Approximation of the heat content by  solutions of one dimensional problems (for a regular boundary)}\label{secAprox}

We denote by $\tilde{G}(s_1,\theta_1,s_2,\theta_2,t)$ the Green
function of the problem~(\ref{eq47})--(\ref{endsys}) in $\del
\Omega\times]-\eps,\eps[$.  Let us fix a boundary point
$(0,\theta_0)$.

We denote by $G^{\theta_0}$ the Green function corresponding to the
following constant coefficient problem, considered as a local trace
problem, i.e. in the sense of the problem~(\ref{pfi1})--(\ref{pfiEnd})
with $\phi_{\theta_i}\equiv 1$ on $B_{i,\eps}$:
\begin{eqnarray}
&&\frac{\del}{\del t}u_+-D_+\left(\frac{\del^2}{\del s^2}+\sum_{i=1}^{n-1}\frac{\del^2}{\del \theta^2_i}\right) 
u_+ =0,\; 0<s<\eps\label{consteq1}\\
&&\frac{\del}{\del t}u_--D_-\left(\frac{\del^2}{\del s^2}+\sum_{i=1}^{n-1}\frac{\del^2}{\del \theta_i^2}\right) 
u_- =0, \; -\eps<s<0\\
&&u_+|_{t=0}=1,\; u_-|_{t=0}=0,\\
&& D_-\frac{\del u_-}{\del s}|_{s=-0}=\lambda(\theta_0)(u_--u_+)|_{s=0},\nonumber \\
&& D_+\frac{\del u_+}{\del s}|_{s=+0}=D_- \frac{\del u_-}{\del s}|_{s=-0}.\label{consteqN}
\end{eqnarray}

Next, let 
\begin{equation*}
G^{\theta_0}_{\R^n}(s_1,\theta_1,s_2,\theta_2,t)=
     \mathds{1}_{\{s_1>0\}} G^{\theta_0}_{++}(s_1,\theta_1,s_2,\theta_2,t)+\mathds{1}_{\{s_1<0\}}
G^{\theta_0}_{-+}(s_1,\theta_1,s_2,\theta_2,t)
\end{equation*}
be the Green function of the constant coefficient problem in the half
space, explicitly obtained in~\ref{SecGreen}.  Then, according to
Ref.~\cite{McKean} p.48--49, due to Varadhan's bound property of
Green functions, in $U_\eps (0,\theta_0)$ the difference between the
Green function $\phi_{\theta_0}G^{\theta_0}$ of the
problem~(\ref{consteq1})--(\ref{consteqN}) and the analogous Green
function in $\R^n$, $G^{\theta_0}_{\R^n}$, is exponentially small:
\begin{equation*}
|(\phi_{\theta_0}G^{\theta_0}-G^{\theta_0}_{\R^n})|_{U_\eps (0,\theta_0)\times U_\eps (0,\theta_0) }|=O\left(e^{-\frac{1}{t^{\delta}}}\right).
\end{equation*}
Therefore, following the ideas of McKean and Singer\cite{McKean}
(p.49), we approximate $\tilde{G}$ by the Green function
$G^{\theta_0}$ with the frozen coefficients on $(0,\theta_0)$, whose
replacement by $G^{\theta_0}_{\R^n}$ yields only an exponentially
small error.

For an abstract operator Cauchy problem
\begin{eqnarray}
\label{opereqCH}
 &&\frac{\del}{\del t}u-A u = \mathcal{R} u,\\
 &&u|_{t=0}=u_0\nonumber
\end{eqnarray}
the solution $u$ can be found by the Duhamel formula
\begin{equation}\label{eqCauchy}
 u(t)=e^{-tA} u_0+\int_0^te^{-(t-\tau)A}\mathcal{R} u (\tau)\dtau.
\end{equation}

Therefore, by the Duhamel formula, locally, we have the following
infinite expansion
\begin{eqnarray}
&& u_+(s,\theta_0,t)=\int_{\Omega_\eps}d\theta_1 d s_1G_{++}^{\theta_0}(s,\theta_0,s_1,\theta_1,t)\nonumber\\
&&+\int_0^t\dtau\int_{\Omega_\eps}d\theta_1 d s_1G_{++}^{\theta_0}(s,\theta_0,s_1,\theta_1,t-\tau)\cdot\nonumber\\
&&\cdot\int_{\Omega_\eps}d\theta_2 d s_2\mathcal{R} G_{++}^{\theta_0}(s_1,\theta_1,s_2,\theta_2,\tau) \nonumber\\
 &&+\int_0^t\dtau\int_{\Omega_\eps}d\theta_1 d s_1 G_{++}^{\theta_0}(s,\theta_0,s_1,\theta_1,t-\tau) \cdot\nonumber\\
&&\cdot\int_0^\tau \dtau_1 \int_{\Omega_\eps}d\theta_2 d s_2 \mathcal{R}G_{++}^{\theta_0}(s_1,\theta_1,s_2,\theta_2,\tau-\tau_1)\cdot
\nonumber\\
 &&\cdot\int_{\Omega_\eps}d\theta_3 d s_3\mathcal{R}G_{++}^{\theta_0}(s_2,\theta_2,s_3,\theta_3,\tau_1)+\ldots+O\left(e^{-\frac{1}{t^{\delta}}}\right)
\label{uCauchyexpansion}
\end{eqnarray}
where the operator $\mathcal{R}$ is defined by  
\begin{eqnarray}
 &&\mathcal{R}=\mathcal{R}_{s_1}(s_1,\theta_1)
+\mathcal{R}_{\theta_1}(s_1,\theta_1),\label{R}\\
&&\mathcal{R}_s(s,\theta)= R(s,\theta)\frac{\del }{\del s}=-D_+
\left(\sum_{i=1}^{n-1}k_i(\theta)+s\sum_{i=1}^{n-1} \frac{k_i^2(\theta)}{1-sk_i(\theta)} \right)\frac{\del }{\del s},\label{Restins}\\
&&\mathcal{R}_{\theta}(s,\theta)=
\sum_{i=1}^{n-1}  \frac{D_+sk_i(\theta)}{1-sk_i(\theta)}\left( 1+\frac{1}{1-sk_i(\theta)}\right)
\frac{\del^2}{\del \theta^2_i}\nonumber\\
&&+\frac{D_+}{|J(s,\theta)|}\sum_{i=1}^{n-1}\frac{\del}{\del \theta_i}\left(\frac{
|J(s,\theta)|}{(1-sk_i(\theta))^2} \right) \frac{\del }{\del \theta_i}.\label{Restinteta}
\end{eqnarray}

We substitute Eq.~(\ref{uCauchyexpansion}) into Eq.~(\ref{nnloc}) with
$\theta=\theta_0$ and prove the following theorem:
\begin{theorem}\label{TH1}
Let 
\begin{equation*}
\hat{u}=\left\{\begin{array}{ll}
                 \hat{u}_+,& 0<s<\eps\\
\hat{u}_-,& -\eps<s<0
                \end{array}
 \right. 
\end{equation*}
be the solution of the one-dimensional problem
\begin{eqnarray}
 &&\frac{\del}{\del t}\hat{u}-D_\pm\frac{\del^2}{\del s^2} \hat{u}=\mathcal{R}_s(s,\theta_0)\hat{u} \quad -\eps<s<\eps, \;\theta\equiv\theta_0,\label{opereqCH1d}\\
&&\hat{u}|_{t=0}=\mathds{1}_{0<s<\eps}(s),\nonumber\\
&& D_- \frac{\del \hat{u}_-}{\del s}|_{s=-0}=\lambda(\theta_0)(\hat{u}_--\hat{u}_+)|_{s=0}, \\
&& D_+\frac{\del \hat{u}_+}{\del s}|_{s=+0}=D_-\frac{\del \hat{u}_-}{\del s}|_{s=-0},\label{opereqCH1dN}
\end{eqnarray}
obtained from~(\ref{eq47})--(\ref{endsys}) setting
$\theta\equiv\theta_0$ ($\mathcal{R}_s(s,\theta_0)$ is given
by~(\ref{Restins})).  Then the heat content $N(t)$, defined
in~(\ref{nnloc}), satisfies
\begin{equation}\label{Nuhat}
  N(t)-\int_{\del \Omega}d\theta_0\int_{[0,\eps]}\ds\;(1- \hat{u}(s,\theta_0,t))|J(s,\theta_0)|=\left\{\begin{array}{ll}
       O(t^\frac{5}{2}),  & 0<\lambda<\infty\\  
         O(t^2), &   \lambda=\infty       \end{array}
 .\right.
\end{equation}
If all principal curvatures of $\del \Omega$ are constant, then
$$N(t)=\int_{\del \Omega}d
\theta_0\int_{[0,\eps]}\ds\;
(1-\hat{u}(s,\theta_0,t))|J(s,\theta_0)|+O(e^{-\frac{1}{t^{\delta}}}).$$
 
Moreover, if $\hat{u}^{hom}$ is the solution of the homogeneous
constant coefficients problem
\begin{eqnarray}
&&\del_t \hat{u}-D_\pm\frac{\del^2}{\del s^2}\hat{u}=0, \quad -\eps<s<\eps, \;\theta\equiv\theta_0,\label{homp01}\\
&&\hat{u}|_{t=0}=\mathds{1}_{0<s<\eps}(s),\nonumber\\
&& D_- \frac{\del \hat{u}_{-}}{\del s}|_{s=-0}=\lambda(\theta_0)(\hat{u}_--\hat{u}_+)|_{s=0}, \\
&& D_+\frac{\del \hat{u}_+}{\del s}|_{s=+0}=D_-\frac{\del \hat{u}_-}{\del s}|_{s=-0},\label{homp0e}
\end{eqnarray} 
then
\begin{equation}
  N(t)-\int_{\del \Omega}d\theta_0\int_{[0,\eps]}\ds\; 
(1-\hat{u}^{hom}(s,\theta_0,t))|J(s,\theta_0)|=\left\{\begin{array}{ll}
       O(t^\frac{3}{2}),  & 0<\lambda<\infty\\ 
          O(t), &   \lambda=\infty                                                                                                                       
\end{array}
 \right.\label{NOt2h}
\end{equation}
\end{theorem}

{F}rom \ref{SecGreen} we get
\begin{equation*}
\begin{split}
G^{\theta_0}(s_1,\theta_1,s_2,\theta_2,t) &=
     \mathds{1}_{\{s_1>0,s_2>0\}} G^{\theta_0}_{++}(s_1,\theta_1,s_2,\theta_2,t) \\
& + \mathds{1}_{\{s_1<0,s_2>0\}}G^{\theta_0}_{-+}(s_1,\theta_1,s_2,\theta_2,t).  \\
\end{split}
\end{equation*}
Due to Eq.~(\ref{nnloc}), we need to know only $G^{\theta_0}_{++}$
\begin{eqnarray*}
 &&G^{\theta_0}_{++}(s_1,\theta_1,s_2,\theta_2,t)=\left(h_+^{\theta_0}(s_1,s_2,t)-f^{\theta_0}_+(s_1,s_2,t)\right)
K(\theta_1,\theta_2,D_+t),
\end{eqnarray*}
with notations
\begin{eqnarray}
&&h_+^{\theta_0}(s_1,s_2,t)= \frac{1}{\sqrt{4\pi D_+t}}\left(\exp\left(-\frac{(s_1-s_2)^2}{4D_+t}\right) \right.\nonumber\\
&&\left.
+a(\lambda,0,\theta_0) \exp\left(-\frac{(s_1+s_2)^2}{4D_+t}\right) \right),\label{h}\\
&&f_+^{\theta_0}(s_1,s_2,t)=b(\lambda,0,\theta_0)\frac{\lambda(\theta_0)}{D_+}
\exp\left(\frac{\lambda(\theta_0) \alpha}{\sqrt{D_+}}(s_1+s_2)+ 
\lambda(\theta_0)^2 \alpha^2 t\right)\nonumber\\
&&\cdot \operatorname{Erfc}\left(\frac{s_1+s_2}{2\sqrt{D_+t}}+\lambda(\theta_0) \alpha\sqrt{ t} \right),\label{fp}
\end{eqnarray}
where 
\begin{eqnarray*}
       &&a(\lambda,0,\theta_0)=\left\{\begin{array}{ll}
                                       1,&\lambda(\theta_0)<\infty ,\\
                                       \frac{\sqrt{D_+}-\sqrt{D_-}}{\sqrt{D_+}+\sqrt{D_-} },&\lambda(\theta_0)=\infty ,
                                      \end{array}
 \right.\\
 &&b(\lambda,0,\theta_0)=\left\{\begin{array}{ll}
                                       1,& \quad \lambda(\theta_0)<\infty,\\
                                       0,& \quad \lambda(\theta_0)=\infty,
                                      \end{array}
 \right.
\end{eqnarray*}
and $K(\theta_1,\theta_2,D_\pm t)$ is the heat kernel in $\R^{n-1}$:
\begin{equation}
 \label{Keq}
 K(\theta_1,\theta_2,D_\pm t)= \frac{1}{(4\pi D_\pm t)^\frac{n-1}{2}}
\exp\left(-\frac{|\theta_1-\theta_2|^2}{4D_\pm t}\right).
\end{equation}
Since 
\begin{eqnarray}
&&N(t)=\int_{\Omega_{\eps}} (1-u_\eps(s,\theta,t))|J(s,\theta)|\ds d\theta+O(e^{-\frac{1}{t^{\delta}}})\nonumber\\
&&=\operatorname{Vol}(\Omega_{\eps})- \int_{\Omega_\eps} \int_{\Omega_\eps} G(s,\theta,s_1,\theta_1,t)
 |J(s,\theta)|\ds d\theta\ds_1 d\theta_1+O(e^{-\frac{1}{t^{\delta}}})\label{NNiloc},
\end{eqnarray}
in what follows we use $P(t)$ for the notation of the principal part
of $N(t)$:
\begin{equation}\label{Piloc}
  P(t)=\int_{\Omega_\eps} \int_{\Omega_\eps} G(s,\theta,s_1,\theta_1,t)|J(s,\theta)|\ds d\theta\ds_1 d\theta_1.
\end{equation}

To prove Theorem~\ref{TH1} we need the following Lemma:

\begin{lemma}\label{LemPrelim}
The principal part $P(t)$ of the heat content for the solution of the
system~(\ref{eq47})--(\ref{endsys}), defined in Eq.~(\ref{Piloc}), is
given by
\begin{eqnarray}
     &&P(t)=\int_{\Omega_\eps}\int_{\Omega_\eps}G^{\theta_0}_{++}(s,\theta_0,s_1,\theta_1,t)
|J(s,\theta_0)|\ds d\theta_0 d\theta_1 d s_1+G^{\theta_0}_{++}\sharp(\mathcal{R}_s^{\theta_0}\nonumber\\
&&+\mathcal{R}_{\theta}^{\theta_0})G^{\theta_0}_{++}+ 
G^{\theta_0}_{++}\sharp(\mathcal{R}_s^{\theta_0}+\mathcal{R}_{\theta}^{\theta_0})G^{\theta_0}_{++}
\sharp(\mathcal{R}_s^{\theta_0}+\mathcal{R}_{\theta}^{\theta_0})u_++O(e^{-\frac{1}{t^\delta}}),\label{CauchyG}
\end{eqnarray}
with notation
\begin{eqnarray}
 &&G^{\theta_0}_{++}\sharp(\mathcal{R}_s^{\theta_0}+\mathcal{R}_{\theta}^{\theta_0})G^{\theta_0}_{++}=
\int_0^t\dtau\int_{\Omega_\eps}\ds d\theta_0 |J(s,\theta_0)|   \int_{\Omega_\eps}  
G^{\theta_0}_{++}(s,\theta_0,s_1,\theta_1,t-\tau)\cdot\nonumber\\
 &&\cdot\int_{\Omega_\eps}  (\mathcal{R}_{s_1}^{\theta_0}(s_1,\theta_1)
+\mathcal{R}_{\theta_1}^{\theta_0}(s_1,\theta_1))  G^{\theta_0}_{++}(s_1,\theta_1,s_2,\theta_2,\tau)  
d\theta_2 d s_2 d\theta_1 d s_1.\label{natation}
\end{eqnarray}
Moreover, the following identities hold
\begin{eqnarray}
 &&M(t)=\int_{\Omega_\eps}\int_{\Omega_\eps}G^{\theta_0}_{++}(s,\theta_0,s_1,\theta_1,t)
|J(s,\theta_0)|\ds d\theta_0 d\theta_1 d s_1\nonumber\\
&&=\int_{\del \Omega}d\theta_0 \int_{[0,\eps]^2}d s_1 \ds \; (h^{\theta_0}_+(s,s_1,t)
-f^{\theta_0}_+(s,s_1,t))|J(s,\theta_0)|,\label{M1}\\
&&G^{\theta_0}_{++}\sharp\mathcal{R}_{\theta}^{\theta_0}G^{\theta_0}_{++} = 
G^{\theta_0}_{++}\sharp\mathcal{R}_{\theta}^{\theta_0}G^{\theta_0}_{++}
\sharp\mathcal{R}_{\theta}^{\theta_0}G^{\theta_0}_{++}=\ldots=0,\label{GRT}\\
 && \int_{\del \Omega} d\theta_1 K(\theta_0,\theta_1,D_+t)\nonumber\\
 &&=\int_{ \del \Omega} d\theta_1    
\int_{\del \Omega} d\theta_2 K(\theta_0,\theta_1,
D_+(t-\tau))K(\theta_1,\theta_2,D_+\tau)\nonumber\\
 &&=\ldots =\mathds{1}_{\del \Omega}(\theta_0). \label{KK}
\end{eqnarray}
\end{lemma}

\textbf{Proof.}
Formula~(\ref{CauchyG}) is the direct corollary of the Duhamel formula
(see~(\ref{eqCauchy}) and~(\ref{uCauchyexpansion})).

Let us start to prove~(\ref{M1}).

Indeed, we find that
\begin{equation*}
\begin{split}
 M(t) & =\int_{\Omega_\eps}\int_{\Omega_\eps}(h^{\theta_0}_+(s_1,s_2,t)-f^{\theta_0}_+(s,s_1,t))
K(\theta_1,\theta_2,D_+(\theta_0)t)|J(s,\theta_0)|d s_1\ds  d\theta_1 d\theta_0  \\
 & =\int_{\R^{n-1}}\int_{\R^{n-1}} d\theta_0 d\theta_1 \frac{1}{(4\pi D_+t)^\frac{n-1}{2}}
\exp\left(-\frac{|\theta_0-\theta_1|^2}{4D_+t}\right)
 \mathds{1}_{\del \Omega}(\theta_0)\mathds{1}_{\del \Omega}(\theta_1)  \Phi(\theta_0,t), \\
\end{split}
\end{equation*}
where 
\begin{eqnarray}\label{EqPhi}
 &&\Phi(\theta_0,t)=\int_{[0,\eps]^2}dsd s_1(h^{\theta_0}_+(s,s_1,t)-f^{\theta_0}_+(s,s_1,t))|J(s,\theta_0)|.
\end{eqnarray}
With the change of variables $\theta^1\mapsto
v=\frac{\theta_0-\theta_1}{\sqrt{4D_+t}}$, $M(t)$ becomes
\begin{eqnarray*}
 &&M(t)=\int_{\R^{n-1}}\int_{\R^{n-1}} \frac{e^{- |v|^2}}{\pi^\frac{n-1}{2}}
\mathds{1}_{\del \Omega}(\theta_0)\mathds{1}_{\del \Omega
+\sqrt{4D_+t} v}(\theta_0)\Phi(\theta_0,t)  d  v d\theta_0 .
\end{eqnarray*}

By our construction, 
\begin{equation*}
\theta_0\in\del \Omega \quad \hbox{ and } \quad\theta_1=
\theta_0-\sqrt{4D_+t} v\in \del \Omega,
\end{equation*}
that implies 
\begin{equation*}
\mathds{1}_{\del \Omega}(\theta_0) -\mathds{1}_{\del \Omega}(\theta_0)
\mathds{1}_{\del \Omega+\sqrt{4D_+t}v}(\theta_0)\equiv 0.
\end{equation*}
It can be interpreted in the following way: if we take a point on the
boundary and move it along the boundary, we obtain another point which
is still a boundary point.

Consequently, we find~(\ref{M1})
\begin{eqnarray*}
 &&M(t)=\int_{\R^{n-1}}\int_{\R^{n-1}} \frac{e^{- v^2}}{\pi^\frac{n-1}{2}}
\mathds{1}_{\del \Omega}(\theta_0)\Phi(\theta_0,t)  d  v d\theta_0 =\int_{\del \Omega}\Phi(\theta_0,t)d\theta_0,
\end{eqnarray*}
which also implies the first part of~(\ref{KK}):
\begin{equation*}
\int_{\del \Omega} d\theta_1K(\theta_0,\theta_1,D_+t)=
\mathds{1}_{\del \Omega}(\theta_0).
\end{equation*}

Let us now prove that in the computation of $P(t)$ all terms
containing the derivatives over the transversal variable
$\theta$ vanish.

For all terms in~(\ref{eq47}) containing a derivative over
$\theta_1$, we calculate (see~(\ref{Restinteta}))
\begin{eqnarray*}
&&\mathcal{R}_{\theta_1}K(\theta_1,\theta_2,D_+t)\\
&&=\sum_{i=1}^{n-1}  \frac{D_+s_1k_i(\theta_1)}{1-s_1k_i(\theta_1)}
\left( 1+\frac{1}{1-s_1k_i(\theta_1)}\right)\\
&&\cdot\frac{1}{2D_+t}\left(\frac{(\theta_1^i-\theta_2^i)^2}{2D_+t}-1 \right)
K\left(\theta_1,\theta_2,D_+t\right)\\
&&-\frac{1}{|J(s_1,\theta_1)|}\sum_{i=1}^{n-1}\frac{\del}{\del \theta_1^i}\left(\frac{D_+
|J(s_1,\theta_1)|}{(1-s_1k_i(\theta_1))^2} \right) \frac{(\theta_1^i-\theta_2^i)}{2D_+t}
K(\theta_1,\theta_2,D_+t).
\end{eqnarray*}

Let us prove Eq.~(\ref{GRT}), noting that
\begin{eqnarray*}
 &&G^{\theta_0}_{++}\sharp\mathcal{R}_{\theta_1}G^{\theta_0}_{++}=\int_0^t\dtau\int_{\Omega_\eps}\ds d\theta_0 |J(s,\theta_0)| 
 \int_{\Omega_\eps}  G^{\theta_0}_{++}(s,\theta_0,s_1,\theta_1,t-\tau)\cdot\\
 &&\cdot\int_{\Omega_\eps} \mathcal{R}_{\theta_1}  G^{\theta_0}_{++}(s_1,\theta_1,s_2,\theta_2,\tau) 
 d\theta_2 d s_2 d\theta_1 d s_1.
\end{eqnarray*}
We can schematically rewrite
$G^{\theta_0}_{++}\sharp\mathcal{R}_{\theta_1}G^{\theta_0}_{++}$ in
the following form:
\begin{eqnarray*}
 &&G^{\theta_0}_{++}\sharp\mathcal{R}_{\theta_1}G^{\theta_0}_{++}=\\
 &&=\int_0^t\dtau 
\int_{\R^{n-1}}d\theta_0 \int_{\R^{n-1}}d\theta_1\int_{\R^{n-1}}d\theta_2
K(\theta_0,\theta_1,D_+(t-\tau))\mathcal{R}_{\theta_1}K(\theta_1,\theta_2,D_+\tau)\\
&& \cdot \mathds{1}_{\del \Omega}(\theta_0)\mathds{1}_{\del \Omega}(\theta_1)
\mathds{1}_{\del \Omega}(\theta_2)
\int_{[0,\eps]^3}\ds \ds_1\ds_2 \phi(s,s_1,s_2,t,\tau,\theta_0).
\end{eqnarray*}
With the change of variables involving $\theta_0$:
\begin{eqnarray}
&& \tilde{\theta}_1=\frac{\theta_0-\theta_0^1}{2\sqrt{D_+(t-\tau)}}, \quad\theta_1=
\theta_0-2\sqrt{D_+(t-\tau)} \tilde{\theta}_1,\label{theta1}\\ 
&& \tilde{\theta}_2=\frac{\theta_0^1-\theta_0^2}{2\sqrt{D_+\tau}}, \quad\theta_2=
\theta_1-2\sqrt{D_+\tau} \tilde{\theta}_2,\hbox{ and so }\nonumber\\
&&\theta_2=\theta_0-2\sqrt{D_+(t-\tau)} \tilde{\theta}_1-2\sqrt{D_+\tau}
 \tilde{\theta}_2,\label{theta2}
\end{eqnarray}
and since for all $\theta_0\in \del \Omega$
\begin{eqnarray*}
&&\mathds{1}_{\del \Omega}(\theta_0)-
\mathds{1}_{\del \Omega}(\theta_0)\mathds{1}_{\del \Omega+2\sqrt{D_+(t-\tau)}
 \tilde{\theta}_1}(\theta_0)\cdot\\
&&\cdot\mathds{1}_{\del \Omega+2\sqrt{D_+(t-\tau)} \tilde{\theta}_1+2\sqrt{D_+\tau}
 \tilde{\theta}_2}(\theta_0)=0,
\end{eqnarray*}
we obtain the separation of variables on $ \tilde{\theta}_2$ from
$(\theta_0,s_1, \tilde{\theta}_1)$:
\begin{eqnarray*}
 &&G^{\theta_0}_{++}\sharp\mathcal{R}_{\theta_1}G^{\theta_0}_{++}=\int_0^t\dtau 
\int_{\del \Omega} d\theta_0  \int_{[0,\eps]^3}\ds \ds_1\ds_2\phi(s,s_1,s_2,t,\tau,\theta_0)\\
&&\cdot  \prod_{i=1}^{n-1}\left[
\int_{\R} d \tilde{\theta}_1^i\frac{e^{-(\tilde{\theta}_1^i)^2}}{\sqrt{\pi } }\int_{\R} d \tilde{\theta}_2^i  
\frac{e^{-(\tilde{\theta}_2^i)^2}}{\sqrt{\pi} }\left(C^i_1(2(\tilde{\theta}_2^i)^2-1)-C^i_2\tilde{\theta}_2^i  \right)\right],
\end{eqnarray*}
where $C^i_1$ and $C^i_2$ are the functions of $s_1,\theta_0,
\tilde{\theta}_1$, but not of $ \tilde{\theta}_2$, and consequently
\begin{eqnarray*}
 &&G^{\theta_0}_{++}\sharp\mathcal{R}_{\theta}G^{\theta_0}_{++}=0.
\end{eqnarray*}
By the same reason we have Eq.~(\ref{GRT}). Changing variables
$\theta_i$ to $ \tilde{\theta}_i$
from~(\ref{theta1})--(\ref{theta2}), we also obtain the last part
of~(\ref{KK}).
$\Box$
Let us know prove Theorem~\ref{TH1}.

\textbf{Proof.}
To find Eq.~(\ref{Nuhat}), we study Eq.~(\ref{CauchyG}) using
proved relations (\ref{M1})--(\ref{KK}). For instance, we have
\begin{eqnarray*}
 &&G^{\theta_0}_{++}\sharp(\mathcal{R}_s^{\theta_0}+\mathcal{R}_{\theta}^{\theta_0})G^{\theta_0}_{++}
=G^{\theta_0}_{++}\sharp\mathcal{R}_s^{\theta_0}G^{\theta_0}_{++}\\
&&=\int_0^t\dtau\int_{\Omega_\eps}\ds d\theta_0 |J(s,\theta_0)|   \int_{\Omega_\eps}  
G^{\theta_0}_{++}(s,\theta_0,s_1,\theta_1,t-\tau)\cdot\nonumber\\
 &&\cdot\int_{\Omega_\eps}  \mathcal{R}_{s_1}^{\theta_0}(s_1,\theta_1)  
G^{\theta_0}_{++}(s_1,\theta_1,s_2,\theta_2,\tau)  d\theta_2 d s_2 d\theta_1 d s_1.
\end{eqnarray*}
As
$G^{\theta_0}_{++}(s,\theta_0,s_1,\theta_1,t)=(h^{\theta_0}_+(s,s_1,t)
-f^{\theta_0}_+(s,s_1,t))K(\theta_0,\theta_1,D_+(\theta_0)t)$, we have
\begin{eqnarray*}
 &&G^{\theta_0}_{++}\sharp\mathcal{R}_s^{\theta_0}G^{\theta_0}_{++}=(h^{\theta_0}_+
-f^{\theta_0}_+)K\sharp\mathcal{R}_s^{\theta_0}(h^{\theta_0}_{+}-f^{\theta_0}_{+})K.
\end{eqnarray*}

Now we perform the change of variables~(\ref{theta1})--(\ref{theta2}).
Since locally $k_i(\theta)\in C^1$, then for all $\theta_0\in \del \Omega$,  for $t\to +0$ we can
develop
\begin{equation*}
k_i(\theta_0-2\sqrt{D_+(t-\tau)} \tilde{\theta}_1)=k_i(\theta_0)-\nabla
k_i(\theta_0)2\sqrt{D_+(t-\tau)} \tilde{\theta}_1+O(t-\tau).
\end{equation*}
Consequently, by definition of $\mathcal{R}_{s_1}(s_1,\theta_1)$
in~(\ref{Restins}), which is a composition of the operator of the
first derivative by $s_1$ and of a multiplication by a function of the
class $C^1$ on $\theta_1$ (locally, in the sense of local variables),
we also have for $t\to +0$
\begin{eqnarray*}
 && \mathcal{R}_{s_1}(s_1,\theta_1)=\mathcal{R}_{s_1}(s_1,\theta)[1+O(t-\tau)] -\nabla_{\theta}  
\mathcal{R}_{s_1}(s_1,\theta)2\sqrt{D_+(t-\tau)} \tilde{\theta}_1.
\end{eqnarray*}
As $2\sqrt{D_+(t-\tau)}\int_\R d \tilde{\theta}_i^1
e^{-(\tilde{\theta}_i^1)^2}\tilde{\theta}_i^1=0$, we obtain
\begin{eqnarray*}
 &&G^{\theta_0}_{++}\sharp\mathcal{R}_s^{\theta_0}G^{\theta_0}_{++}=
\int_0^t\dtau\int_{\del\Omega} d\theta_0 \int_{[0,\eps]} \ds |J(s,\theta_0)|   \int_{[0,\eps]} \ds_1 \int_{[0,\eps]}\ds_2\\
&&\cdot (h^{\theta_0}_+-f^{\theta_0}_+)(s,s_1,t-\tau)
\mathcal{R}_{s_1}(s_1,\theta_0)\left[1+O(t-\tau)\right](h^{\theta_0}_{+}-f^{\theta_0}_{+})(s_1,s_2,\tau)\\
&&\cdot
\prod_{i=1}^{n-1}\left[
\int_{\R} d \tilde{\theta}^1_i\frac{e^{-(\tilde{\theta}^1_i)^2}}{\sqrt{\pi } }\int_{\R} d \tilde{\theta}_i^2  
\frac{e^{-(\tilde{\theta}^2_i)^2}}{\sqrt{\pi} }\right]
=\int_0^t\dtau\int_{\del\Omega} d\theta_0 \int_{[0,\eps]} \ds |J(s,\theta_0)|   \int_{[0,\eps]} \ds_1 \\
&&\cdot\int_{[0,\eps]}\ds_2 (h^{\theta_0}_+-f^{\theta_0}_+)(s,s_1,t-\tau)\mathcal{R}_{s_1}(s_1,\theta_0)
\left[1+O(t-\tau)\right](h^{\theta_0}_{+}-f^{\theta_0}_{+})(s_1,s_2,\tau),
\end{eqnarray*}
from which it follows
\begin{eqnarray}
  &&P(t)= \int_{\del \Omega}d\theta_0\int_{[0,\eps]}\ds \int_{[0,\eps]}\ds_1 
(h^{\theta_0}_+(s,s_1,t)-f_+^{\theta_0}(s,s_1,t))|J(s,\theta_0)| \nonumber\\ 
  &&+\left[1+O(t)\right]\int_0^t\dtau\int_{\del\Omega} d\theta_0    
\int_{[0,\eps]}\ds|J(s,\theta_0)|\int_{[0,\eps]}\ds_1 (h^{\theta_0}_{+}-f^{\theta_0}_{+})(s,s_1,t-\tau)\nonumber\\
  &&
 \cdot\int_{[0,\eps]}\ds_2 \mathcal{R}_{s_1}(s_1,\theta_0) (h^{\theta_0}_{+}-f^{\theta_0}_{+})(s_1,s_2,\tau) \nonumber\nonumber\\
  &&+\left[1+O(t)\right]^2\int_{\del\Omega} d\theta_0  (h^{\theta_0}_{+}-f^{\theta_0}_{+}) 
\sharp  \mathcal{R}_{s}(s,\theta_0)(h^{\theta_0}_{+}-f^{\theta_0}_{+})\sharp  \mathcal{R}_{s}(s,\theta_0)u^{+}_\eps\label{Mm}.
\end{eqnarray}
We notice that the solution $\hat{u}(s,\theta_0,t)$ of the
one-dimensional system~(\ref{opereqCH1d})--(\ref{opereqCH1dN}) is
given by
\begin{eqnarray*}
  &&\hat{u}(s,\theta_0,t)=\int_{[0,\eps]}\ds_1 (h^{\theta_0}_+(s,s_1,t)-f_+^{\theta_0}(s,s_1,t))\nonumber\\ 
  &&+\int_0^t\dtau\int_{[0,\eps]}\ds_1 (h^{\theta_0}_{+}-f^{\theta_0}_{+})(s,s_1,t-\tau)\int_{[0,\eps]}\ds_2 
\mathcal{R}_{s_1}(s_1,\theta_0) (h^{\theta_0}_{+}-f^{\theta_0}_{+})(s_1,s_2,\tau) \nonumber\nonumber\\
  &&+  (h^{\theta_0}_{+}-f^{\theta_0}_{+}) \sharp \mathcal{R}_{s}(s,\theta_0)(h^{\theta_0}_{+}-f^{\theta_0}_{+})\sharp  
\mathcal{R}_{s}(s,\theta_0)\hat{u}.
\end{eqnarray*}

To obtain~(\ref{Nuhat}) of Theorem~\ref{TH1} from formula~(\ref{Mm}),
we estimate
\begin{equation}
\label{NNloc}
\begin{split}
  NN^2(t) &=
 O(t)\int_0^t\dtau\int_{\del\Omega} d\theta_0    \int_{[0,\eps]}\ds|J(s,\theta_0)| \\
  &
 \cdot\int_{[0,\eps]} \hspace*{-1mm} \ds_1 (h^{\theta_0}_{+}-f^{\theta_0}_{+})(s,s_1,t-\tau)\int_{[0,\eps]} \hspace*{-1mm}  \ds_2 
  \mathcal{R}_{s_1}(s_1,\theta_0)(h^{\theta_0}_{+}-f^{\theta_0}_{+})(s_1,s_2,\tau) . \\
\end{split} 
\end{equation}
In fact, from~(\ref{NOt2h}), proven in what follows, it holds
(see~(\ref{defNN1})  for the definition of
$NN^1(t)$)
\begin{equation*}
NN^2(t)=O(t)NN^1(t)=O(t)\left\{\begin{array}{ll}
       O(t^\frac{3}{2}),  & 0<\lambda<\infty\\  
         O(t), &   \lambda=\infty \end{array}
 \right. = \left\{\begin{array}{ll}
       O(t^\frac{5}{2}),  & 0<\lambda<\infty\\  
         O(t^2), &   \lambda=\infty \end{array}
 .\right.
\end{equation*}
   
To conclude, we note that if all principal curvatures $k_j(\theta)$ on
$\Omega_\eps$ are constant, then for all $\theta\in
\del\Omega$
\begin{equation*}
\mathcal{R}_{s}(s,\theta)\equiv \mathcal{R}_{s}(s,\theta_0),
\end{equation*}
and thus 
\begin{equation*}
N(t)= \int_{\del \Omega}d\theta_0\int_{[0,\eps]}\ds\; (1-\hat{u}(s,\theta_0,t))
|J(s,\theta_0)|+O(e^{-\frac{1}{t^{\delta}}}). 
\end{equation*}

To show~(\ref{NOt2h}), we need to estimate
\begin{eqnarray}
&&NN^j(t)=\sum_{l=1}^j\Gamma^{\theta_0}_{++}\sharp  \mathcal{R}_{s}\Gamma^{\theta_0}_{++}\sharp \ldots \sharp 
\mathcal{R}_{s}\Gamma^{\theta_0}_{++} \quad (l-\hbox{fold}),\label{defNN1}
 \end{eqnarray}
where
\begin{equation*}
\Gamma^{\theta_0}_{++}=(h^{\theta_0}_{+}-f^{\theta_0}_{+}).
\end{equation*}
More precisely we want to prove that for all $j\ge 1$
\begin{equation}\label{NNj}
 |NN^j(t)|\le C\left\{\begin{array}{ll}
       t^{\frac{1+j}{2}}\mu(\del \Omega,\sqrt{4D_+t}),  & \quad 0<\lambda<\infty\\ 
        t^{\frac{j}{2}}\mu(\del \Omega,\sqrt{4D_+t}), &  \quad \lambda=\infty   \end{array}
 \right..
\end{equation}

Due to Lemma~\ref{LemPrelim}, we start with (see~(\ref{Mm}))
\begin{eqnarray*}
 &&NN^1(t)=\int_0^t\dtau\int_{\del\Omega} d\theta_0    \int_{[0,\eps]}\ds|J(s,\theta_0)|\nonumber\\
 &&
\cdot\int_{[0,\eps]}\ds_1 (h^{\theta_0}_{+}-f^{\theta_0}_{+})(s,s_1,t-\tau)\int_{[0,\eps]}\ds_2 
\mathcal{R}_{s_1}(s_1,\theta_0) (h^{\theta_0}_{+}-f^{\theta_0}_{+})(s_1,s_2,\tau). 
\end{eqnarray*}
Therefore, we have to estimate four terms:
\begin{eqnarray*}
 &&NN^1(t)=\sum_{j=1}^4 MM_j(t),
\end{eqnarray*}
where 
\begin{eqnarray*}
 &&MM_1(t)=\int_0^t\dtau\int_{\del\Omega} d\theta_0    \int_{[0,\eps]}\ds|J(s,\theta_0)|\nonumber\\
 &&
\cdot\int_{[0,\eps]}\ds_1 h^{\theta_0}_{+}(s,s_1,t-\tau)\int_{[0,\eps]}\ds_2 \mathcal{R}_{s_1}(s_1,\theta_0) h^{\theta_0}_{+}(s_1,s_2,\tau),\\
&&MM_2(t)=-\int_0^t\dtau\int_{\del\Omega} d\theta_0    \int_{[0,\eps]}\ds|J(s,\theta_0)|\nonumber\\
 &&
\cdot\int_{[0,\eps]}\ds_1 f^{\theta_0}_{+}(s,s_1,t-\tau)\int_{[0,\eps]}\ds_2 \mathcal{R}_{s_1}(s_1,\theta_0) h^{\theta_0}_{+}(s_1,s_2,\tau),\\
&&MM_3(t)=-\int_0^t\dtau\int_{\del\Omega} d\theta_0    \int_{[0,\eps]}\ds|J(s,\theta_0)|\nonumber\\
 &&
\cdot\int_{[0,\eps]}\ds_1 h^{\theta_0}_{+}(s,s_1,t-\tau)\int_{[0,\eps]}\ds_2 \mathcal{R}_{s_1}(s_1,\theta_0) f^{\theta_0}_{+}(s_1,s_2,\tau),\\
&&MM_4(t)=\int_0^t\dtau\int_{\del\Omega} d\theta_0    \int_{[0,\eps]}\ds|J(s,\theta_0)|\nonumber\\
 &&
\cdot\int_{[0,\eps]}\ds_1 f^{\theta_0}_{+}(s,s_1,t-\tau)\int_{[0,\eps]}\ds_2 \mathcal{R}_{s_1}(s_1,\theta_0) f^{\theta_0}_{+}(s_1,s_2,\tau).
\end{eqnarray*}
We aim to approximate
$\mathcal{R}_{s_1}(s_1,\theta_0)=R(s_1,\theta_0)\del_{s_1}$ from
Eq.~(\ref{Restins}) near the point~$(s,\theta_0)$.  For $t\to+0$ and
$0<s_1<\eps=O(\sqrt{t})$, we find that
\begin{equation*}
\frac{1}{1-s_1k_i(\theta_0)}=1+s_1k_i(\theta_0)+O(s_1^2),
\end{equation*}
which gives
\begin{eqnarray*}
 &&R(s_1,\theta_0)=-D_+
\left(\sum_{i=1}^{n-1}k_i(\theta_0)+s_1\sum_{i=1}^{n-1}k_i^2(\theta_0) +O(s_1^2) \right).
\end{eqnarray*}
Introducing the notations
\begin{eqnarray*}
 &&C_\pm=s_1\pm s_2,\quad I_{s_1\pm s_2}(\tau)= \exp\left(-\frac{(s_1\pm s_2)^2}{4D_+\tau}\right), \quad 
c= \frac{1}{8\pi D_+^2(0,\theta_0)\sqrt{(t-\tau)}\tau^\frac{3}{2}},
\end{eqnarray*}
we find
\begin{eqnarray*}
&&h^{\theta_0}_+(s,s_1,t-\tau)\mathcal{R}_{s_1}(s_1,\theta_0)h^{\theta_0}_+(s_1,s_2,\tau)=
-cR(s_1,\theta_0)\left(C_-[I_{s-s_1}(t-\tau)I_{s_1-s_2}(\tau)\right.\\
&&\left.+a(\lambda,0,\theta_0) I_{s+s_1}(t-\tau)I_{s_1-s_2}(\tau)]\right.\\
&&\left. +a(\lambda,0,\theta_0)C_+[I_{s-s_1}(t-\tau)I_{s_1+s_2}(\tau)+a(\lambda,0,\theta_0) I_{s+s_1}(t-\tau)I_{s_1+s_2}(\tau)]\right).
\end{eqnarray*}

We now change $s_1$ to $z_1$ and $s_2$ to $z_2$ by the following
change of variables:
\begin{itemize}
\item for $P_{s_1\mp s_2}(\tau)$: $z_2=\frac{s_1\mp s_2}{2\sqrt{D_+\tau}}$ and $s_2=\pm s_1\mp 2\sqrt{D_+\tau}z_2$,
\item for $P_{s\mp s_1}(t-\tau)$: $z_1=\frac{s\mp s_1}{2\sqrt{D_+(t-\tau)}}$ and $s_1=\pm s\mp 2\sqrt{D_+(t-\tau)}z_1$.
\end{itemize}
Let us notice that $t$ is a constant parameter and, as $\tau$ takes
its values between $0$ and $t$, hence, $z_1$ and $z_2$ are in $\R^+$
or $\R$.  But at the same time $2\sqrt{D_+(t-\tau)}z_1=s\pm s_1$ and
$2\sqrt{D_+\tau}z_2= s_1\pm s_2$ are bounded to the interval
$[-\eps,2\eps]$ and hence are of the order of $O(\sqrt{t})$.  In what
follows, we suppose that $\tau$ and $t-\tau$ have the same order of
smallness as $t$:
\begin{equation*}
O(t)=O(\tau)=O(t-\tau).
\end{equation*}

Therefore, for $0<s_1=\pm
s\mp2\sqrt{D_+(t-\tau)}z_1<\eps$ we have
\begin{eqnarray*}
&&\mathcal{R}_{s_1}(s_1,\theta_0) =
[\phi(\theta_0)\mp\psi(s,z_1,\theta_0)]\frac{\del}{\del s_1},
\end{eqnarray*}
where
\begin{eqnarray*}
 &&\phi(\theta_0)=-D_+
\sum_{i=1}^{n-1}k_i(\theta_0),\\
&&\psi(s,z_1,\theta_0)= (s-2\sqrt{D_+(t-\tau)}z_1)D_+\sum_{i=1}^{n-1}k_i^2(\theta_0)+O\left(t\right).
\end{eqnarray*}

If we develop $R$ in the neighborhood of
$(0,\theta_0)$, we find
\begin{eqnarray}
&&R(s_1,\theta_0)=R(0,\theta_0)+O(\sqrt{t})= \phi(\theta_0)+O(\sqrt{t}).\label{Rest0S1}
\end{eqnarray}

For $\lambda=\infty$ on $\del \Omega$, we simply have
\begin{equation*}
MM_2(t)=MM_3(t)=MM_4(t)=0,
\end{equation*}
and
\begin{eqnarray*}
 &&|NN^1(t)|=|MM_1(t)|=|h^{\theta_0}_+\sharp\mathcal{R}_{s}h^{\theta_0}_+|\nonumber\\
&&\le C\left|\int_0^t\dtau\frac{1}{\sqrt{\tau}}\int_{\del\Omega} d\theta_0 
\int_{0}^\epsilon \ds |J(s,\theta_0)|\right|\le C\sqrt{t}\mu(\del \Omega,\sqrt{4D_+t}).\label{etoile}
\end{eqnarray*}
By iteration of the proof, we show that $|NN^j(t)|\le
Ct^{\frac{j}{2}}\mu(\del \Omega,\sqrt{4D_+t})$ for $j\ge 1$.

Now, for $\lambda<\infty$,
\begin{eqnarray*}
&&MM_1(t)=h^{\theta_0}_+\sharp\mathcal{R}_{s}h^{\theta_0}_+=-\int_0^t\dtau\frac{1}{\sqrt{\tau}}\int_{\del\Omega} d\theta_0 
\frac{\pi}{\sqrt{D_+}}\int_{\R} \ds |J(s,\theta_0)|  \int_{\R} \dz_1 e^{-z^2_1}\\
&&\cdot\left[\int_{\R}\dz_2 \phi(\theta_0) z_2e^{-z^2_2}\right] 
\sum_{i=1}^4\chi_i(s,z_1,z_2)-\int_0^t\dtau\frac{1}{\sqrt{\tau}}\int_{\del\Omega} d\theta_0 \frac{\pi}{\sqrt{D_+}}
\int_{\R} \ds |J(s,\theta_0)|   \\
&& \cdot \int_{\R} \dz_1 e^{-z^2_1} \cdot\left[\int_{\R}\dz_2 \psi(s,z_1,\theta_0) z_2e^{-z^2_2}\right] \cdot 
\left(\sum_{i=1}^2\chi_i(s,z_1,z_2)-\sum_{i=3}^4\chi_i(s,z_1,z_2) \right).
\end{eqnarray*}
Here for $v_+(z_1)= 2\sqrt{D_+(t-\tau)}z_1\mathds{1}_{\R^+}(z_1)$,
$v(z_1)= 2\sqrt{D_+(t-\tau)}z_1$,
$w^+(z_2)=2\sqrt{D_+\tau}z_2\mathds{1}_{\R^+}(z_2)$ and
$w(z_2)=2\sqrt{D_+\tau}z_2$
\begin{eqnarray*}
&&\chi_1(s,z_1,z_2)= \mathds{1}_{[0,\eps]}(s)\mathds{1}_{[0,\eps]}(s-v(z_1))\mathds{1}_{[0,\eps]}(s-v(z_1)-w(z_2)),\\
 &&\chi_2(s,z_1,z_2)=\mathds{1}_{[0,\eps]}(s)\mathds{1}_{[0,\eps]}(s-v(z_1))\mathds{1}_{[0,\eps]}(-s+v(z_1)+w^+(z_2)),\\
&&\chi_3(s,z_1,z_2)=\mathds{1}_{[0,\eps]}(s)\mathds{1}_{[0,\eps]}(-s+v_+(z_1))\mathds{1}_{[0,\eps]}(-s+v_+(z_1)-w(z_2)),\\                                            
&&\chi_4(s,z_1,z_2)=\mathds{1}_{[0,\eps]}(s)\mathds{1}_{[0,\eps]}(-s+v_+(z_1))\mathds{1}_{[0,\eps]}(s-v_+(z_1)+w^+(z_2)).
\end{eqnarray*}

Considering two formulas:
\begin{equation*}
\eta=\mathds{1}_{[0,\eps]}(s)-\mathds{1}_{[0,\eps]}(s-v)\mathds{1}_{[0,\eps]}(s), \quad \hbox{ and }
\zeta=\mathds{1}_{[0,\eps]}(s)\mathds{1}_{[0,\eps]}(-s+v_+),
\end{equation*}
we find that
\begin{eqnarray}
 &&\eta\ne 0 \quad \Longleftrightarrow \quad \left\{\begin{array}{ll}
                                                             0<v<\eps&0<s<v\\
                                                             -\eps<v<0&\eps+v<s<\eps
                                                            \end{array}
 \right.\label{eta}\\
 &&\zeta\ne 0 \quad \Longleftrightarrow \quad \left\{\begin{array}{ll}
                                                             0<v<\eps&0<s<v\\
                                                             \eps<v<2\eps&v-\eps<s<\eps
                                                            \end{array} \right.\label{zeta}
\end{eqnarray}
It means that for $0<v=v_+<\eps$, it holds 
\begin{equation*}
\eta(s)= \zeta(s)=\mathds{1}_{[0,v_+]}(s)
\end{equation*}
and for $-\eps<v<0$ and $\eps<v_+=v+2\eps<2\eps$ it holds
\begin{equation*}
\eta(s)= \zeta(s)=\mathds{1}_{[v_+-\eps,\eps]}(s)=\mathds{1}_{[\eps+v,\eps]}(s).
\end{equation*}
Consequently, we found the formula
\begin{equation}\label{relation}
 \mathds{1}_{[0,\eps]}(s)-\mathds{1}_{[0,\eps]}(s-v)\mathds{1}_{[0,\eps]}(s)=\mathds{1}_{[0,\eps]}(s)\mathds{1}_{[0,\eps]}(-s+v_+),
\end{equation}
from which it follows
\begin{eqnarray*}
&&\chi_1= \mathds{1}_{[0,\eps]}(s)\mathds{1}_{[0,\eps]}(s-v)-\chi_2,\quad \chi_4= \mathds{1}_{[0,\eps]}(s)\mathds{1}_{[0,\eps]}(-s+v_+)-\chi_3.
\end{eqnarray*}

Therefore, we have
\begin{equation*}
 \sum_{i=1}^4\chi_i(s,z_1,z_2)=\mathds{1}_{[0,\eps]}(s) \quad \hbox{and } \chi_1+\chi_2-\chi_3-\chi_4
=\mathds{1}_{[0,\eps]}(s)-2\cdot\mathds{1}_{[0,v_+(z_1)]}(s),
\end{equation*}
which are independent of $z_2$.  Since
\begin{equation*}
\int_{\R} z_2 e^{-z_2^2}d z_2=0\quad \hbox{ and } \quad \int_{\R} (2z_2^2-1) e^{-z_2^2}d z_2=0,
\end{equation*}
we obtain exactly
\begin{eqnarray*}
&&|MM_1(t)|=|h\sharp\mathcal{R}_{s}h|=0.
\end{eqnarray*}

For $MM_2$ we find in completely analogous way
\begin{eqnarray*}
 &&MM_2(t)=\frac{2}{\sqrt{\pi}}\int_0^t\dtau\frac{\sqrt{t-\tau}}{\sqrt{\tau}}\int_{\del\Omega} d\theta_0 
\frac{\lambda(\theta_0)}{D_+}\int_{\R} \ds |J(s,\theta_0)| \\
&&\cdot \int_{\R^+} \dz_1  e^{2\lambda(\theta_0) \alpha z_1\sqrt{t-\tau}
+\lambda(\theta_0)^2\alpha^2(t-\tau)}\operatorname{Erfc}(z_1+ \lambda(\theta_0) \alpha\sqrt{t-\tau})\\
&&\cdot\int_{\R}\dz_2  \left[\phi(\theta_0) z_2e^{-z^2_2} 
\right](\chi_3(s,z_1,z_2)+\chi_4(s,z_1,z_2))\\
&&+  \frac{2}{\sqrt{\pi}}\int_0^t\dtau\frac{\sqrt{t-\tau}}{\sqrt{\tau}}\int_{\del\Omega} 
d\theta_0\frac{\lambda(\theta_0)}{D_+}\int_{\R} \ds |J(s,\theta_0)|  \int_{\R^+} \dz_1 \\
&&\cdot e^{2\lambda(\theta_0) \alpha z_1\sqrt{t-\tau}+\lambda(\theta_0)^2\alpha^2(t-\tau)}
\operatorname{Erfc}(z_1+ \lambda(\theta_0) \alpha\sqrt{t-\tau})\\
&&\cdot \int_{\R}\dz_2\left[\psi(s,z_1,\theta_0) z_2e^{-z^2_2} \right]  
(\chi_3(s,z_1,z_2)-\chi_4(s,z_1,z_2)).
\end{eqnarray*}
Since 
\begin{eqnarray*}
&&\chi_3(s,z_1,z_2)+\chi_4(s,z_1,z_2) = \mathds{1}_{[0,\eps]}(s)\mathds{1}_{[0,\eps]}(-s+v_+(z_1))=\mathds{1}_{[0,v_+(z_1)]}(s),\\
 &&\chi_3(s,z_1,z_2)-\chi_4(s,z_1,z_2)=2\chi_3(s,z_1,z_2) -\mathds{1}_{[0,\eps]}(s)\mathds{1}_{[0,\eps]}(-s+v_+(z_1)),
\end{eqnarray*}
the parts of $MM_2$, which contain the integration over $s$ on
$[0,2\sqrt{D_+(t-\tau)}z_1]$, are equal to zero.  In addition, for
$\ell=0,1$
\begin{eqnarray*}
 && \left |\int_{\R^+}\dz_1 z_1^\ell e^{2\lambda(\theta_0) \alpha z_1\sqrt{t-\tau}
+\lambda(\theta_0)^2\alpha^2(t-\tau)}\operatorname{Erfc}(z_1+ \lambda(\theta_0) 
\alpha\sqrt{t-\tau})\right|\le C.
\end{eqnarray*}
As $\psi$ is of the order $O(\sqrt{t})$ and linear on $z_1$, and
$\eps=O(\sqrt{t})$, we directly obtain
\begin{eqnarray*}
 &&|MM_2(t)|=\left|\frac{2}{\sqrt{\pi}}\int_0^t\dtau\frac{\sqrt{t-\tau}}{\sqrt{\tau}}\int_{\del\Omega} 
d\theta_0\frac{\lambda(\theta_0)}{D_+}\int_{\R} \ds |J(s,\theta_0)|  \int_{\R^+} \dz_1 \right.\\
&&\left.\cdot e^{2\lambda(\theta_0) \alpha z_1\sqrt{t-\tau}
+\lambda(\theta_0)^2\alpha^2(t-\tau)}\operatorname{Erfc}(z_1+ \lambda(\theta_0) 
\alpha\sqrt{t-\tau}) \right.\\
&&\left.\cdot\int_{\R}\dz_2\psi(s,z_1,\theta_0) z_2e^{-z^2_2} 2\chi_3(s,z_1,z_2)\right|\\
&&\le C\int_0^t\dtau \sqrt{\tau}\int_{\del\Omega} 
d\theta_0 \int_{0}^{\eps}\ds |J(s,\theta_0)|  \le Ct^\frac{3}{2}\mu(\del \Omega,\sqrt{4D_+t}).
\end{eqnarray*}
Since $\mu(\del \Omega,\sqrt{4D_+t})=C\sqrt{t}$ for a regular
boundary, then $|MM_2(t)|\le Ct^2.$

To estimate $MM_3$ we find
\begin{eqnarray*}
 \del_{s_1} f^{\theta_0}_+(s_1,s_2,\tau)&&=\frac{\lambda(\theta_0) \alpha}{\sqrt{D_+}}
f^{\theta_0}_+(s_1,s_2,\tau)-\frac{\lambda(\theta_0)}{D_+}\frac{1}{\sqrt{\pi D_+ \tau}}
\exp\left(-\frac{(s_1+s_2)^2}{4D_+\tau}\right).
\end{eqnarray*}
In our notations, using~(\ref{Rest0S1}), we have
\begin{eqnarray*}
 &&h^{\theta_0}_+(s,s_1,t-\tau)\mathcal{R}_{s_1}(s_1,\theta_0)f^{\theta_0}_+(s_1,s_2,\tau)=
\frac{P_{s-s_1}+P_{s+s_1}}{\sqrt{4\pi D_+ (t-\tau)}}\cdot\\
 &&\cdot(\phi(\theta_0)+O(\sqrt{t}))\left\{\frac{\lambda(\theta_0) 
\alpha}{\sqrt{D_+}}f^{\theta_0}_+(s_1,s_2,\tau)
-\frac{\lambda(\theta_0) }{D_+}\frac{1}{\sqrt{\pi D_+\tau}}P_{s_1+s_2}\right\}.
\end{eqnarray*}
Changing variables $s_1$ to $z_1$ and $s_2$ to $z_2$, we obtain
$\chi_2\pm\chi_4$ for the area of $s$, which gives intervals
(linearly) depending on the values of $z_1$ and $z_2$. Thus, we
majorate $s$ by $\eps$ and estimate $MM_3$:
\begin{eqnarray*}
 &&|MM_3(t)|\le C|\int_0^t\dtau\sqrt{\tau}\int_{\del\Omega} d\theta_0 \int_{0}^\eps \ds |J(s,\theta_0)| \\
&&\cdot \int_{\R} \dz_1 e^{-z_1^2}\int_{\R^+}\dz_2  (\phi(\theta_0)+O(\sqrt{t}))f(z_2,\tau)|,
\end{eqnarray*}
where
\begin{eqnarray*}
 &&f(z_2,\tau)=\frac{\lambda(\theta_0) \alpha}{\sqrt{D_+}}
\exp\left(2\lambda(\theta_0) \alpha z_2\sqrt{\tau}+\lambda(\theta_0)^2\alpha^2\tau\right)\cdot\\
 &&\operatorname{Erfc}(z_2+ \lambda(\theta_0) \alpha\sqrt{\tau})-\frac{1}{\sqrt{\pi D_+ \tau}}e^{-z^2_2}.
\end{eqnarray*}
We see that 
\begin{eqnarray*}
 &&\sqrt{\tau}\left|\int_{\R^+}\dz_2  (\phi(\theta_0)+O(\sqrt{t}))f(z_2,\tau)\right|\le C.
\end{eqnarray*}
Therefore, we have
\begin{eqnarray*}
  &&|MM_3(t)| \le Ct \mu(\del \Omega,\sqrt{4D_+t}).
\end{eqnarray*}

In the same way, since $\chi_4$ depends on $z_1$ and $z_2$ at the same
time, we have
\begin{eqnarray*}
 &&|MM_4(t)|\le C\left|\int_0^t\dtau\sqrt{\tau(t-\tau)}\int_{\del\Omega} d\theta_0\int_{0}^\eps \ds |J(s,\theta_0)| \right.\\
&&\left.\cdot \int_0^{+\infty} \dz_1 e^{2\lambda(\theta_0) \alpha z_1\sqrt{t-\tau}
+\lambda(\theta_0)^2\alpha^2(t-\tau)}\operatorname{Erfc}(z_1+ \lambda(\theta_0) \alpha\sqrt{t-\tau})\right.\\
&&\left.\cdot  \int_0^{+\infty}\dz_2 (\phi(\theta_0)+O(\sqrt{t}))f(z_2,\tau)\right|\le Ct^\frac{3}{2}\mu(\del \Omega,\sqrt{4D_+t}).
\end{eqnarray*}
By iteration of the proof, we show for $j\ge
1$ that 
\begin{equation*}
|NN^j(t)|\le Ct^{\frac{1+j}{2}}\mu(\del \Omega,\sqrt{4D_+t}).
\end{equation*}

$\Box$

\section{Relation of the heat content expansion with the interior Minkowski sausage}\label{subs2}

Let us start with a heat problem with just a discontinuous initial
condition.

\subsection{Particular case $D_+=D_-=const$}\label{SSpartic}

\begin{lemma}
Let $\Omega \subset\R^n$ be a compact connected bounded domain with a
connected boundary $\del \Omega$ of the Hausdorff dimension $d$ and
$u$ is the solution of the following problem:
\begin{eqnarray}
 && \del_t u-D \triangle u=0 \quad x\in \R^n, \; t>0,\label{pch}\\
&& u|_{t=0}=\mathds{1}_{\Omega },\label{ps}
\end{eqnarray}
Then for $t\to+0$ we have
\begin{eqnarray}
  &&N (t)=\int_{0}^2\frac{e^{-z^2}}{\sqrt{\pi}}\mu(\del\Omega,2\sqrt{D t}z)\dz + o\left(t^\frac{n-d}{2}\right)\label{resV}.
\end{eqnarray}
Moreover, it can be approximated by 
\begin{eqnarray}
  &&N (t)= \beta_{n-d}~ \mu(\del\Omega,2\sqrt{D t})  + o\left(t^\frac{n-d}{2}\right),\label{resVV} 
\end{eqnarray}
with the prefactor
\begin{equation}  \label{eq:betax}
\beta_x \equiv \int_0^{2}\frac{z^{x} e^{-z^2}}{\sqrt{\pi}} \dz = \frac{1}{2\sqrt{\pi}}~ \gamma\left(\frac{x+1}{2},4\right)
\end{equation}
is expressed through the incomplete Gamma function.
\end{lemma}

\textbf{Proof.}
Let us prove formula~(\ref{resV}).  By definition
\begin{eqnarray*}
 &&N (t)=\int_{\R^n\setminus\overline{\Omega}}\int_{\R^n}G(x,y,t)\mathds{1}_{\Omega }\dx\dy,
\end{eqnarray*}
where this time $G$ is the heat kernel in $\R^n$
\begin{equation*}
G(x,y,t)= (4D\pi t)^{-\frac{n}{2}} \exp\left(-\frac{|x-y|^2}{4Dt}\right).
\end{equation*}
Therefore, we have
\begin{eqnarray*}
 &&N (t) =\mathrm{Vol}(\Omega )-\int_{\R^n}\int_{\R^n}\frac{1}{(4\pi D t)^\frac{n}{2}}e^{-\frac{|x-y|^2}{4Dt}}
    \mathds{1}_{\Omega }(x)\mathds{1}_{\Omega }(y)\dx\dy\\
 &&=\mathrm{Vol}(\Omega )-\int_{\R^n}\frac{1}{\pi^\frac{n}{2}}e^{-|v|^2}  
    \left( \int_{\R^n}\mathds{1}_{\Omega }(x)\mathds{1}_{\Omega }(x+2\sqrt{Dt}v)\dx\right)\dv\\
&&=\int_{\R^n}\frac{1}{\pi^\frac{n}{2}}e^{-|v|^2}  \left[\int_{\Omega }\left(\mathds{1}_{\Omega }(x)-\mathds{1}_{\Omega -2\sqrt{Dt}v}(x)\right)\dx\right]\dv,
\end{eqnarray*}
where $\mathds{1}_{\Omega-2\sqrt{Dt}v }(x)=\mathds{1}_{\Omega
}(x+2\sqrt{Dt}v)$ and the notation $\Omega-2\sqrt{Dt}v$ means that
$\Omega$ is shifted by the vector $-2\sqrt{Dt}v\in \R^n$.

Let us firstly suppose that $\del \Omega$ is regular, i.e of the class
$C^3$.  We see that for all points $x\in \Omega$ for which $d(x,\del
\Omega )\ge 2\sqrt{Dt}\|v\|$, it holds $(x+2\sqrt{Dt}v)\in\Omega $.
Thus, it follows that for $\eps=2\sqrt{Dt}\|v\|$,
\begin{equation*}
\mathds{1}_{\Omega }(x)\left(\mathds{1}_{\Omega }(x)-\mathds{1}_{\Omega -2\sqrt{Dt}v}(x)\right)=0 \hbox{ for all } x\in\Omega \setminus\Omega_\eps.
\end{equation*}
Therefore, only $x$ belonging to $\Omega_\eps$ with
$\|v\|<\frac{\eps}{2\sqrt{Dt}}$ contribute to $N (t)$ and we can
write:
\begin{eqnarray*}
 &&N (t)=\int_{\R^n}\frac{1}{\pi^\frac{n}{2}}e^{-|v|^2}  \left[\int_{\Omega_\eps }\left(\mathds{1}_{\Omega_\eps }(x)-
\mathds{1}_{\Omega_\eps -2\sqrt{Dt}v}(x)\right)\dx\right]\dv+O\left(e^{-\frac{1}{t^{\delta}}}\right),
\end{eqnarray*}
where the exponentially small error with a $\delta>0$ is defined by
the integral
\begin{equation*}
\int_{\|v\|>\frac{\eps}{2\sqrt{Dt}}}\frac{1}{\pi^\frac{n}{2}}e^{-|v|^2}\dv.
\end{equation*}

Since $\del \Omega$ is regular, we introduce (see Section~\ref{secLC})  the local coordinates
$x=(\theta,s)$ and thus have $\hat{x}(\theta)\in
\del \Omega$ and $x\in \Omega_\eps$ iff $0<s<\eps$.  In this case,
$\chi_{2\sqrt{Dt},v}(x)=\mathds{1}_{\Omega_\eps }(x)-\mathds{1}_{\Omega_\eps -2\sqrt{Dt}v}(x)\ne 0$ 
iff $x\in \Omega_\eps$ and $\hat{x}(\theta)-s
n(\theta)+2\sqrt{Dt}v\notin \Omega.$ Moreover, with the notation
$(v,n)$ for the Euclidean inner product of two vectors in $\R^n$,
\begin{equation*}
(\hat{x}(\theta)-s n(\theta)+2\sqrt{Dt}v)\cdot n(\theta)=-s+2\sqrt{Dt}(v,n).
\end{equation*}
We deduce that 
\begin{equation*}
\chi_{2\sqrt{Dt},v}(x)\ne 0 \hbox{ iff } s-2\sqrt{Dt}(v,n)<0.
\end{equation*}
Consequently, if $(v,n)<0$, as $s>0$, it is not possible to have
$s-2\sqrt{Dt}(v,n)<0$.  In turn, if $0<(v,n)$ then
$s\in]0,2\sqrt{Dt}(v,n)[$.  Considering only $(v,n)>0$, we can define
$\eps= 2\sqrt{Dt}(v,n)$ and, since $v=\frac{x-y}{\sqrt{4Dt}}$ and
$x,y\in \Omega_{2\sqrt{Dt}(v,n)}$, we have $0<(v,n)<2$.  Thus, the
vector $v$ can be locally decomposed in two parts:
$v=((v,n),(v,\hat{x}))=(v_n,v_{\hat{x}})$.  Thus, returning to $N
(t)$, we obtain with the error $O(t)$ which comes from the Jacobian
approximation (see $|J(s,\theta)|$ in Section~\ref{secLC})
\begin{eqnarray*}
 &&N (t)=\int_{\R^{n-1}}\frac{1}{\pi^\frac{n-1}{2}}e^{-|v_{\hat{x}}|^2}d v_{\hat{x}}\int_0^2\frac{1}{\sqrt{\pi}}e^{-|v_n|^2}  
\left( \int_{\Omega_\eps}\chi_{2\sqrt{Dt},v_n}(x)\dx\right)d v_n+O(t)\\
&&=\int_{0}^2\frac{e^{-z^2}}{\sqrt{\pi}} \mu(\del\Omega,2\sqrt{Dt}z)dz+o(t^\frac{n-d}{2}).\quad 
\end{eqnarray*}
If $\del \Omega$ is regular, then $d=n-1$ and
$o(t^\frac{n-d}{2})=o(\sqrt{t})$, which, as it was mentioned, is actually $O(t)$.  The last formula that depends only on a
volume of the interior Minkowski sausage, holds for all types of
connected boundaries described in Subsection~\ref{SubsExFrac}.

The formula~(\ref{resVV}) follows from Eq.~(\ref{resV}) and the
relation
\begin{equation}\label{aprsocic}
 \mu(\del\Omega,\eps z)=z^{n-d}\mu(\del\Omega,\eps)+O(\eps^{2(n-d)}), 
\end{equation}
which, for a fixed $z$ and $\eps\to+0$, is evident for the regular
case and can be proved by approximating the fractal volume by a
converging sequence of the volumes for smooth boundaries.  For $d=n-1$
in Eq.~(\ref{resVV}), one has $\beta_1 = \frac{1-e^{-4}}{2\sqrt{\pi}}
\approx 0.2769$.
$\Box$

\begin{figure} [!htb] 
\begin{center}
\includegraphics[width=62mm]{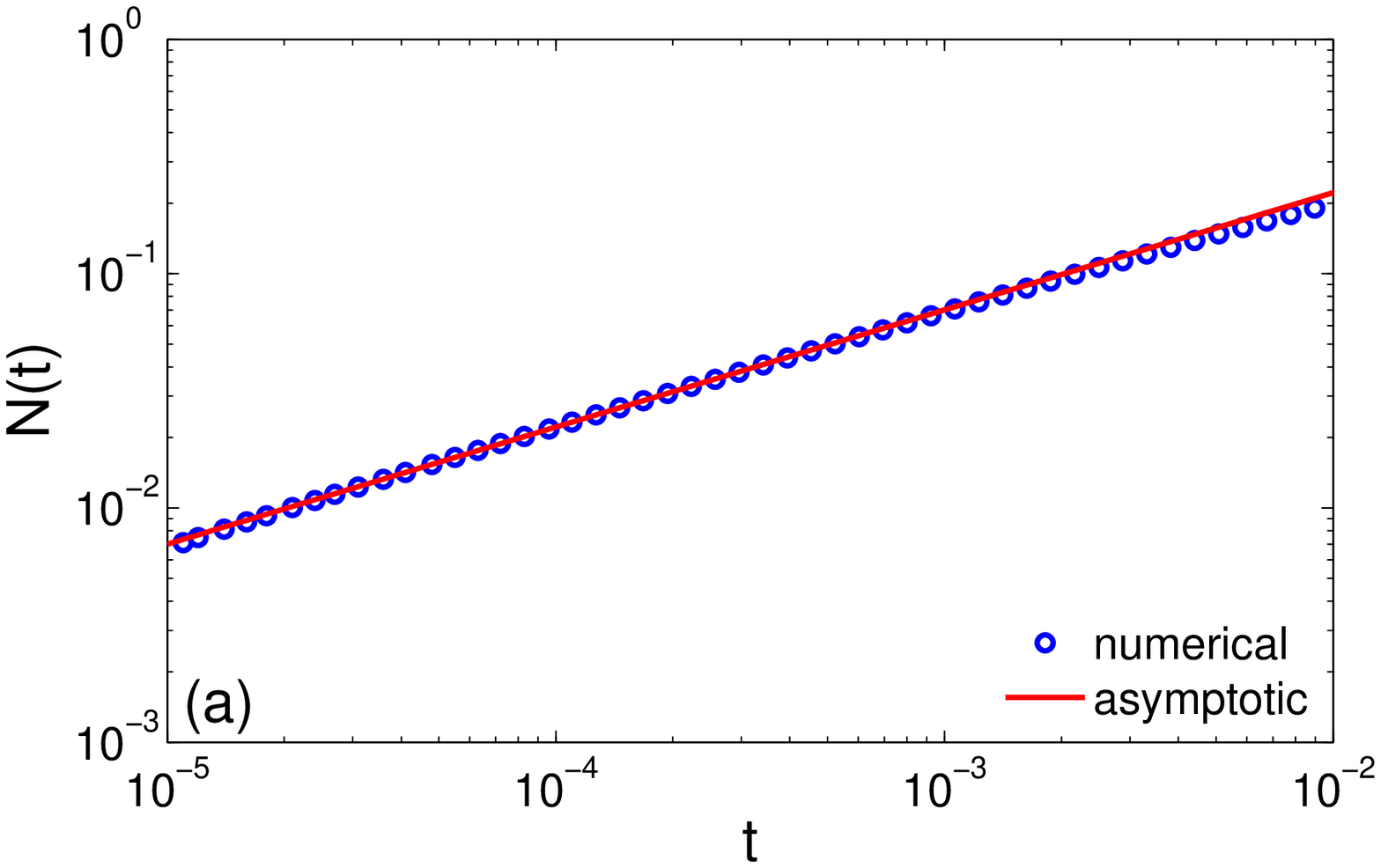}  
\includegraphics[width=62mm]{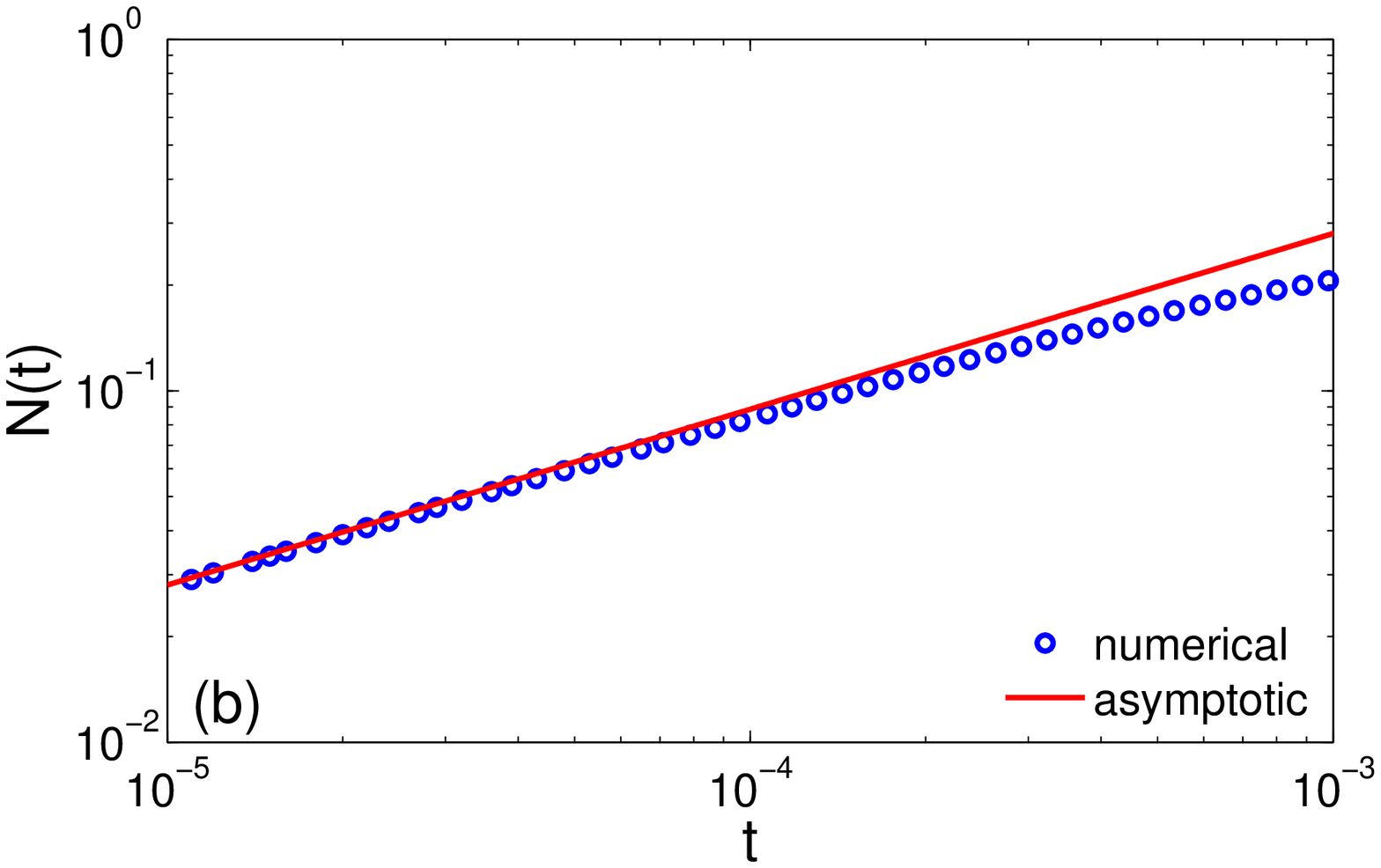}  
\end{center}
\caption{
Comparison between the asymptotic formula (\ref{resVV}) (solid line)
and a FreeFem++ numerical solution of the problem
(\ref{pch})--(\ref{ps}) (circles) for two domains: (a) the unit square
(with $\operatorname{Vol}(\del \Omega)=4$) and (b) the second
generation of the Minkowski fractal, with $\operatorname{Vol}(\del
\Omega)=2^2\cdot 4$.  We set $D_+=D_-=D=1$.}\label{FigCarInf}
\end{figure}

A comparison between the asymptotic formula~(\ref{resVV}) and a
numerical solution of the problem~(\ref{pch})--(\ref{ps}) is
illustrated in Fig.~\ref{FigCarInf} (for a square and a prefractal
domain).

\subsection{General case}

Let us come back to the problem~(\ref{prb1})--(\ref{prb1end}).

According to Theorem~\ref{TH1} (Eq.~(\ref{NOt2h})), the heat content
can be found up to the terms either $t^\frac{3}{2}$, or $t$ (depending
on values of $\lambda$), by integrating over all boundary points
$\theta$ of the solution $\hat{u}^{hom}$ of the homogeneous
problem~(\ref{homp01})--(\ref{homp0e}) with constant coefficients
taken at a boundary point $(0,\theta)$.  Obviously, Eq.~(\ref{NOt2h})
is valid only for regular boundaries.  Let us reformulate it to allow
an explicit calculation of the heat content for all types of
boundaries mentioned in Section~\ref{secWP}.

For this purpose, given $\eps=O(\sqrt{t})$, $\eps>\sqrt{4D_+t}$, we
divide $\del \Omega$ (which is still supposed to be regular) into $J$
disjoint parts $B_j$ ($j=1,\ldots,J$) of the size $\delta^{n-1}$ with
$0<\delta\le \eps$ such that $\del \Omega=\sqcup_{j=1}^{J} B_j$.

For $t\to +0$, $\delta\to 0$ and thus, due to regularity of $\del
\Omega$ on each $B_j\times ]0,\eps[$ the local change of variables
from Section~\ref{secLC} is a $C^1$-diffeomorphism.  In addition,
since $u$ continuously depends on $\lambda$ (see
Theorem~\ref{Thulamc}), $\hat{u}^{hom}$, considered as a function of
$\theta$, by the continuity of $\lambda$, is continuous on $\theta$.
Therefore, by the mean value theorem and due to the positivity of
$|J(s,\theta)|$, we deduce that for all $j=1,\ldots J$ there exists $
\theta_0^j\in \overline{B_j}$ such that
\begin{equation*}
\int_{B_j}d\theta \int_{[0,\eps]} 
\ds (1-\hat{u}^{hom}(s,\theta,t))|J(s,\theta)|=\int_{[0,\eps]} 
\ds (1-\hat{u}^{hom}(s,\theta_0^j,t)) \int_{B_j}d\theta |J(s,\theta)|.
\end{equation*}
From Eq.~(\ref{NNj}), Eq.~(\ref{NOt2h}) becomes
\begin{eqnarray*}
 &&N(t)-\sum_{j=1}^J\int_{[0,\eps]} 
\ds (1-\hat{u}^{hom}(s,\theta_0^j,t)) \int_{B_j}d\theta |J(s,\theta)| \\
&&=\left\{\begin{array}{ll}
       O(t~\mu(\del \Omega,\sqrt{t})),  & \quad  0<\lambda<\infty\\
         O(\sqrt{t}~\mu(\del \Omega,\sqrt{t})), & \quad  \lambda=\infty   \end{array}
 \right..
\end{eqnarray*}
Hence we prove
\begin{theorem} \label{FinLem} 
The heat content for the solution of the
problem~(\ref{prb1})--(\ref{prb1end}) can be explicitly found for all
types of boundaries $\del \Omega$ (a connected boundary of a compact
domain described in Subsection~\ref{SubsExFrac}) using the following
expressions:
\begin{enumerate}
 \item for $\lambda<\infty$ on $\del \Omega$:
\begin{eqnarray}
 &&N(t)=\frac{2\sqrt{t}}{\sqrt{D_+}\operatorname{Vol}(\del \Omega)} 
\left[\mu(\del \Omega,\sqrt{4D_+t})\int_{\del \Omega}d \sigma\lambda(\sigma)\int_1^{2}\dz  f(\sigma,z,t) \right.\nonumber\\
  &&\left. 
 - \int_1^{2} \dz \mu(\del \Omega,\sqrt{4D_+t}(z-1)) \int_{\del \Omega} d \sigma \lambda(\sigma) f(\sigma,z,t)
 \right.\nonumber\\
&&\left.-
 \int_0^{1}  \dz \mu(\del \Omega,\sqrt{4D_+t}z) \int_{\del \Omega} d \sigma \lambda(\sigma)f(\sigma,z,t)  \right]+O(t\mu(\del \Omega,\sqrt{t})),\label{NilocFint2C}
\end{eqnarray}
where $\alpha=\frac{1}{\sqrt{D_-}}+\frac{1}{\sqrt{D_+}}$ and
\begin{equation}\label{eqf}
f(\sigma,z,t)= \exp\left(2\lambda(\sigma) \alpha \sqrt{t}z+\lambda(\sigma)^2\alpha^2t\right) 
\operatorname{Erfc}(z+\lambda(\sigma) \alpha\sqrt{t}).
\end{equation}
\item for $\lambda=\infty$ on $\del \Omega$:
\begin{equation}
 N(t)=\frac{2\sqrt{D_-}}{\sqrt{D_-}+\sqrt{D_+}}\int_0^{2}\frac{e^{-z^2}}{\sqrt{\pi}}\mu(\del \Omega,\sqrt{4D_+t} z) \dz
+O(\sqrt{t}~\mu(\del \Omega,\sqrt{t})).\label{NilocFint2InfC}
 \end{equation}
\end{enumerate}
Formulas~(\ref{NilocFint2C}) and~(\ref{NilocFint2InfC}) can be
approximated by
\begin{enumerate}
 \item for $\lambda<\infty$ on $\del \Omega$:
\begin{eqnarray}
 &&N(t)=\frac{2\sqrt{t}~\mu(\del \Omega,\sqrt{4D_+t})}{\sqrt{D_+}\operatorname{Vol}(\del \Omega)} 
\left[\int_{\del \Omega}d \sigma\lambda(\sigma)\int_1^{2}\dz  f(\sigma,z,t) \right.\nonumber\\
  &&\left. 
 - \int_1^{2} \dz (z-1)^{n-d} \int_{\del \Omega} d \sigma \lambda(\sigma) f(\sigma,z,t)
 \right.\nonumber\\
&&\left.-
 \int_0^{1}  \dz z^{n-d} \int_{\del \Omega} d \sigma \lambda(\sigma)f(\sigma,z,t)  \right]+O(\sqrt{t}~\mu(\del \Omega,\sqrt{t})^2),\label{NFint2Ca}
\end{eqnarray}
\item for $\lambda=\infty$ on $\del \Omega$:
\begin{eqnarray}
 &&N(t)=\frac{2\sqrt{D_-}~ \beta_{n-d}}{\sqrt{D_-}+\sqrt{D_+}} \mu(\del \Omega,\sqrt{4D_+t}) 
+O(\mu(\del \Omega,\sqrt{t})^2), \label{NFint2InfCa}
\end{eqnarray}
where $\beta_x$ was defined in Eq. (\ref{eq:betax}).
\end{enumerate}
\end{theorem}

\textbf{Proof.}
Using Eqs.~(\ref{NNiloc})--(\ref{M1}), $N(t)$ becomes
\begin{eqnarray*}
 &&N(t)-\mu(\del \Omega,\eps) +\sum_{j=1}^J\int_{[0,\eps]^2} d s_1 \ds \; (h^{\theta_0^j}_+(s,s_1,t)-f^{\theta_0^j}_+(s,s_1,t))\int_{B_j}d\theta |J(s,\theta)|\\
 &&=\left\{\begin{array}{ll}
       O(t~\mu(\del \Omega,\sqrt{4D_+t})),  & \quad 0<\lambda<\infty\\ 
         O(\sqrt{t}~\mu(\del \Omega,\sqrt{4D_+t})), & \quad  \lambda=\infty   \end{array}
 \right..
\end{eqnarray*}

Let us calculate it explicitly.  We start with the part
\begin{equation*}
Nh_j(t)= \int_{[0,\eps]^2}d s_1 \ds \; h^{\theta_0^j}_+(s,s_1,t)\int_{B_j}d\theta |J(s,\theta)|.
\end{equation*}
Changing variables as in the proof of Theorem~\ref{TH1}, $Nh_j(t)$
becomes
\begin{eqnarray*}
 &&Nh_j(t)=\int_{B_j}d\theta \int_\R \frac{e^{-z^2}}{\sqrt{\pi}}
\left(\int_\R \mathds{1}_{[0,\eps]}(s)\mathds{1}_{[0,\eps]}(s-\sqrt{4D_+t}z)|J(s,\theta)| \ds\right)\dz\\
 &&+ a(\lambda,0,\theta_0^j) \int_{B_j}d\theta \int_\R \frac{e^{-z^2}}{\sqrt{\pi}}
\left(\int_\R \mathds{1}_{[0,\eps]}(s)\mathds{1}_{[0,\eps]}(-s+\sqrt{4D_+t}z) |J(s,\theta)|\ds\right)\dz.
\end{eqnarray*}

Therefore, we obtain
\begin{eqnarray*}
&&Nh_j(t)=\int_{B_j}d\theta \left[\int_{[0,\eps]}|J(s,\theta)|\ds\right.\\
&&-\left.\int_\R \frac{e^{-z^2}}{\sqrt{\pi}}\left(\int_{[0,\eps]} (\mathds{1}_{[0,\eps]}(s)-\mathds{1}_{[0,\eps]+
\sqrt{4D_+t}z}(s) ) |J(s,\theta)|\ds\right)\dz\right.\\
 &&+a(\lambda,0,\theta_0^j)\left. \int_\R \frac{e^{-z^2}}{\sqrt{\pi}}\left(\int_\R \mathds{1}_{[0,\eps]}(s)\mathds{1}_{[-\eps,0]+
\sqrt{4D_+t}z}(s) |J(s,\theta)|\ds\right)\dz\right] .
\end{eqnarray*}
Applying formula~(\ref{relation}) with $v=\sqrt{4D_+t}~ z$ (see also
Subsection~\ref{SSpartic}), we find
\begin{eqnarray*}
 &&Nh_j(t)=\mu(B_j,\eps)-(1-a(\lambda,0,\theta_0^j))\int_{0}^{ 2}\frac{e^{-z^2}}{\sqrt{\pi}} \mu(B_j,\sqrt{4D_+t}z)\dz.
\end{eqnarray*}
Thus, for $\lambda<\infty$ $Nh_j(t)=\mu(B_j,\eps)$ since $a=1$.

We treat the second part in the same way,
\begin{equation*}
Nf_j(t)=- \int_{[0,\eps]^2}d s_1 \ds \; f^{\theta_0^j}_+(s,s_1,t)\int_{B_j}d\theta|J(s,\theta)|,
\end{equation*}
which is equal to zero for $\lambda=\infty$.  For $f(\theta_0^j,z,t)$
from Eq.~(\ref{eqf}), we find that
\begin{eqnarray*}
&&Nf_j(t)=- \frac{2\lambda(\theta_0^j)\sqrt{t}}{\sqrt{D_+}}\int_{\R^2}\ds \dz   
\mathds{1}_{[0,\eps]}(s)\mathds{1}_{[0,\eps]}(-s+2\sqrt{D_+t}z)
 f(z,t) \int_{B_j}|J(s,\theta)|  d\theta\\
&&=-\frac{2\lambda(\theta_0^j)\sqrt{t}}{\sqrt{D_+}}\left[\int_1^{2} \dz f(\theta_0^j,z,t) 
\int_{B_j}\int_{(z-1)\sqrt{4D_+t}}^{\sqrt{4D_+t}}|J(s,\theta)|\ds d\theta\right.\\
&&\left.+ \int_0^{1} \dz f(\theta_0^j,z,t) \int_{B_j}\int_{0}^{\sqrt{4D_+t}z}|J(s,\theta)|\ds d\theta\right] \\ 
&&=-\frac{2\lambda(\theta_0^j)\sqrt{t}}{\sqrt{D_+}}\left[ \mu(B_j,\sqrt{4D_+t})\int_1^{2}  
f(\theta_0^j,z,t)\dz  \right. \\
&& \left. -\int_1^{2}f(\theta_0^j,z,t) \mu(B_j,\sqrt{4D_+t}(z-1))\dz 
+  \int_0^{1}  f(\theta_0^j,z,t)\mu(B_j,\sqrt{4D_+t}z)\dz  \right] .
\end{eqnarray*}
Putting two results together, we obtain the following approximations
for $N(t)$:
\begin{enumerate}
 \item for $\lambda<\infty$ on $\del \Omega$:
\begin{eqnarray}
 &&N(t)=\sum_{j=1}^J\mu(B_j,\sqrt{4D_+t}) 
 \frac{2\lambda(\theta_0^j)\sqrt{t}}{\sqrt{D_+}}\int_1^{2}  f(\theta_0^j,z,t)\dz\nonumber\\
  && 
 -\sum_{j=1}^J\frac{2\lambda(\theta_0^j)\sqrt{t}}{\sqrt{D_+}}\left[\int_1^{2}f(\theta_0^j,z,t)
\mu(B_j,\sqrt{4D_+t}(z-1))\dz \right.\nonumber\\
&&\left.-
 \int_0^{1}  f(\theta_0^j,z,t)\mu(B_j,\sqrt{4D_+t}z)\dz  \right]+ O(t~\mu(\del \Omega,\sqrt{t})),\label{NilocFint2}
\end{eqnarray}
\item for $\lambda=\infty$ on $\del \Omega$:
\begin{eqnarray}
\nonumber
 N(t) &=&\frac{2 \sqrt{D_-}}{\sqrt{D_-}+\sqrt{D_+}}\sum_{j=1}^J\int_0^{2}\frac{e^{-z^2}}{\sqrt{\pi}}\mu(B_j,\sqrt{4D_+t}~ z) \dz \\
&+& O(\sqrt{t}~\mu(\del \Omega,\sqrt{t})).\label{NilocFint2Inf}
 \end{eqnarray}
\end{enumerate}
It means that if the formulas for $\mu(B_j,\delta)$ are known, we get
the approximation of $N(t)$ up to terms of the order of
$t^\frac{n-d+2}{2}$ for $\lambda<\infty$, and of the order of
$t^\frac{1+n-d}{2}$ for $\lambda=\infty$.  Moreover, this
approximation, depending only on the volume of $\del \Omega$, holds
for all types of boundaries, even fractals (see
Subsection~\ref{SubsExFrac} and p. 378 of Ref.~\cite{Fleckinger}
for a similar conclusion).

Let us now change the sum over $j$ with the integral over $z$ and make
$J\to+\infty$:
\begin{equation*}
\lim_{J\to +\infty}\sum_{j=1}^J C(z,t,\theta_0^j)\mu(B_j,\sqrt{4D_+t} z)= \int_{\del \Omega} C(z,t,\sigma)\operatorname{dist}(\sigma,\sqrt{4D_+t} z) d \sigma,
\end{equation*}
where $d\sigma$ is understood in the sense of the Hausdorff measure ($d$-measure)
defined on $\del \Omega$.  Thus, again with the help of the mean value
theorem, we have
\begin{equation*}
\int_{\del \Omega} C(z,t,\sigma)\operatorname{dist}(\sigma,\sqrt{4D_+t} z) d \sigma=
\frac{\mu(\del \Omega,\sqrt{4D_+t})}{\operatorname{Vol}(\del \Omega)}\int_{\del \Omega} C(z,t,\sigma)d \sigma,
\end{equation*}
from which Eqs.~(\ref{NilocFint2C}) and~(\ref{NilocFint2InfC}) follow.
We use Eq.~(\ref{aprsocic}) to obtain formulas~(\ref{NFint2Ca})
and~(\ref{NFint2InfCa}).
$\Box$

\section{Regular case}\label{SecR}

In the case of a regular boundary we provide the asymptotic expansion
of the heat content up to the third-order term.

In this case, we can approximate the solution of the
system~(\ref{eq47})--(\ref{endsys}) by the solution $v$ of the
following problem (instead of~(\ref{consteq1})--(\ref{consteqN}), as
previously)
\begin{eqnarray}
&&\frac{\del}{\del t}u_+-D_+\left(\frac{\del^2}{\del s^2}+\sum_{i=1}^{n-1}\frac{\del^2}{\del \theta^2_i}\right) 
u_+ +D_+\sum_{i=1}^{n-1}k_i(\theta_0)\frac{\del u_+}{\del s}=0,\; 0<s<\eps\label{consteq1r}\\
&&\frac{\del}{\del t}u_--D_-\left(\frac{\del^2}{\del s^2}+\sum_{i=1}^{n-1}\frac{\del^2}{\del \theta_i^2}\right) 
u_-  +D_-\sum_{i=1}^{n-1}k_i(\theta_0) \frac{\del u_-}{\del s}=0, \; -\eps<s<0\\
&&u_+|_{t=0}=1,\; u_-|_{t=0}=0,\\
&& D_-\frac{\del u_-}{\del s}|_{s=-0}=\lambda(\theta_0)(u_--u_+)|_{s=0},\nonumber \\
&& D_+\frac{\del u_+}{\del s}|_{s=+0}=D_- \frac{\del u_-}{\del s}|_{s=-0}.\label{consteqNr}
\end{eqnarray}
In this approximation the remainder terms of the
system~(\ref{eq47})--(\ref{endsys}) contain only the coefficients of
the order $\sqrt{t}$ (to compare with~(\ref{Rest0S1})):
\begin{eqnarray*}
&&R(s_1,\theta_0)=s_1D_+\sum_{i=1}^{n-1} k_i^2(\theta_0) +O(s_1^2) ,
\end{eqnarray*}
that gives
\begin{eqnarray}
&&R(s_1,\theta_0)=(s\mp2\sqrt{D_+(t-\tau)}z_1)D_+\sum_{i=1}^{n-1} k_i^2(\theta_0)+O(t)= O(\sqrt{t})\label{Rest0S1r}.
\end{eqnarray}

The basis of the parametrix is the Green function given by (see
Section~\ref{SecGF2})
\begin{eqnarray}
&&h_+^{\theta_0}(s_1,s_2,t)= \frac{1}{\sqrt{4\pi D_+t}}\left(\exp\left(-\frac{(s_1-s_2-tD_+\gamma(\theta_0))^2}{4D_+t}\right) \right.\nonumber\\
&&\left.
+a(\lambda,0,\theta_0) \exp\left(-\frac{(s_1+s_2-tD_+\gamma(\theta_0))^2}{4D_+t}\right) \right),\label{hr}
\end{eqnarray}
\begin{eqnarray}
&&f_+^{\theta_0}(s_1,s_2,t)=b(\lambda,0,\theta_0)\frac{\lambda(\theta_0)}{D_+}\cdot\nonumber\\
&&\cdot
\exp\left(\frac{\lambda(\theta_0) \alpha}{\sqrt{D_+}}(s_1+s_2-tD_+\gamma(\theta_0))+ 
\lambda(\theta_0)^2 \alpha^2 t\right)\cdot\nonumber\\
&&\cdot \operatorname{Erfc}\left(\frac{s_1+s_2-tD_+\gamma(\theta_0)}{2\sqrt{D_+t}}+\lambda(\theta_0) \alpha\sqrt{ t} \right),\label{fpr}
\end{eqnarray}
where $\gamma(\theta_0)=\sum_{i=1}^{n-1}k_i(\theta_0)$.  The
estimate~(\ref{NNj}) becomes
\begin{equation}
  N(t)-\int_{\del \Omega}d\theta_0\int_{[0,\eps]}\ds\; (1-\hat{u}^{hom}_\eps(s,\theta_0,t))
|J(s,\theta_0)|=\left\{\begin{array}{ll}
       O(t^2),  & 0<\lambda<\infty\\                                                                                                                            
         O(t^\frac{3}{2}), &   \lambda=\infty   \end{array}
 \right.\label{NOt2hN}
\end{equation}
Consequently, for the regular case we have 
\begin{theorem}\label{ThFinalr}
Let $\Omega$ be a compact domain of $\R^n$ with a connected boundary
$\del \Omega\in C^\infty(\R^n)$.  Then for $\lambda=\infty$ we have
\begin{equation}\label{NtregLI}
N(t)=2 \frac{1-e^{-4}}{\sqrt{\pi}} ~ \frac{\sqrt{D_+ D_-}}{\sqrt{D_+}+\sqrt{D_-}} ~ \operatorname{Vol}(\del \Omega)~ \sqrt{t} +O(t^{\frac{3}{2}}).
\end{equation} 
In the case of $0<\lambda<\infty$, we have
\begin{eqnarray}
   &&N(t)=4 C_0t \int_{\del \Omega}\lambda(\sigma) d \sigma - 
\frac{2}{3}C_1t^\frac{3}{2}\biggl[2 \biggl(\frac{1}{\sqrt{D_+}}+\frac{1}{\sqrt{D_-}}\biggr) \int_{\del \Omega}\lambda^2(\sigma)d \sigma  \nonumber\\
   &&  - \sqrt{D_+}(n-1)\int_{\del \Omega}\lambda(\sigma) H(\sigma)d \sigma\biggr]+O(t^2),\label{NtregLF}
\end{eqnarray}
where $H$ is the mean curvature, and 
\begin{eqnarray}
\label{eq:C0}
 &&C_0=1+\frac{3}{2}\operatorname{erf}(1)-\frac{9}{4}\operatorname{erf}(2)+\frac{1}{\sqrt{\pi}}\left(\frac{1}{e}-\frac{1}{e^4} \right) \approx 0.2218,\\
 &&C_1=\frac{1}{\sqrt{\pi}} -6+\frac{5e^{-4}-4e^{-1}}{\sqrt{\pi}}-5\operatorname{erf}(1)+11\operatorname{erf}(2) \approx 0.5207.
\end{eqnarray}
\end{theorem}

\textbf{Proof.}
Let us consider the case $\lambda=\infty$.  Using the Green function
given in Eqs.~(\ref{hr}) and~(\ref{fpr}), we obtain
\begin{equation} \label{NilocFint2r} 
 N(t) =\sum_{j=1}^J\beta\int_0^{2-\sqrt{tD_+}\frac{\gamma(\theta_0^j)}{2}}\frac{e^{-z^2}}{\sqrt{\pi}}
 \mu(B_j,\sqrt{4D_+t} z+tD_+\gamma(\theta_0^j)) dz   
 +O(t^\frac{3}{2}), 
\end{equation}
where $\beta=\frac{2\sqrt{D_-}}{\sqrt{D_-}+\sqrt{D_+}}$.  In
Eq.~(\ref{NilocFint2r}) the remainder term also contains the integrals
$\int^0_{-\sqrt{tD_+}\frac{\gamma(\theta_0^j)}{2}}dz$.  From
Eq.~(\ref{NilocFint2r}) we find
\begin{eqnarray*}
 N(t) && =\sum_{j=1}^J \beta\int_0^{2-\sqrt{tD_+}\frac{\gamma(x_j)}{2}}\frac{e^{-z^2}}{\sqrt{\pi}}\int_{B_j} 
d \theta \int_0^{2\sqrt{D_+t}z+ \gamma(x_j)D_+t}ds (1-s (n-1) H(\theta)) \\
&& +O(t^\frac{3}{2}).
\end{eqnarray*}
Therefore, we have
\begin{eqnarray}
 &&N(t)=\sqrt{t}\left(2C \sqrt{D_+}\sum_{j=1}^J \operatorname{Vol}(B_j)\right)\nonumber\\
 &&- t(n-1)\left(\sum_{j=1}^J  \xi\left[\operatorname{Vol}(B_j)H(x_j)-\int_{B_j}H(\sigma)d \sigma\right] \right)+O(t^{\frac{3}{2}}),\label{resultreg2czz}
\end{eqnarray}
where  
\begin{eqnarray*}
 && C=\frac{1-e^{-4}}{\sqrt{\pi}} ~ \frac{\sqrt{D_-}}{\sqrt{D_+}+\sqrt{D_-}} , \\
 &&\xi=\left(4\frac{e^{-4}}{\sqrt{\pi}}-\operatorname{erf}(2)\right)  \frac{D_+\sqrt{D_-}}{\sqrt{D_+}+\sqrt{D_-}}.
\end{eqnarray*}
In addition, for all $\sigma\in B_j$, the distance between $x_j$
(which also belongs to $B_j$) and $\sigma$ goes to $0$ as $J\to
+\infty$.  Thus, since
\begin{equation*}
|H(x_j)-H(\sigma)|\le H'(\sigma) |x_j-\sigma|\le C\operatorname{Vol}(B_j),
\end{equation*}
we have
\begin{equation*}
\lim_{J\to+\infty}\sum_{j=1}^J\left|\operatorname{Vol}(B_j)H(x_j)-\int_{B_j}H(\sigma)d \sigma\right|=0.
\end{equation*}
Hence, from Eq.~(\ref{resultreg2czz}) we obtain Eq.~(\ref{NtregLI}).

The case $0<\lambda<\infty$ can be treated in the similar way using in
Eq.~(\ref{NilocFint2}) the expansion of the $f(\sigma,t,z)$:
\begin{eqnarray*}
 &&f(\sigma,t,z)=\exp\left(2\lambda(\sigma)\alpha\sqrt{t}z+ \lambda^2(\sigma) \alpha^2(\sigma) t\right)\operatorname{Erfc}(z+\lambda(\sigma)\alpha\sqrt{t})\\
 &&=\operatorname{Erfc}(z)-2\lambda(\sigma)\alpha\sqrt{t}\left(\frac{1}{\sqrt{\pi}}e^{-z^2}- z \operatorname{Erfc}(z)\right)+O(t).
\end{eqnarray*} 
$\Box$

\section*{Acknowledgment}

We thank Fr\'ed\'eric Abergel and Fran\c{c}ois Golse for helpful
discussions and Jean-Baptiste Apoung Kamga for advises on FreeFem++
simulations.

\appendix

\section{Definitions of Besov spaces on fractals}~\label{ApBS}

Let us define the Besov space $B^{2,2}_{\beta}(\del \Omega)$ on a
$d$-set $\del \Omega$ (see Ref.~\cite{JW} p.135 and
Ref.~\cite{Wallin}).

There are many equivalent definitions\cite{JWBesov,Triebel} of Besov
spaces.  To give one of them, we introduce\cite{JW,JWBesov} a net
$\mathcal{N}$ with mesh $2^{-\nu}$, $\nu\in \N$, i.e. a division of
$\R^n$ into half-open non-overlapping cubes $W$ with edges of length
$2^{-\nu}$, obtained by intersecting $\R^n$ with hyperplanes
orthogonal to the axes.  In addition, we denote by
$\mathcal{P}_k(\mathcal{N})$ the set of functions which on each cube
$W$ in the net $\mathcal{N}$ coincide with a polynomial of degree at
most $k$.

\begin{definition}\textbf{(Besov space $B^{p,q}_\beta(\Gamma),$ $\beta>0$, see Ref.~\cite{JW})}
Let $\Gamma$ be a closed subset of $\R^n$ which is a $d$-set
preserving Markov's inequality for $0<d\le n$ and let $m_d$ be a fixed
$d$-measure on $\Gamma$.  We say that $f\in B^{p,q}_\beta(\Gamma)$,
$\beta>0$, $1\le p,q\le +\infty$, if $f\in L^p(m_d)$ and there is a
sequence $B=(B_\nu)_{\nu\in\N}\in \ell^q$ such that for every net
$\mathcal{N}$ with mesh $2^{-\nu}$, $\nu\in \N$ there exists a
function $s(\mathcal{N})\in \mathcal{P}_{[\beta]}(\mathcal{N})$ (by
$[\beta]$ is denoted the integer part of $\beta$) satisfying
\begin{equation*}
\|f-s(\mathcal{N})\|_{L^p(m_d)}\le 2^{-\nu \beta} B_\nu.
\end{equation*}
The norm of $f$ in $B^{p,q}_\beta(\Gamma)$ is given by the formula
\begin{equation*}
\|f\|_{B^{p,q}_\beta(\Gamma)}=\|f\|_{L^p(m_d)}+\inf_{B}\|B\|_{\ell^q},
\end{equation*}
where the infimum is over all such sequences $B$.
\end{definition}

The dual Besov space $(B^{2,2}_{\beta}(\del
\Omega))'=B^{2,2}_{-\beta}(\del \Omega)$ is introduced in
Ref.~\cite{JWBesov}. To give the definition of the Besov space
$B^{2,2}_{-\beta}(\del \Omega)$ we need to define the atoms:
\begin{definition}\textbf{(Atom\cite{JWBesov})}
Let $\beta>0$, $1\le p\le \infty$, and let $W$ with $W\cap \Gamma\ne
\varnothing$ be a cube with edge length $2^{-\nu}$, $\nu\in \N$.  A
function $a=a_W\in L^p(m_d)$ is a $(-\beta,p)$-atom associated with
$W$ if
\begin{enumerate}
 \item $\operatorname{supp} a\subset 2W$, where $2W$ is the cube obtained by expanding $W$ twice from its center,
 \item $\int x^\gamma a(x) d m_d=0$ for $|\gamma|\le [\beta]$ if $\nu>0$,
 \item $\|a\|_{L^p(m_d)}\le 2^{\nu \beta}$.
\end{enumerate}
\end{definition}
Let $\mathcal{N}_\nu(\Gamma)=\{W\in \mathcal{N}_\nu|\; W\cap \Gamma
\ne \varnothing\}$ with the notation $\mathcal{N}_\nu$ of the net with
mesh $2^{-\nu}$ such that the origin is a corner of some cube in the
net.  Then we can define the Besov space with a negative parameter
$-\beta$, $B^{2,2}_{-\beta}(\del \Omega)$, which is
actually\cite{JWBesov} the dual Besov space of $B^{2,2}_{\beta}(\del
\Omega)$:
\begin{definition}\textbf{(Besov space $B^{p,q}_{-\beta}(\Gamma),$ $\beta>0$, see Ref.~\cite{JW})}
The space $B^{p,q}_{-\beta}(\Gamma),$ $\beta>0$, $1\le p,q \le \infty$
consists of functions $f\in \mathcal{D}'(\R^n)$ which are given by
\begin{equation*}
\forall \phi\in \mathcal{D}(\R^n) \quad \langle f,\phi \rangle=\sum_{\nu\in \N} \sum_{W\in \mathcal{N}_\nu(\Gamma)}s_W\int a_W\phi d m_d ,
\end{equation*}
where $a_W$ are $(-\beta,p)$-atoms and $s_W$ are numbers such that
$S=(S_\nu)_{\nu \in \N}\in \ell^q$ and $S_\nu$ is defined by
\begin{equation*}
S_\nu=\left(\sum_{W\in \mathcal{N}_\nu(\Gamma)}|s_W|^p \right)^\frac{1}{p}.
\end{equation*}
The norm of $f$ is defined by
\begin{equation*}
\|f\|_{B^{p,q}_{-\beta}(\Gamma)}=\inf\|S\|_{\ell^q},
\end{equation*}
where the infimum is taken over all possible atomic decompositions of
$f$:
\begin{equation*}
f=\sum_{\nu\in \N}\sum_{W\in \mathcal{N}_\nu(\Gamma)} s_Wa_W.
\end{equation*}
\end{definition}

\section{Explicit computations for half space problem with constant coefficients}~\label{SecGreen}

\subsection{Case $\lambda=\infty$}\label{SecGF2}

The Green function of the one-dimensional
problem~(\ref{homp01})--(\ref{homp0e}) with $\lambda=\infty$ and $s\in
\R$ was treated in Ref.~\cite{PirN,Carslaw} and it is given by
\begin{equation*}
\Gamma(s,s_1,t)=  \mathds{1}_{\{s>0,\;s_1>0\}} \Gamma_{++}(s,s_1,t)+ 
\mathds{1}_{\{s<0,\; s_1>0\}}\Gamma_{-+}(s,s_1,t)
\end{equation*}
with
\begin{eqnarray}
 && \Gamma_{++}(s,s_1,t)=\frac{1}{\sqrt{4\pi D_+t}}\left( \exp\left(-\frac{(s-s_1)^2}{4D_+t}\right) + 
A  \exp\left(-\frac{(s+s_1)^2}{4D_+t}\right)\right),\label{GammaP}\\
 && \Gamma_{-+}(s,s_1,t)=B\frac{1}{\sqrt{\pi D_+t}}  
\exp\left(-\frac{\left(s-s_1\sqrt{\frac{D_-}{D_+}}\right)^2}{4D_-t}\right),\label{GammaM}
\end{eqnarray}
where $A=\frac{\sqrt{D_+}-\sqrt{D_-}}{\sqrt{D_+}+\sqrt{D_-} }$ and
$B=\frac{\sqrt{D_+}}{\sqrt{D_+}+\sqrt{D_-} }$.

Let us use this result to find the Green function
$\Gamma^{reg}(s,s_1,t)$ of the following one-dimensional problem
\begin{eqnarray}
 &&\frac{\del}{\del t}u_+-D_+\frac{\del^2}{\del s^2} u_+ +D_+k\frac{\del}{\del s} u_+=0,\; s>0 \label{es1}\\
&&\frac{\del}{\del t}u_--D_-\frac{\del^2}{\del s^2} u_- +D_-k\frac{\del}{\del s} u_-=0, \; s<0\\
&&u_+|_{t=0}=1,\; u_-|_{t=0}=0,\\
&&u_+|_{s=+0}=u_-|_{s=-0}, \quad D_+\frac{\del }{\del s} u_+|_{s=+0}=D_-\frac{\del }{\del s}u_-|_{s=-0},\label{esn}
\end{eqnarray}

The constant coefficient problem
\begin{eqnarray}
 &&\frac{\del}{\del t}u-D\frac{\del^2}{\del s^2} u +D k\frac{\del}{\del s} u=0,\; s\in \R,~\label{fdseq}\\
&&u|_{t=0}=u_0,
\end{eqnarray}
has the Green function of the form
\begin{equation*}
K(s,s_1,t)=\frac{1}{\sqrt{4\pi Dt}}e^{-\frac{(s-y-tDk)^2}{4Dt}},
\end{equation*}
that means that the change of variables $s-tDk=X$
transforms~(\ref{fdseq}) to
\begin{equation*}
\frac{\del}{\del t}u-D\frac{\del^2}{\del X^2} u=0
\end{equation*}
with the Green function
\begin{equation*}
K_0(X,Y,t)=\frac{1}{4\pi Dt}e^{-\frac{(X-Y)^2}{4Dt}}.
\end{equation*}
In addition\cite{PirN}, we know (see~(\ref{GammaP})--(\ref{GammaM}))
the Green function for the constant coefficient problem
\begin{eqnarray*}
 &&\frac{\del}{\del t}u_+-D_+\frac{\del^2}{\del X^2} u_+=0 \quad X>0,\\
&&\frac{\del}{\del t}u_--D_-\frac{\del^2}{\del X^2} u_-=0 \quad X<0,\\
&&u_+|_{t=0}=1,\quad u_-|_{t=0}=0,\\
&&u_+|_{X=+0}=u_-|_{X=-0}, \quad D_+ \frac{\del }{\del X} u_+|_{X=+0}=D_-\frac{\del }{\del X}u_-|_{X=-0}.
\end{eqnarray*}
Consequently, we perform the following change of variables in
Eqs. (\ref{es1})--(\ref{esn}):
\begin{eqnarray*}
 &&X=\mathds{1}_{s>0}(s)\left(s+tD_+k\right)+\mathds{1}_{s<0}(s)\left(\frac{D_-}{D_+}s-tD_-k\right)
\end{eqnarray*}
and obtain for $z=\frac{D_+}{D_-}X$ that
\begin{eqnarray*}
 &&\frac{\del}{\del t}u_+-D_+\frac{\del^2}{\del X^2} u_+=0 \quad X>tD_+k,\\
&&\frac{\del}{\del t}u_--D_-\frac{\del^2}{\del z^2} u_-=0 \quad z<tD_+k,\\
&&u_+|_{t=0}=1,\quad u_-|_{t=0}=0,\\
&&u_+|_{X=+tD_+k}=u_-|_{z=-tD_+k}, \quad D_+ \frac{\del }{\del X}u_+|_{X=+tD_+k}=D_-\frac{\del }{\del z}u_-|_{z=-tD_+k}.
\end{eqnarray*}
Thus,
\begin{eqnarray*}
 \Gamma_{++}^{reg}(s,s_1,t) && = \frac{1}{\sqrt{4\pi D_+t}}\left(\exp\left(-\frac{(s-s_1-tD_+k)^2}{4D_+t}\right) \right. \\
 && \left. +A \exp\left(-\frac{(s+s_1-tD_+k)^2}{4D_+t}\right)\right) 
\end{eqnarray*}
for $s,s_1>0$, and
\begin{equation*}
 \Gamma_{-+}^{reg}(s,s_1,t) = \frac{1}{ \sqrt{\pi D_+}t}B \exp\left(-\frac{([D_+/D_-] s-s_1\sqrt{D_+/D_-}+tD_+k)^2}{4D_- t}\right)
\end{equation*}
for $s<0,~ s_1>0$.
Now, to obtain the Green function of the multidimensional problem
\begin{eqnarray*}
 &&\frac{\del}{\del t}u_+-D_+ \left(\frac{\del^2}{\del s^2}+\sum_{i=1}^{n-1}\frac{\del^2}{\del \theta_i^2}\right) u_+ 
+D_+k\frac{\del}{\del s} u_+=0,\; s>0,\; \theta_i\in \R, \label{es1m}\\
&&\frac{\del}{\del t}u_--D_- \left(\frac{\del^2}{\del s^2}+\sum_{i=1}^{n-1}\frac{\del^2}{\del \theta_i^2}\right)  u_- 
+D_-k\frac{\del}{\del s} u_-=0, \; s<0,\;  \theta_i\in \R\\
&&u_+|_{t=0}=1,\; u_-|_{t=0}=0,\\
&&u_+|_{s=+0}=u_-|_{s=-0}, \quad D_+\frac{\del }{\del s} u_+|_{s=+0}=D_-\frac{\del }{\del s}u_-|_{s=-0},\label{esnm}
\end{eqnarray*}
we apply the Fourier transform in $s_i$ variables and, due to the
boundary conditions depending only on $s$, we obtain that
$\hat{G}_{\pm +}$, the Fourier transform of the Green function
$G_{\pm+}$, can be found by the formula
\begin{equation*}
\hat{G}_{\pm +}(s,s_1,\xi, t)=e^{-D_{\pm}|\xi|^2 t} \Gamma_{\pm +}(s,s_1,t),
\end{equation*}
where $\Gamma_{\pm +}(s,s_1,t)$ is the Green function of the
corresponding one-dimensional problem.  This implies
\begin{eqnarray*}
 &&\Gamma_{++}^{reg}(s,\theta,s_1,\theta_1,t)= \frac{1}{(4\pi D_+t)^\frac{n}{2}}\left( \exp\left(-\frac{(s-s_1-tD_+k)^2}{4D_+t}\right) \right. \\
 && \left. +A  \exp\left(-\frac{(s+s_1-tD_+k)^2}{4D_+t}\right)\right)
 \exp\left( -\frac{|\theta-\theta_1|^2}{4D_+t}\right)\quad\hbox{ for }  s,s_1>0,\; \theta, \theta_1\in \R^{n-1},
\end{eqnarray*}
\begin{eqnarray*}
&&\Gamma_{-+}^{reg}(s,\theta,s_1,\theta_1,t)= \frac{B}{(4\pi D_-t)^\frac{n-1}{2} \sqrt{\pi D_+}t}
\exp\left(-\frac{(\frac{D_+}{D_-} s-s_1\sqrt{\frac{D_+}{D_-}}+tD_+k)^2}{4D_- t}\right) \\
&& \cdot \exp\left(-\frac{|\theta-\theta_1|^2}{4D_-t}\right) \quad\hbox{for } s<0, \; s_1>0,\;\theta, \theta_1\in \R^{n-1}.
\end{eqnarray*}

For $\Gamma$ and $\Gamma^{reg}$ we also have Varadhan's
bounds\cite{McKean} for $s\ne s_1$
\begin{eqnarray*}
 &&\lim\limits_{t\to 0+} t\ln \Gamma_{++}(s,s_1,t)=\lim\limits_{t\to 0+} t\ln \Gamma_{++}^{reg}(s,s_1,t)=-\frac{d(s,s_1)^2}{4D_+},\\
&&\lim\limits_{t\to 0+} t\ln \Gamma_{-+}(s,s_1,t)=\lim\limits_{t\to 0+} t\ln \Gamma_{-+}^{reg}(s,s_1,t)=-\frac{d\bigl(s,s_1\sqrt{\frac{D_-}{D_+}}\bigr)^2}{4D_-},
\end{eqnarray*}
where $d(s,s_1)$ is the Riemannian distance between $s$ and $s_1$,
which is equal here to the Euclidean distance, since $D_+$ and $D_-$
are constant in $\Omega_+$ and $\Omega_-$ respectively.

\subsection{Case $0<\lambda<\infty$}\label{SubsecFGLF}

Let us consider the one-dimensional
problem~(\ref{homp01})--(\ref{homp0e}) with
$\lambda\equiv\lambda(\theta_0)$ and $s\in \R$.  The associated
problem for the heat kernel is then given by
\begin{eqnarray}
 &&\left(\del_t -D_\pm \del^2_{s}\right) G(s,s_1,t)=0,\nonumber\\
&&G|_{t=0}=\delta(s,s_1) \quad \hbox{for } s>0,\nonumber\\
&&D_- \frac{\del}{\del s} G(-0,s_1,t)=\lambda(G(-0,s_1,t)-G(+0,s_1,t)),\label{BC1}\\
&& D_+  \frac{\del}{\del s} G(+0,s_1,t)=D_- \frac{\del}{\del s} G(-0,s_1,t).\label{BC2}
\end{eqnarray}
We search the explicit solution of the problem\cite{PirN} with
\begin{equation*}
  G(s,s_1,t)=\left\{\begin{array}{c}
              G_{-+}, \quad s<0,s_1>0\\
	      G_{++}, \quad s>0,s_1>0
             \end{array}
 \right..
\end{equation*}
We seek for $G_{-+}$ and $G_{++}$ in terms of free heat kernel
$K(s,s_1,D_\pm t)$ (see Eq.~(\ref{Keq})) and single layer heat
potentials for $s_1 > 0$:
\begin{eqnarray*}
 &&G_{++}(s,s_1,t)=K(s,s_1,D_+ t)+D_+\int_0^t K(s,0, D_+(t-\tau))\alpha_+(s_1,\tau)  \dtau\; \quad   (s>0),  \\
&&G_{-+}(s,s_1,t)=D_-\int_0^t K(s,0, D_-(t-\tau))\alpha_-(s_1,\tau)  \dtau\;  \quad (s<0),
\end{eqnarray*}
where $\alpha_\pm(s_1,\tau)$ are unknown densities to be determined.
Considering the boundary conditions~(\ref{BC1})--(\ref{BC2}) and the
jumps of the first derivatives of $G_{\pm+}$,
\begin{eqnarray*}
 &&\frac{\del}{\del s}G_{++}|_{s=+0}=-\frac{1}{2}\alpha_+(s_1,t)+\frac{\del}{\del s}K(0,s_1,D_+ t),\\
 &&\frac{\del}{\del s}G_{-+}|_{s=-0}=-\frac{1}{2}\alpha_-(s_1,t),
\end{eqnarray*}
we obtain two relations
\begin{eqnarray*}
 D_-\alpha_-(s_1,t) &&= -D_+\alpha_+(s_1,t)+2 D_+\frac{\del}{\del  s}K(0,s_1,D_+t),\\
 D_-\alpha_-(s_1,t) && =2\lambda K(0,s_1,D_+t)+\lambda\frac{\sqrt{D_+}}{\sqrt{\pi}}\int_0^t\frac{\alpha_+(s_1,\tau)}{\sqrt{t-\tau}}  \dtau \\
 && -\lambda  \frac{\sqrt{D_-}}{ \sqrt{\pi}}\int_0^t \frac{\alpha_-(s_1,\tau)}{\sqrt{t-\tau}}  \dtau.
\end{eqnarray*}

Following the method from Ref.~\cite{PirN}, we solve the system
corresponding to $\alpha_-(s_1,t)$ and $\alpha_+(s_1 ,t)$:
\begin{eqnarray*}
 &&D_-\alpha_-(s_1,t)+D_+\alpha_+(s_1,t)=2 D_+\frac{\del}{\del  s}K(0,s_1,D_+t),\\
&&D_-\alpha_-(s_1,t)+ \frac{\lambda\sqrt{D_-}}{ \sqrt{\pi}}\left(1+\sqrt{\frac{D_-}{D_+}} \right)  
\int_0^t \frac{\alpha_-(s_1,\tau)}{\sqrt{t-\tau}}  \dtau=4\lambda K(0,s_1,D_+t).
\end{eqnarray*}
We obtain therefore the Abel integral equation of the second kind for
$\alpha_-(s_1,t)$
\begin{eqnarray*}
&& \alpha_-(s_1,t)+\gamma \int_0^t\frac{\alpha_-(s_1,\tau) }{\sqrt{t-\tau}}\dtau=\frac{4\lambda}{D_-}K(0,s_1,D_+t),
\end{eqnarray*}
where $\gamma=\frac{\lambda}{\sqrt{\pi
D_-}}\left(1+\sqrt{\frac{D_-}{D_+}} \right)$.  Consequently,
\begin{eqnarray*}
 &&\alpha_-(s_1,t)=\frac{4\lambda}{D_-} K(0,s_1,D_+t)-\gamma\frac{4\lambda}{D_-}\int_0^t \frac{K(0,s_1,D_+\tau) }{\sqrt{t-\tau}}\dtau\\
&&+\pi \gamma^2\frac{4\lambda}{D_-}\int_0^te^{\pi \gamma^2(t-\tau)}\left(K(0,s_1,D_+\tau)-\gamma \int_0^\tau \frac{K(0,s_1,D_+s) }{\sqrt{\tau-s}}ds\right)\dtau.
\end{eqnarray*}
Using the Laplace transform yields, after simplifications:
\begin{eqnarray*}
G_{++}(s,s_1,t) &=& \frac{1}{\sqrt{4\pi D_+t}}\left(\exp\left(-\frac{(s-s_1)^2}{4D_+t}\right) + \exp\left(-\frac{(s+s_1)^2}{4D_+t}\right)\right) \\
&-& \frac{\lambda}{D_+} \exp\left(\frac{\lambda \alpha}{\sqrt{D_+}}(s+s_1)+ \lambda^2 \alpha^2 t\right)
\operatorname{Erfc}\left(\frac{s+s_1}{2\sqrt{D_+t}}+\lambda \alpha\sqrt{ t} \right),
\end{eqnarray*}
where $\alpha=\frac{1}{\sqrt{D_-}}+\frac{1}{\sqrt{D_+}}$.  By the same
way,
\begin{eqnarray*}
 G_{-+}(s,s_1,t) &&=\frac{\lambda}{\sqrt{D_-D_+}} \exp\left(\frac{\lambda \alpha}{\sqrt{D_-}}\left(-s+s_1\sqrt{\frac{D_-}{D_+}}\right)+ \lambda^2 \alpha^2 t\right) \\
&& \cdot \operatorname{Erfc}\left(\frac{-s+s_1\sqrt{\frac{D_-}{D_+}}}{2\sqrt{D_-t}}+\lambda \alpha\sqrt{ t} \right).
\end{eqnarray*}

We see that the Green function $G_{++}$ for $\lambda=0$ becomes the
Green function of the problem with the Neumann boundary conditions and
in this case $N(t)=0$, as $u_-\equiv 0$.  This property, $N(t)=0$, can
be also directly found using the Green function.
  
In $\R^n$ for $x=(s,\theta)$ and $y=(s_1,\theta_1)\in
\R\times\R^{n-1}$ we have
\begin{eqnarray*}
 &&G_{++}(s,\theta,s_1,\theta_1,t)_{\R^n}=G_{++}(s,s_1,t)_{\R}K(\theta,\theta_1,D_+t)_{\R^{n-1}},\\
&&G_{-+}(s,\theta,s_1,\theta_1,t)_{\R^n}=G_{-+}(s,s_1,t)_{\R}K(\theta,\theta_1,D_-t)_{\R^{n-1}}.
\end{eqnarray*}
Therefore in $\R^n$ for Varadhan's bounds with $x\ne y$ we have
\begin{eqnarray*}
 &&\lim\limits_{t\to 0+} t\ln G_{++}(x,y,t)_{\R^n} =-\frac{d(x,y)^2}{4D_+},\\
&&\lim\limits_{t\to 0+} t\ln G_{-+}(x,y,t)_{\R^n} =-\frac{d\bigl(s,s_1\sqrt{\frac{D_-}{D_+}}\bigr)^2+d(\theta,\theta_1)^2}{4D_-}.
\end{eqnarray*}

\begin{remark}\label{remarkGFl}
Applying this framework to the same system but with the transmittal
boundary condition for $0<\lambda<\infty$, we obtain
\begin{eqnarray*}
  &&G_{++}(s,\theta,s_1,\theta_1,t)=\frac{1}{(4\pi D_+t)^\frac{n}{2}}\left(\exp\left(-\frac{(s-s_1-tD_+k)^2}{4D_+t}\right) \right. \\
 && \left. + \exp\left(-\frac{(s+s_1-tD_+k)^2}{4D_+t}\right)\right)
\exp\left(-\frac{d(\theta,\theta_1)^2}{4D_+t}\right) \\
&&-\frac{1}{(4\pi D_+t)^\frac{n-1}{2}}\frac{\lambda}{D_+} \exp\left(\frac{\lambda \alpha}{\sqrt{D_+}}(s+s_1-tD_+k)+ \lambda^2 \alpha^2 t\right)\cdot\\
&&\cdot
\operatorname{Erfc}\left(\frac{s+s_1-tD_+k}{2\sqrt{D_+t}}+\lambda \alpha\sqrt{ t} \right) \exp\left(-\frac{d(\theta,\theta_1)^2}{4D_+t}\right),\\
&&G_{-+}(s,\theta,s_1,\theta_1,t)=\frac{1}{(4\pi D_-t)^\frac{n-1}{2}}\frac{\lambda}{\sqrt{D_-D_+}} \\
&& \cdot \exp\left(\frac{\lambda \alpha}{\sqrt{D_-}}\left(-\frac{D_+}{D_-}s+s_1\sqrt{\frac{D_+}{D_-}}+tD_+k\right)+ \lambda^2 \alpha^2 t\right) \\
&&\cdot\operatorname{Erfc}\left(\frac{-\frac{D_+}{D_-}s+s_1\sqrt{\frac{D_+}{D_-}}+tD_+k}{2\sqrt{D_-t}}
+\lambda \alpha\sqrt{ t} \right)\exp\left(-\frac{d(\theta,\theta_1)^2}{4D_-t}\right).
 \end{eqnarray*}
We also notice that for a fixed $t>0$ for $\lambda\to +\infty$ we obtain 
\begin{eqnarray*}
 &&G_{++}(s,s_1,t)\to \Gamma_{++}(s,s_1,t)\quad \hbox{and} \quad G_{-+}(s,s_1,t)\to \Gamma_{-+}(s,s_1,t).
\end{eqnarray*}

\end{remark}

\end{document}